\documentclass[3p,super,authoryear]{elsarticle}

\usepackage{amssymb,amsmath,mathtools}
\usepackage{amsthm}
\usepackage{xcolor}
\definecolor{darkgreen}{RGB}{0,80,0}  % simile a DarkGreen CSS
\usepackage{stmaryrd}
\usepackage{multirow}
\usepackage[hyphens]{url}
\usepackage[colorlinks=false]{hyperref}
\usepackage{booktabs}
\usepackage{graphbox}
\usepackage{longtable}
\usepackage{adjustbox}
\usepackage{geometry}
\usepackage[normalem]{ulem}
\usepackage{float}

\hypersetup{
  linkcolor=black,      % link interni (sezioni, ecc.)
  urlcolor=blue,        % link web
  citecolor=orange      %  colore citazioni
}

\newcommand{\revB}[1]{\textcolor{black}{#1}}
\newcommand{\revA}[1]{\textcolor{black}{#1}}
\newcommand{\sbal}[1]{\textcolor{black}{#1}}
\usepackage{geometry}
\newcommand*{\vcenteredhbox}[1]{\begingroup
\setbox0=\hbox{#1}\parbox{\wd0}{\box0}\endgroup}
\usepackage[scaled]{helvet}
 
\usepackage[T1]{fontenc}

\usepackage{titlesec}
\usepackage{tikz}
\newcommand{\whiteline}[1]{%
  \begin{tikzpicture}[baseline=(text.base)]
    \node[inner sep=0pt, outer sep=0pt] (text) {\strut #1};
    \draw[white, line width=2pt] (text.south west) -- (text.south east);
  \end{tikzpicture}%
}

\titleformat{\section}
  {\normalfont\large\bfseries}
  {\thesection}
  {1em}
  {\whiteline}

\titleformat{\subsection}
  {\normalfont\bfseries\itshape}
  {\thesubsection}
  {1em}
  {\whiteline}

\titleformat{\subsubsection}
  {\normalfont\itshape}
  {\thesubsubsection}
  {1em}
  {}

\usepackage[nomarkers,noheads]{endfloat}

\journal{Netherlands Journal of Geosciences}

\begin{document}

\begin{frontmatter}

%% Title, authors and addresses

%% use the tnoteref command within \title for footnotes;
%% use the tnotetext command for theassociated footnote;
%% use the fnref command within \author or \address for footnotes;
%% use the fntext command for theassociated footnote;
%% use the corref command within \author for corresponding author footnotes;
%% use the cortext command for theassociated footnote;
%% use the ead command for the email address,
%% and the form \ead[url] for the home page:
%% \title{Title\tnoteref{label1}}
%% \tnotetext[label1]{}
%% \author{Name\corref{cor1}\fnref{label2}}
%% \ead{email address}
%% \ead[url]{home page}
%% \fntext[label2]{}
%% \cortext[cor1]{}
%% \affiliation{organization={},
%%             addressline={},
%%             city={},
%%             postcode={},
%%             state={},
%%             country={}}
%% \fntext[label3]{}

\title{Unexpected fault activation due to underground gas storage in produced reservoirs. Part II: Definition of safe operational bandwidths}

%% use optional labels to link authors explicitly to addresses:
%% \author[label1,label2]{}
%% \affiliation[label1]{organization={},
%%             addressline={},
%%             city={},
%%             postcode={},
%%             state={},
%%             country={}}
%%
%% \affiliation[label2]{organization={},
%%             addressline={},
%%             city={},
%%             postcode={},
%%             state={},
%%             country={}}

\author[UNIPD]{Selena Baldan\corref{correspondingauthor}}
\ead{selena.baldan@phd.unipd.it}

\author[UNIPD]{Massimiliano Ferronato}
\ead{massimiliano.ferronato@unipd.it}

\author[UNIPD]{Andrea Franceschini}
\ead{andrea.franceschini@unipd.it}

\author[UNIPD]{Carlo Janna}
\ead{carlo.janna@unipd.it}

\author[UNIPD]{Claudia Zoccarato}
\ead{claudia.zoccarato@unipd.it}

\author[M3E]{Matteo Frigo}
\ead{m.frigo@m3eweb.it}

\author[M3E]{Giovanni Isotton}
\ead{g.isotton@m3eweb.it}

\author[UNIRO]{Cristiano Collettini}
\ead{cristiano.collettini@uniroma1.it}

\author[POLITO]{Chiara Deangeli}
\ead{chiara.deangeli@polito.it}

\author[POLITO]{Vera Rocca}
\ead{vera.rocca@polito.it}

\author[POLITO]{Francesca Verga}
\ead{francesca.verga@polito.it}

\author[UNIPD]{Pietro Teatini}
\ead{pietro.teatini@unipd.it}

\affiliation[UNIPD]{organization={Department of Civil, Environmental and Architectural Engineering, University of Padua},
                    city={Padua},
                    postalcode={35131},
                    country={Italy}}

\affiliation[M3E]{organization={M3E S.r.l.},
                  city={Padua},
                  postalcode={35121},
                  country={Italy}}

\affiliation[UNIRO]{organization={Department of Earth Sciences, Sapienza University of Rome},
                    city={Rome},
                    postalcode={00185},
                    country={Italy}}

\affiliation[POLITO]{organization={Department of Environment, Land and Infrastructure Engineering, Politecnico di Torino},
                  city={Turin},
                  postalcode={10129},
                  country={Italy}}

\cortext[correspondingauthor]{Corresponding author}

\begin{abstract}
Underground gas storage is a versatile tool for managing energy resources and addressing pressing environmental concerns. While natural gas is stored in geological formations since the early 20th century, hydrogen has recently been considered as a potential candidate toward a more flexible and sustainable energy infrastructure. Furthermore, these formations can additionally capture gases that contribute to climate change, such as CO$_2$. When such operations are implemented in faulted basins, however, safety concerns may arise due to the potential reactivation of pre-existing faults, which could trigger (micro)-seismicity events. In the Netherlands, it has been recently noted that fault reactivation can occur “unexpectedly” during the life of an underground gas storage (UGS) site, even when stress conditions are not expected to cause a failure. The present two-part work aims to develop a modeling framework to investigate the physical mechanisms causing such occurrences in previously produced gas reservoirs and define a safe operational bandwidth for pore pressure variation for UGS operations in the faulted reservoirs of the \revB{Upper Rotliegend Group}, the Netherlands. This follow-up paper investigates in detail the mechanisms and crucial factors that result in fault reactivation at various stages of a UGS. The mathematical and numerical model described in Part I is used, also considering how the presence of stored gases may influence the mechanical properties of the reservoir and caprock, in particular the Young modulus. The study investigates the hazard of fault activation caused by the storage of different fluids for various purposes, such as long-term CO$_2$ sequestration, CH4 and H$_2$ injection and extraction cycles, and N$_2$ injection as cushion gas. The results show how geological configuration, geomechanical properties, and reservoir operating conditions may increase the hazard of fault reactivation at various UGS stages. Furthermore, the analysis indicates that reservoir pressure near critical levels and reactivation during primary production may significantly contribute to potential fault instabilities. Operational guidelines for improving secure and effective storage operations are thereby presented.
\end{abstract}

%%%Graphical abstract
%\begin{graphicalabstract}
%%\includegraphics{grabs}
%\end{graphicalabstract}

%%Research highlights
%\begin{highlights}
%\item Research highlight 1
%\item Research highlight 2
%\end{highlights}

\begin{keyword}
Carbon Capture and Sequestration \sep Critical pressure \sep
Fault reactivation \sep Safety guidelines
%% keywords here, in the form: keyword \sep keyword

%% PACS codes here, in the form: \PACS code \sep code

%% MSC codes here, in the form: \MSC code \sep code
%% or \MSC[2008] code \sep code (2000 is the default)

\end{keyword}

\end{frontmatter}

%% \linenumbers

%% main text
\section{Introduction}
The development of underground gas storage sites \revA{(UGS)} has been crucial in the global energy infrastructure since the early 20th century. The first operational underground storage system (USS) site dates back to 1915 \revA{in the Welland gas field (Canada)} \citep{alshafi_2023_a} . Over time, this technology has expanded with more than 600 facilities worldwide \citep{Fou_etal18} to effectively manage gas supply and demand on both seasonal and daily bases \citep{verga_2018}. While USS has typically focused on storing natural gas (i.e., mainly CH$_4$) through UGS, it possesses the potential to store a variety of gases for diverse purposes, thus becoming a flexible tool for managing energy resources and addressing environmental concerns. This is especially important as Europe moves toward a net-zero greenhouse gas emission energy system \citep{netzero2050}. The scientific community has shown a growing interest in underground hydrogen storage (UHS) as an alternative to natural gas. Having the potential to accommodate significant volumes, ranging from tens of millions of cubic meters in caverns to potentially billions of cubic meters in depleted gas fields \citep{groenenberg_2020_largescale}, hydrogen versatility and absence of carbon emissions during production and utilization make it a sustainable option. Additionally, carbon capture and storage (CCS) initiatives have identified underground storage of CO$_2$ as a crucial component in mitigating greenhouse gas emissions and contrasting climate change.

However, despite the widespread use of USS, one of the important aspects to be considered is %inherent risks are associated with
the potential reactivation of existing faults and consequent induced seismicity \citep{ellsworth_2013_injectioninduced,Fou_etal18,KerWei18}. Although such events are statistically rare \citep{Fou_etal18}, in terms of the number of affected sites rather than the frequency of events within a single reservoir, they demand significant attention due to their social and economic implications \revA{\citep{van2015social}}.
In the Netherlands, natural gas is currently stored in four UGS facilities: Bergermeer, Alkmaar, Norg, and Grijpskerk. They are located in the \revB{Upper Rotliegend group} that widely extends in Central Europe and is considered one of the most extensively explored petroleum systems worldwide \citep{Gau03}. Over the past few decades, seismic activity associated to USS has been observed in this region \citep{Wee_etal14,Uta17}, particularly in three of the four UGS facilities (i.e., Bergermeer, Norg, and Grijpskerk). 
\revA{Aside from the Castor in Spain \citep{Ces_etal14,Vil_etal21} and Hutubi in China \citep{Jia_etal21,Liu_etal23}}, these Dutch reservoir represent the only cases where UGS activities have been definitively associated to induced seismicity events~\citep{Vil_etal21}.

Recently, the Netherlands have also focused attention on the feasibility of storing different gases in depleted gas fields and salt caverns as a crucial measure in the ongoing energy transition. The same technology used in UGS is being applied to offshore CO$_2$ storage through the Porthos project \citep{PorthosCO2}, while research on hydrogen storage is in progress and primarily focused on salt caverns, although the potential for storage in abandoned gas fields is also being considered for the future. Additionally, besides being used as cushion gas \citep{shoushtari2023utilization}, nitrogen is utilized to convert high-calorific gas to low-calorific gas for residential heating and cooking purposes in Dutch households. Currently, nitrogen is stored in a salt cavern near Heiligerlee for this purpose \citep{nlog,mun_etal22}.

Given the future role that USS will play, it is crucial to understand why faults can become active and which mechanisms are potentially prone to \revA{cause} instability. During the cyclic or permanent storage of these fluids, activation can occur when the natural stress regime on a fault surface is altered due to changes in pore pressure within a reservoir, \revB{although additional mechanisms may also contribute}. \revB{The combined effects of human-induced pressure and stress changes at depth, the pre-existing stress field, and the frictional and mechanical properties of the rocks and faults determine the initiation, amount of slip, and extent of reactivation} \citep{Seg_etal94, Het_etal00, Can_etal19}.

\begin{figure}
    \centering
    \includegraphics[width=0.9\linewidth]{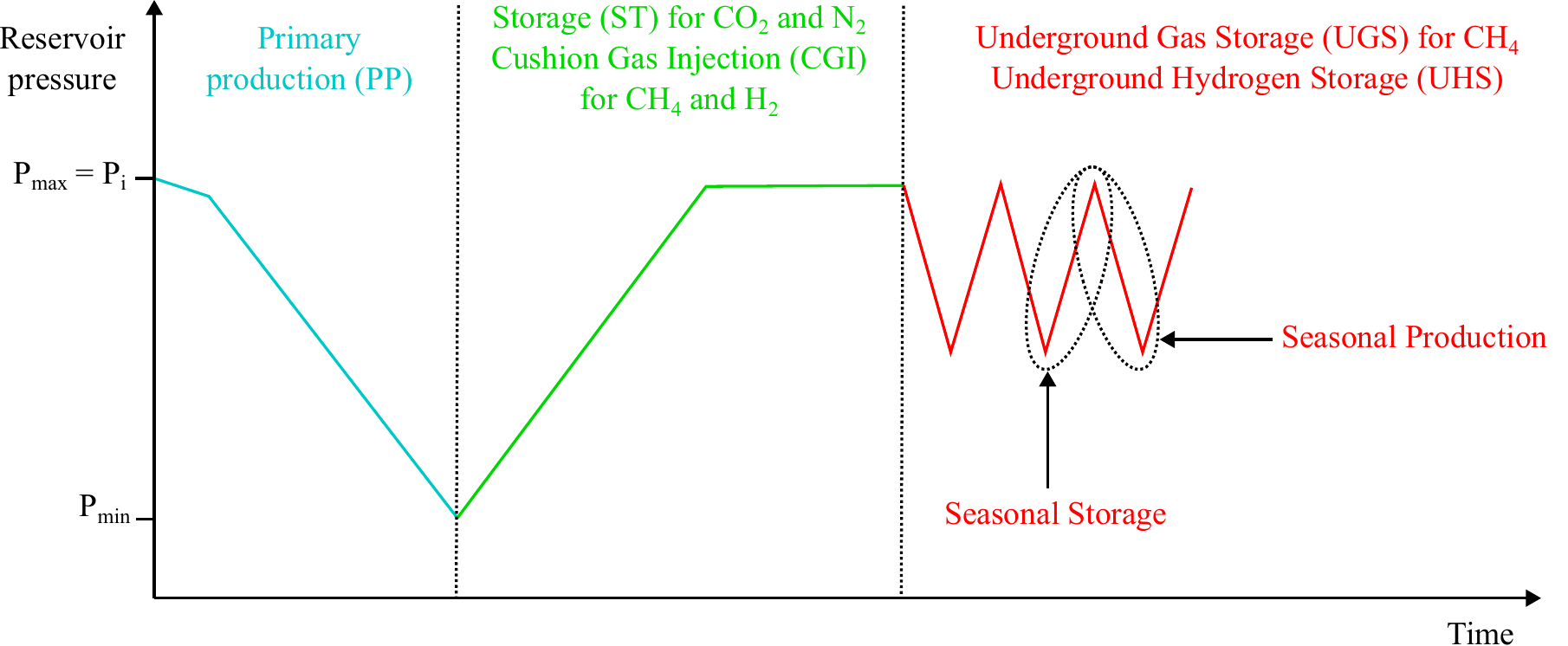}
    \caption{Sketch of the main stages during the reservoir lifespan in terms of pressure change over time for CH$_4$, CO$_2$, H$_2$ and N$_2$: Primary Production (PP) for CH$_4$ production, Storage (ST) for CO$_2$ and N$_2$ storage, Cushion Gas Injection (CGI) for H$_2$ and CH$_4$ injection, and Underground H$_2$ Storage (UHS) and Underground Gas Storage (UGS) cycles for seasonal pressure fluctuations during CH$_4$/H$_2$ injection and withdrawal.
    Notice that, according to the Netherlands regulation, $P$ must remain below $P_i$ regardless of the stage under consideration \citep{minEZK_grijpskerk2022,tno_ondergrondse2018}}.
    \label{fig:sketch_p}
\end{figure}

Over a short- to mid-term timeframe (e.g., from days to years), fault reactivation can be interpreted through a geomechanical approach involving pressure variations following injection/extraction activities \revA{that can last for decades}. 
In this context, most human-induced seismic events can be easily explained, \revA{as they are triggered by fluid injection that increases pore pressure beyond its initial reservoir pressure, thereby driving shear stress on the fault surface to its failure limit \citep{walsh_2016, SodM1}}. However, a subset of recorded events cannot be accounted for by this mechanism. These events, referred to as ``unexpected'' seismic events, \revB{typically exhibit small magnitudes} and occur when the pressure falls within the pressure range already experienced during primary production, \revA{that is, pressures lower than the initial value $P_i$ but still above the minimum pressure generally experienced by the reservoir at the end of depletion $P_{\min}$. In this scenario, fault reactivation occurs not only during primary production or gas storage at pressures exceeding the initial value $P_i$ \citep{Def_etal95, Ces_etal14, Zho_etal19}, but also during injection or producing and storing phases in which pore pressures remain within the previously experienced pressure interval \citep{MIT09, Kra_etal13, TNO15, NAM16}.} A typical time-behavior of the fluid pressure variations expected in \revA{USS} applications is sketched in Figure~\ref{fig:sketch_p}.
The recorded seismic events \revA{in UGS} are unexpected because they occur at a stress state already experienced by the reservoir and surrounding faults during previous production phases. \revB{However, the mechanical response of the reservoir–fault system may evolve over time, meaning that previously experienced stress conditions do not necessarily lead to identical system behavior \citep[e.g.][]{pijnenburg2019inelastic}}. Unraveling the underlying mechanisms behind these ``unexpected'' seismic events is crucial for ensuring the safe and efficient operation of underground gas storage facilities.

On the other hand, long-term processes, spanning hundreds or even thousands of years, primarily involve fluid-rock interactions that can alter fault rock fabric, mineralogy, frictional properties, and cohesion. Currently, the understanding of the specific impacts of different gases on the mechanical properties of faults and rock is still poorly understood, but the potential risk associated with their storage has to be considered. Indeed, the performance, efficiency, and safety of any storage greatly depends on fluid-rock interaction.

In this follow-up paper of %Franceschini et al. (2024)
\cite{SodM1}, we conduct an investigation of the phenomenon of ``unexpected'' seismic events in underground gas storage facilities. By understanding the factors contributing to such events, we can establish enhanced guidelines and best practices for operating underground gas facilities during different storage phases. We present a comprehensive study to reassess the fundamental geomechanical causes of induced seismicity and specific risk factors, including the potential impact of gas type on reservoir and caprock mechanical weakening, particularly with regard to changes in elastic properties.

This work has four specific aims. The first is to improve the understanding of the physics-based processes responsible for induced seismicity during both cyclic and permanent gas storage. \revB{To this end, the analysis does not attempt to reproduce individual field events, but rather to identify the mechanical conditions and stress-path evolutions under which fault reactivation may occur during storage operations.} \revA{The second aim is to investigate how the peculiarities of the Rotliegend reservoirs, such as graben-bounding faults and intra-field faults responsible for reservoir compartmentalization, affect fault reactivation.} 
\revB{The third aim is to assess how different injected or stored fluids (CH$_4$, CO$_2$, H$_2$, and N$_2$), through their associated effects on rock elastic properties such as the Young’s modulus, impact the physical mechanisms of fault reactivation}.  \revB{The final aim is to provide, based on the modeling outcomes, recommendations for operational guidelines defining safe pressure bandwidths for reservoir storage conditions.} These guidelines are relevant to pressure cycling (and noncycling) during the underground storage of CH$_4$, CO$_2$, H$_2$, and N$_2$ in depleted gas fields with configurations similar to those of the reservoirs located in the Rotliegend Group, and are intended to assist risk assessment and mitigation strategies in underground gas storage operations. A few preliminary outcomes regarding UGS are already reported in %Teatini et al. (2019,2020)
\cite{Tea_etal19,Tea_etal20} and %Franceschini et al. (2024)
\cite{SodM1}.

It is important to emphasize that this study specifically addresses the hazard of \revA{fault reactivation in UGS sites that are hosted by depleted hydrocarbon reservoirs, where the mechanical behavior of the formation is already well constrained}. Gas storage in other subsurface systems, such as saline aquifers, salt caverns and unexploited reservoirs, may differ substantially in terms of knowledge of geomechanical setting, and are not addressed by this work. As such, the findings and modeling framework presented here are not intended to be generalized to those settings. 

The analysis that follows is accomplished by using the mathematical and numerical framework developed in Part I of this work \citep{SodM1}, where a Finite Element (FE) geomechanical simulator is implemented with a Lagrange multiplier-based treatment of the frictional contact conditions for reproducing the fault behavior.
The simulator is supplemented with
%through an advanced finite elements - interface elements (FE-IE)
visco-elasto-plastic constitutive laws~\cite{Iso_etal19}, and fault activation is ruled by the Coulomb frictional criterion. Pressure change within the faults, variation of Coulomb's parameters due to slip-weakening, and the rheology of the caprock are properly accounted for. Because of the main aims listed above, the modeling study is not focused on a specific UGS field, but on a conceptual reservoir and on a fault system that realistically represents the main geologic features of the Upper Rotliegend group.

The paper is organized as follows. The typical setting of the Rotliegend reservoirs and a few examples of seismicity recorded during UGS are initially presented. Then, fluids under investigation are described, and the modeling set-up is outlined. The various scenarios are then developed to understand the geomechanical behavior of the faulted system, first analyzing the parameters involved in UGS and then comparing the scenarios using different fluids. Modeling results are presented, and the mechanisms responsible for fault reactivation during the different phases are identified. The modeling outcomes are discussed, highlighting the peculiarities of the storage scenarios with respect to induced seismicity during primary production as addressed in previous studies, ranking the natural features and anthropogenic factors prone to cause fault reactivation, and illustrating general operational guidelines to reduce the probability of reactivation occurrences.

\subsection{Induced seismicity in Rotliegend UGS reservoirs}
\begin{figure}
    \centering
    \includegraphics[width=1.0\linewidth]{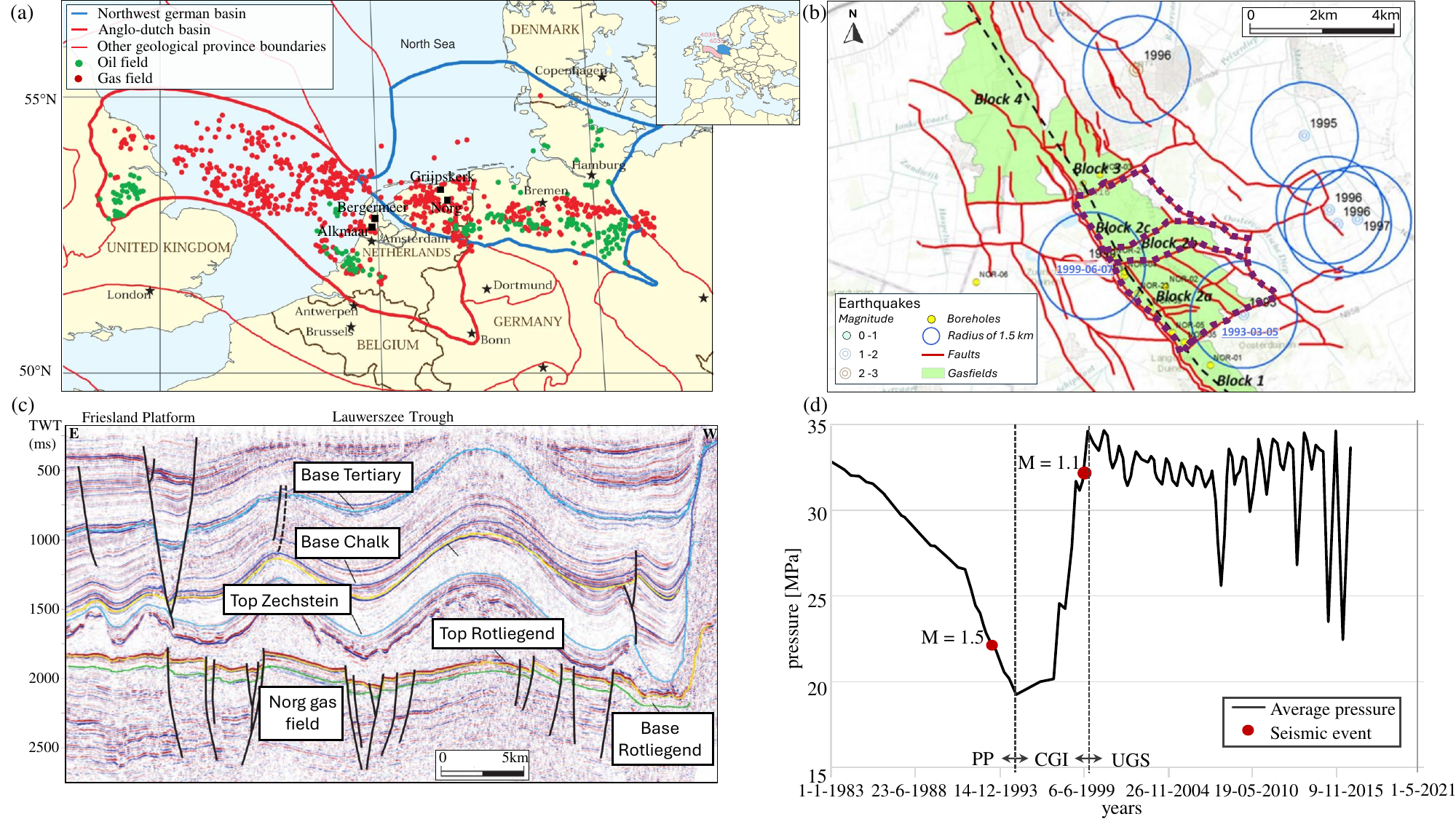}
    \caption{Key geological and seismological features of the Rotliegend formation. a): Tectonic provinces of the Carboniferous–Rotliegend Total Petroleum System with the hydrocarbon field centerpoints distinguished between oil and gas fields (modified after %Gautier (2003)
    \cite{Gau03}).  b): Top view of the Norg UGS site (in green), with bounding fault traces (red lines) and recorded seismic events (blue circles). The dashed purple line surrounding Blocks 2a, 2b, and 2c illustrates an example of reservoir compartmentalization (modified after %TNO Technical Report (2015)
    \cite{TNO15}). c): Interpreted seismic profile across the northeast Netherlands, showing the Norg gas field located in fault-bounded blocks at Rotliegend level (modified after %de Jager and Geluk (2007)
    \cite{JAG07}). d)t: Average pressure evolution at the Norg UGS reservoir over time, with annotated seismic events from 1993 and 1999 (modified after %Nederlandse Aardolie Maatschappij BV (2016)
    \cite{NAM16}).}
\label{fig:geol-seismic}
\end{figure}

The Rotliegend Group extends over a wide region that includes North-Western Germany and the Netherlands, with an offshore portion into the southern North Sea and parts of the UK (Figure~\ref{fig:geol-seismic}~a). The \revB{Rotliegend} Group sandstones \revA{were deposited} in a variety of arid, terrestrial environments, among which ephemeral fluvial (wadi) systems, various types of eolian deposits, desert-lake environments, and adjacent sabkhas \revA{were} dominant. The sealing nature of the \revB{Zechstein Group salt} is the main reason for the existence of the large number of traps that retained gas migrated from the underlying coals and carbonaceous shales of the \revB{Carboniferous Limburg} Group for approximately 150 million years \citep{Gau03}.
\revB{In the Dutch sector, the reservoirs of the Upper Rotliegend Group} spans a depth range between 2,000 and 4,700~m, and is generally \revA{composed of permeable sandstones} (average 100-200 mD), with high net-to-gross and porosity (average 15-20\%). %, considering its burial depth should have induced significant diagenetic degradation.
 The thickness of the reservoir sands is generally over 50~m and rarely exceeds 300~m. Faulting is a characteristic feature of this geological unit. \revB{Most present-day faults represent reactivations of inherited Mid-Palaeozoic basement structures that were deformed during the Variscan orogeny and subsequently reactivated during Permian, Mesozoic, and Cenozoic tectonic phases \citep{geluk2005stratigraphy}. During the Mesozoic, successive rifting episodes further reactivated and enhanced the extensional normal-fault pattern that characterizes the region \citep{ziegler1990,deJ07}}. As a consequence, many fields are formed by a number of rhomboid-shaped dipping fault blocks, positioned one next to each other. The trapping mechanism \revA{of gas} is structural, by the juxtaposition of the Rotliegend reservoir blocks against the thick Zechstein formation. Moreover, in several Rotliegend fields, the hydraulic connectivity between fault-bounded blocks is limited, causing a clear compartmentalization of the reservoirs. This behaviour results from a combination of fault sealing, \revB{fault offsets that produce structural barriers through lithological juxtaposition}, and stratigraphic and diagenetic heterogeneities \citep{Hul10}. \revB{Nevertheless, some Dutch Rotliegend reservoirs (e.g., Groningen) do contain transmissive faults~\citep{de2017}. In such settings, compartmentalization is mainly structural rather than hydraulic.}
 Figure~\ref{fig:geol-seismic}~c sketches a representative example of the structural map of the Norg UGS reservoir.

%>\begin{figure}
%    \centering
%    \includegraphics[height=0.33\linewidth]{figs/Model/Norg_press}
%    \caption{Left: Map of the Norg UGS site (in blue) with traces of the bounding faults (black lines) and localization of the recorded seismic events (red cirles). Right: average pressure evolution over time at the Norg UGS field. The location of the two seismic events in 1993 and 1999 are reported (modified after \cite{NAM16}).}
%    \label{fig:geol-seismic}
%\end{figure}

%\begin{figure}
%    \centering
%    %\includegraphics[width=0.8\linewidth]{figs/Model/Norg}
%    \null\hfill
%    \includegraphics[height=0.28\linewidth]{figs/Model/Norg_Zechstein}
%    \hfill
%    \includegraphics[height=0.28\linewidth]{figs/Model/Grijpskerk_press}
%    \hfill\null
%    \caption{Left: Map of the Norg UGS site (in blue) with traces of the bounding faults (black lines) and localization of the recorded seismic events (red cirles). Right: average pressure evolution over time at the Norg UGS field. The location of the two seismic events in 1993 and 1999 are reported (modified after \cite{NAM16}).}
%    \label{fig:geol-seismic}
%\end{figure}

The reservoir is bounded by normal faults with a \revA{significant throw (up to a 250~m)} and consist of a few compartments that are separated by internal faults. The gas fields are 2,500-3,500~m deep, with the Rotliegend reservoir rock characterized by an average net thickness of 150-200~m (Figure~\ref{fig:geol-seismic}~c. 
As usual for UGS plants, the reservoir is a partially-depleted gas field that was converted \revA{into an UGS} after a %10-year long conventional
primary production (PP) period. PP is generally followed by a quick gas injection (CGI) phase when fluid pressure in the reservoir is risen to almost the original value, and then by the UGS injection and withdrawal cycles. Gas is injected at temperatures similar to, or slightly below, those at reservoir depth, thus minimizing thermal stresses. Figure~\ref{fig:geol-seismic}~d summarizes the evolution of the average pressure in Norg. The timing and location of the seismic events recorded in the vicinity of the UGS reservoir is also shown in Figure~\ref{fig:geol-seismic}~b,d. In Norg, the time of the first earthquake \revA{($M=1.5$)}  coincided with low pressure level during PP. The second seismic event \revA{$M=1.1$} occurred at the end of the CGI phase, with the pressure in the reservoir close to the original pressure level. In Grijpskerk, where pressure differences between reservoir blocks can be as high as about 7~MPa, a first seismic event occurred at the end of PP in 1997 \revB{at a pressure about 11~MPa below the initial value}, and a second one in 2015, near the end of a UGS production phase, \revB{at about 12~MPa below the initial conditions \citep{TNO15}}. Four seismic events with magnitudes $M \in [3.0 , 3.5]$ were detected in 1994 and 2001 during PP in the Bergermeer reservoir, \revB{when reservoir pressures had decreased by about 18–20~MPa from the initial pressure conditions}. The events were located at the tip of the central fault separating the two main blocks \citep{TNO15}. A down-hole microseismic array characterized by a magnitude of completeness below $M=0$ was established during the CGI phase. The array recorded a large number of microseismic events ($M<0$) on the faults at the reservoir depth and a main event with $M \approx 0.8$ along the central fault when the pressure difference between the two blocks peaked at about 4~MPa in early 2013 \citep{Bai_etal16}.
\revB{Across these Dutch Rotliegend UGS sites, injection-related events are of very small magnitude (typically $M < 1.5$), and thus contribute only a negligible fraction of the cumulative seismic moment compared to the $M$ $3 - 3.5$ events produced during primary depletion.}
Table~\ref{tab:seism} presents an overview of induced seismicity cases across the different reservoir sites, detailing the operational phase, activated fault geometry, and event magnitudes. Table~\ref{tab:regime} follows with a summary of the natural stress regimes, fault orientations, and key geological and operational characteristics of the reservoirs.
\revB{Taken together, these cases show a common pattern: fault reactivation can occur not only at high pore pressures but also within pressure ranges that were previously experienced without seismicity. This indicates that the key control on stability is the evolving stress path of the reservoir–fault system, influenced by differential compaction, fault-block displacement, and changes in fault-zone hydraulic properties. These shared features form the basis of our geomechanical modelling, which aims to investigate the physical mechanisms capable of triggering such “unexpected’’ seismicity during UGS operations.}

\begin{table}[]
\centering
\footnotesize
\begin{tabular}{|p{3.5cm}|p{1.2cm}|p{4.0cm}|p{2.5cm}|p{0.7cm}|p{0.7cm}|p{0.7cm}|p{0.5cm}|}
\hline
\multirow{2}{*}{Reservoir} & \multirow{2}{*}{Type / gas}  & \multirow{2}{*}{Activated fault geometry} & \multirow{2}{*}{\revB{$\Delta P$ at event (MPa)}} & \multicolumn{3}{c|}{Magnitude} \\
\cline{5-7}
                           &                             &                        &                                           & $\leq$1 & 1--3 & $\geq$3 \\ 
\hline

Norg \cite{NAM16,TNO15} $\quad\quad$ The Netherlands 
  & UGS (CH$_4$) 
  & Slip motion on a sub-vertical NW--SE fault (PP) and a sub-vertical compressive SW--NE fault (UGS) 
  & \begin{tabular}[c]{@{}l@{}}\\PP (1993) $\sim-9$\\ UGS (1999) $\sim-2$ \end{tabular} 
  &  & \begin{tabular}[c]{@{}l@{}}x\\ (2)\end{tabular} & \\ 
\hline

Grijpskerk \cite{TNO15} $\quad\quad$ The Netherlands 
  & UGS (CH$_4$) 
  & Large uncertainty in locating the active fault(s) 
  & \begin{tabular}[c]{@{}l@{}}\\UGS (1997) $\sim-11$ \\ UGS (2015) $\sim-12$ \end{tabular} 
  &  & \begin{tabular}[c]{@{}l@{}}x\\ (2)\end{tabular} & \\ 
\hline

Bergermeer \cite{TNO15,Orl_etal13} $\quad\quad$ The Netherlands 
  & UGS (CH$_4$) 
  & Central fault with strike SSW and dip 60--65$^\circ$ 
  & \begin{tabular}[c]{@{}l@{}}\\PP (1994): $\sim -18$ \\ PP (2001): $\sim -20$\\ CG(2013): $\sim -15$\end{tabular}
  & \begin{tabular}[c]{@{}l@{}}x\\ (1--CG)\end{tabular} 
  &  & \begin{tabular}[c]{@{}l@{}}x\\ (1--PP)\end{tabular} \\ 
\hline

\end{tabular}
\caption{Inventory of induced seismicity cases in the Netherlands, categorized by activated fault geometry, pressure conditions, and event magnitude.
\revB{PP = primary production, CG = cushion gas injection, UGS = underground gas storage.
$\Delta P$ denotes the pressure difference at the time of the event with respect to the initial reservoir pressure $P_i$.
The notation ``x (M--Phase)'' indicates a recorded seismic event of approximate magnitude $M$ during the specified operational phase.}}
\label{tab:seism}
\end{table}

\begin{table}[]
\centering
\footnotesize
\begin{tabular}{|p{2.5cm}|p{3.5cm}|p{6.0cm}|p{3.0cm}|}
\hline
Reservoir   & Stress regime & Reservoir info & Fault orientation \\ \hline
Norg \cite{NAM16,TNO15}  & Normal faulting; $\sigma_v > \sigma_H > \sigma_h$; SHmax $\approx$ 1.92 bar/10m (171° E of N), $\sigma_h$ $\approx$ 1.44 bar/10 m, and $\sigma_v$ $\approx$ 2.2 bar/10 m & Depth: 2670 m; Thickness: 140 m. Primary production started in 1983 until 1995 with a max DP of about 117 bar; UGS started in 1995–1997 with a DP for cycle of about 57 bar & Steeply dipping NW–SE and SW–NE; mostly non-sealing \\ \hline
Grijpskerk \cite{TNO15} & Normal faulting; $\sigma_v > \sigma_H > \sigma_h$; SHmax oriented $\approx$ N170°E & Depth: 3300 m; Thickness: 220 m. Primary production started in 1993 until 1994 with a max DP of about 53 bar; UGS started in 1997 with a DP for cycle of about 127 bar & Steeply dipping NW–SE and NE–SW; mostly non-sealing \\ \hline
Bergermeer \cite{TNO15,Orl_etal13} & $\sigma_v$ = 25 MPa; $\sigma_{h\,min}$ = 9 MPa (effective stress); initial direction of the minimum horizontal stress is oriented NE–SW & Depth: 2200 m; Thickness: 200 m. Primary production started in 1971 until 2005 with a max DP of about 213 bar; no activity between 2007 and 2010; cushion gas injection between 2010 and 2012; UGS started in 2013 & Not specified here; typically NW–SE in region \\ \hline
\end{tabular}
\caption{Overview of induced seismicity cases in relation to the stress regime, seismic monitoring, reservoir characteristics, and production history.}
\label{tab:regime}
\end{table}

\subsection{Chemo-mechanical effects by CO$_2$, H$_2$ and N$_2$ on rocks and faults}\label{sec:fluids}
 \subsubsection{\revA{Effect on CO$_2$}}
The interactions between gas, brine, and rock can lead to significant changes in the reservoir and fault properties that, in turn, can greatly impact the performance, efficiency, and safety of the storage process. However, the impact of injected gas type on the mechanical response of sedimentary rocks is still not well understood. Existing research predominantly focuses on CO$_2$ injection \citep{kim_2022_short}, with limited investigation into N$_2$ injection \citep{hu_2016_a,fuchs_2019_geochemical}.
There is some ongoing development regarding H$_2$ injection into sedimentary rocks \citep{borello_2024_underground, vasile_2024_innovative}. The latest comprehensive review on H$_2$ injection and withdrawal has been compiled by %Miocic et al. (2023)
\cite{Mic_etal23}.  %\af{Non so se ti riferisci a questo punto quando parli di ref mancanti, ma qui ne servirebbero un paio, una per gas almeno.}
From the extensive available literature on CO$_2$-fluid-rock interactions, %particularly studies by
e.g., %Rohmer et al. (2016)~
\cite{rohmer_2016_mechanochemical} and %Peter et al. (2022)
\cite{peter_2022_a}, %and Huang et al. (2020)(?),
two contrasting viewpoints emerge:  CO$_2$ can either have a negligible influence on the rock mechanical parameters, or degrade them because of mineral dissolution and modifications in pore size distribution. %induced by its presence.
%It has been observed that in certain types of geological formations, CO$_2$ injection and its interactions with the reservoir rock do not significantly impact the mechanical properties. This stability
A stable
condition is often observed when the mineral composition and pore structure remains unchanged during CO$_2$ injection \citep{rimmel_2010_evolution, bolourinejad_2015_chemical}, with the main elastic properties %such as Young modulus \citep{mikhaltsevitch_2014_measurements}, Poisson ratio, and mechanical strength remain
largely unaffected by the CO$_2$ presence\citep{mikhaltsevitch_2014_measurements,hu_2016_a}.
Conversely, in cases where a weakening effect is observed after CO$_2$ injection, the mechanical strength of the rock decreases, leading to increased permeability and reduced rock stiffness \citep{hol2018,manj2023}. The dissolution of carbonate minerals, particularly calcite, reduces the cohesion of the rock matrix and weakens grain-to-grain contacts \citep{kim_2022_short}. This effect has been observed in experiments on sedimentary rocks like sandstones and carbonates. %Bolourinejad et al., 2014(?).
Generally, carbonate-rich rocks are more susceptible to mechanical weakening due to CO$_2$-fluid-rock interactions.

While several experimental evidences support these findings, uncertainties still persist in the results. For instance, CO$_2$-acidified brine injection can precipitate minerals, potentially increasing the strength and stiffness of the rock depending on the specific material and its mineral composition. CO$_2$ can react with certain minerals like calcite, inducing mineral precipitation and cementation, thereby enhancing the mechanical integrity of the rock \citep{dnicolasespinoza_2018_co2}. Sandstones, for example, may experience both increase and decrease in the deformation moduli \citep{HANGX2013, shi2019, tar2020}.
The variability in results regarding the effects on rock stiffness and strength highlights the need for additional experimental data to address existing uncertainties.
%% Giustificazione 30%
%\sout{Regarding the geochemical effects of CO$_2$ injection in the Upper Rotliegend formation, \cite{bolourinejad_2015_chemical} reported an increase in permeability ranging from 10\% to 30\%. This suggests changes in porosity and deformation moduli, with potential variation of dynamic deformation modulus up to $\pm$30\% \citep{harbert_2020}.} 

Regarding the geochemical effects of CO$_2$ injection in the Upper Rotliegend formation, literature \citep{bolourinejad2015effects,bolourinejad_2015_chemical} reported an increase in permeability ranging from 10\% to 30\%. The increase in permeability induced by CO$_2$ injection is often accompanied by microstructural changes, such as a change in porosity, which leads to a variation in the Young modulus, reaching up to ±30\% in some cases \citep{harbert_2020}.
% Based on experimental studies on sandstones \citep{HANGX2013, shi2019, tar2020}, a range of stiffness changes has been observed, with variations in Young modulus reaching up to ±30\% in some cases \citep{harbert_2020}.}

\subsubsection{\revA{Effect on H$_2$ and N$_2$}}
Regarding N$_2$ and H$_2$ injection,  %into sedimentary rocks,
there is a limited investigation. N$_2$ is used mainly to determine the porosity and permeability of the rock or as a control gas \citep{fuchs_2019_geochemical}. Therefore, the exposure of rocks to N$_2$ is expected to have negligible effects, although further investigations are required. %(Hu et al., 2019?).
In the case of H$_2$ injection, rock degradation is expected. Hydrogen injection in UHS porous reservoirs may trigger the dissolution of carbonates and sulfates, as well as grain crushing and local compaction, possibly altering the reservoir porosity and mechanical and/or flow properties \citep{alyaseri_2023_experimental}. The response will be time-dependent, with pH equilibration to higher values potentially leading to further compaction and stiffening of the reservoir. However, this effect is likely to stabilize after a few cycles. Consequently, the poroelastic response of the reservoir will heavily depend on the chemical environment in the early stages of operations and become more stable and predictable in later stages. %[KEM39-WP2].
Since no data are available for Young modulus variations in the Rotliegend reservoirs, we used the reported variability for CO$_2$ systems as a reference and assumed a ±30\% change in the modeling simulations.

% \subsection{Effect on faults}
From a geochemical point of view, gas injection may also affect faults, which represent the key factors responsible for fluid migration and containment. Gas injection can affect fault permeability, reactivation potential, and overall stability, which could lead to unintended fluid pathways or seismic activity.
Concerning CO$_2$, during the storage phase over long periods (of the order of thousands of years), fluid-rock interactions can lead to changes in fault rock fabric and mineralogy. This may result in fault frictional properties and cohesion alterations due to dissolution and cementation processes \citep{collettini_2008_fault}.
However, laboratory friction tests comparing fault rocks exposed to CO$_2$ for several hours to samples taken from natural analogs exposed to CO$_2$ over millions of years, show limited friction reduction due to CO$_2$ exposure \citep{samuelson_2012_fault}. Moreover, according to the author's knowledge, no specific laboratory studies in the Rotliegend formation show how CO$_2$-enhanced dissolution and re-precipitation processes may impact on fault cohesion. Based on the literature, the chemical effects on fault mechanical properties after CO$_2$ injection are considered as negligible.

Finally, there is a lack of documentation on the effect of H$_2$ and N$_2$ injection on faults. %frictional properties and cohesion.
Experience with underground hydrogen storage in porous geological formations is minimal, with practical applications restricted to the storage of town gas, i.e., a mixture of gases containing 25-60\% hydrogen, along with smaller amounts of CH$_4$, CO, and CO$_2$ \citep{heinemann_2021_enabling}. Injecting hydrogen into a porous reservoir can change the chemical equilibrium of the formation pore water, leading to fluid-assisted grain-scale processes such as cement dissolution, clay mineral sorption/desorption within grain boundaries, stress corrosion cracking, dissolution-precipitation, and/or intergranular frictional slip \citep{heinemann_2021_enabling}. Although these grain-scale mechanisms are well-studied, little is known about the specific effects of hydrogen and nitrogen on their rates.
\section{Methods and materials}\label{sec:model}

We use the mathematical and numerical framework developed, discussed, and tested in Part~I of this work \citep{SodM1}.
The modeling approach consists of one-way coupled Finite Element (FE) hydro-poromechanical simulations in fractured geological media. 
\revA{For ease of reference, the main modelling features relevant to the present work are briefly summarized here. The porous medium is treated as a quasi‐static, linear elastic, saturated continuum under small strains, whereas faults are represented as zero–thickness internal surfaces governed by Coulomb friction with optional slip–weakening. Normal and tangential contact tractions are enforced by Lagrange multipliers, leading to a mixed finite element formulation in which displacement is approximated with first–order hexahedral elements and fault tractions with piecewise constant interface elements. The non–linear contact problem (stick–slip–open) is solved by an active–set strategy combined with an exact Newton method and a preconditioned Krylov solver. Pore pressure fields are computed beforehand with a multiphase flow simulator and then imposed in a one–way coupled fashion, acting on the solid skeleton through Terzaghi–Biot effective stress. The same hexahedral grid is shared by the flow and geomechanical models, ensuring conservative transfer of pore pressures to the fault elements.} For all the details, the reader is referred to Part I \cite{SodM1}.

\begin{figure}%[htbp]
    \centering
    \null\hfill
    \includegraphics[width=1\textwidth]{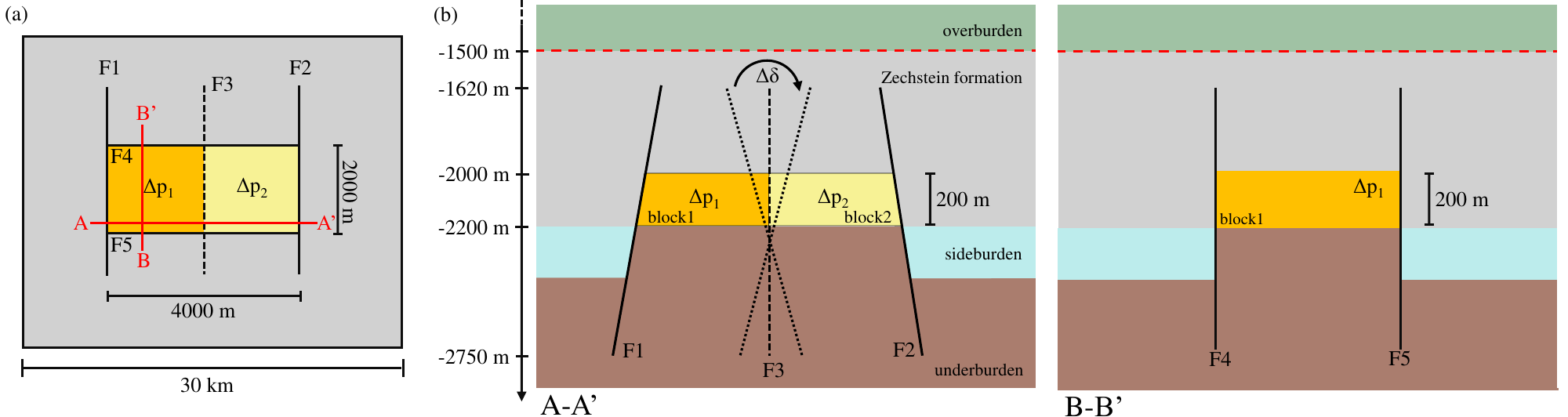}\hfill
    \caption{(a) horizontal cross-section of the conceptual model at 2100~m depth. Faults F1 and F2 are inclined and aligned along the y-axis; F4 and F5 are vertical and aligned along the x-axis. Fault F3 separates Block 1 and Block 2, which may experience different pressure histories ($\Delta p_1$, $\Delta p_2$).  (b) vertical sections of the conceptual model along the trace A-A and B-B. $\Delta \delta$ represents the possible variation in the dip angle of fault F3. This figure is not to scale to highlight the local features within the reservoir.}
    \label{fig:conceptual-3D}
\end{figure}

\subsection{Conceptual reference model}
\begin{figure}
    \centering
    \includegraphics[width=0.9\textwidth]{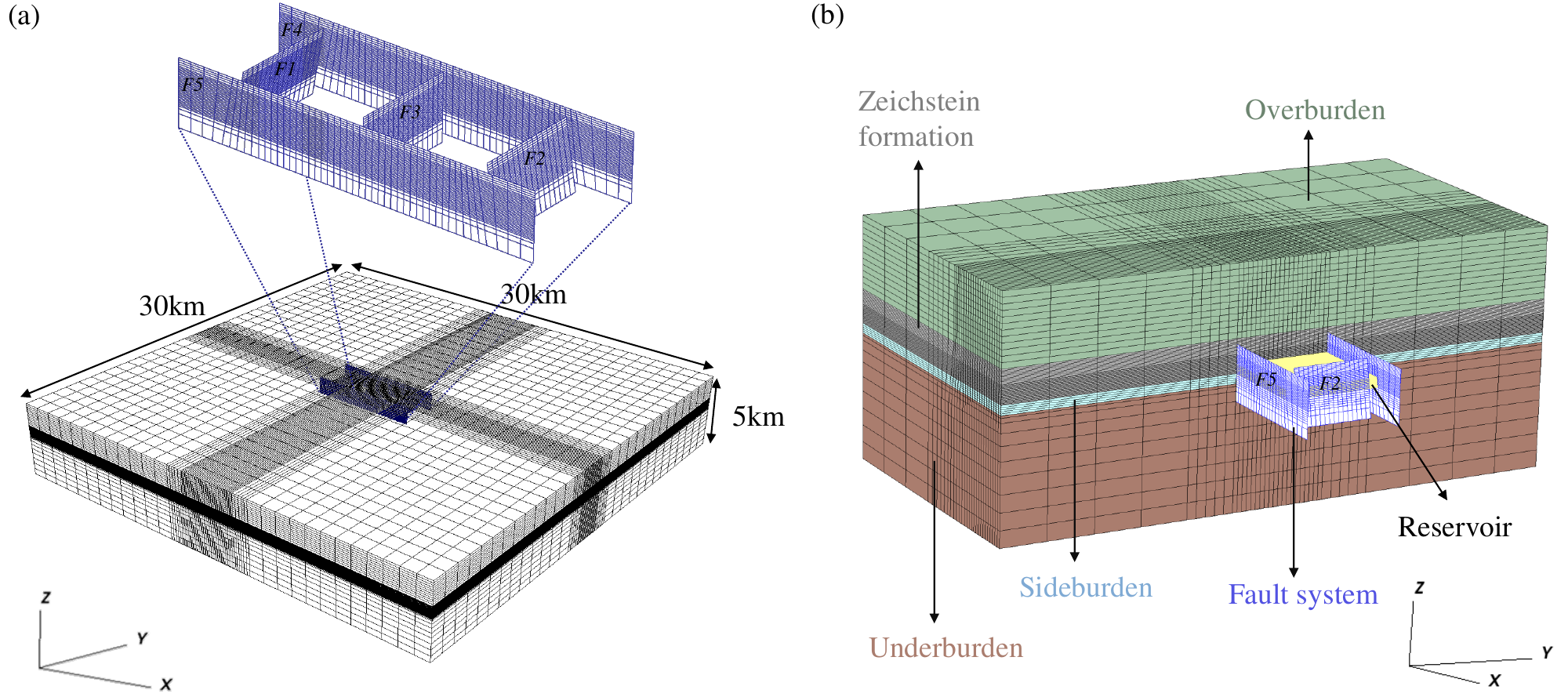}
    \caption{%On the left:
    (a): axonometric view of the 3D computational grid used for the geomechanical simulation and the embedded 2D grid used to represent the fault system. Right: cut of the model along a vertical plane of symmetry with the FE mesh grid on the back part and the discretization of fault planes with IE in the front part. The various colors represent the different portions of the domain in agreement with Figure~\ref{fig:conceptual-3D}.}
    \label{fig:axonometric_view}
\end{figure}

The domain geometry conceptually reproduces the compartmentalization commonly observed in \revB{many} UGS fields of the Rotliegend formation. The geometry and mechanical parameters were defined in agreement with the \cite{sodm_web}. Detailed information regarding the conceptual model is available in Franceschini et al. %Franceschini et al. (2024)
\cite{SodM1}.
The reservoir compartments are approximately located at the center of a 30$\times$30$\times$5~km domain %, to minimize the influence of boundary conditions on the area of interest.
(Figure~\ref{fig:conceptual-3D}~a, Figure~\ref{fig:axonometric_view}~a). %illustrates a view of the model.
The reservoir is 200-m thick, located at a depth of 2,000-m, and consists of two adjacent blocks measuring 2$\times$2~km each. % (Figure~\ref{fig:geol-seismic}~c). 
The two reservoir compartments are separated by a fault (F3) and confined laterally by two sets of orthogonal faults (F1-F2 along the y-axis, F4-F5 along the x-axis, as shown in Figure\ref{fig:conceptual-3D}~a). Depending on the sealing properties of fault F3, the two compartments can have a partial hydraulic connection. %Therefore, the distribution of pore pressure in space and time may differ.
The faults extend from 3,000-m to 1,600-m depth and terminate in the Zeichestein salt formation sealing the reservoir on top (Figure~\ref{fig:geol-seismic}~c and Figure~\ref{fig:conceptual-3D}~b). Faults F1 and F2 %are inclined respect to the vertical z-axis with
have a dip angle of $\pm 85^{\circ}$, while faults F4 and F5 are vertical. The dip angle of fault F3 can vary from $+65^{\circ}$ to $+90^{\circ}$ (vertical) and from $-90^{\circ}$~(vertical) to $-65^{\circ}$. %In the vertical direction,
Block~2 can be shallower or deeper, with a maximum offset of 200~m, with respect to the Block~1. %, equivalent to half of the entire reservoir thickness.
%The elevation of the reservoir and caprock formations remains constant in the direction perpendicular to faults F1, F2, and F3, (i.e., along the B-B direction of Fig~\ref{fig:conceptual-3D}). %\af{i.e., along the x-axis. Corretto???}.

%\subsection{FE-IE model}
The 3D domain is discretized by 8-node hexahedra, with 253,165 nodes and 236,208 elements. A finer discretization is used between 1,800 and 2,200~m depth. Here, the reservoir element characteristic size is 100$\times$100$\times$20~m.
%Based on the Lagrangian formulation proposed by \cite{franceschini2016novel}, 4-node quadrilateral interface elements (IE) have been integrated into the 3D geomechanical FE simulator to investigate the mechanical reactivation of faults.
The fault system embedded in the 3D grid is discretized by 5,215 zero-thickness contact elements. Standard boundary conditions with zero displacements on the outer and bottom surfaces are prescribed, whereas the top surface, representing the ground, has a stress-free condition. Figure~\ref{fig:axonometric_view}~b shows a view of the discretized domain with the embedded fault system. %embedded in the continuous 3D grid, in the reference case where F3 is vertical.

%The model has been initially applied to the so-called ``reference'' scenario for the fluids CH$_4$, CO$_2$, H$_2$, and N$_2$.
The setting of the reference scenario is established on the basis of the hydro-geomechanical properties of the Rotliegend formation~\citep{Bai_etal16,Bui_etal17,Hau_etal18,NAM16,Wee_etal14} and the frictional properties of faults intercepting the Rotliegend together with immediate overburden and underburden formations~\citep{hun2021,Hun_etal17,Hun_etal20}. Fault F3 is kept vertical, %configuration is characterized by a vertical central fault F3,
with no offset between the reservoir blocks.
\sbal{In agreement with available data of the Rotliegend formation \citep{pijn19,hol2018,Zhao24}, the entire domain is modeled with a poroelasto-plastic constitutive law using the material properties listed in Table \ref{tab:ref_geomechanics}. Notice that overburden, underburden, and caprock behave elastically, as no significant stress change is computed.} % Cambiare in elastico? 
%Overburden, caprock, reservoir, and underburden are assumed to behave %linear
%elastically, with the %geomechanical
%properties reported in Table~\ref{tab:ref_geomechanics}.
%The failure criterion governing the fault reactivation is characterized by
Faults are governed by Coulomb's failure criterion with
cohesion $c = 2$~MPa and static friction angle $\varphi_s = 30^\circ$.
A linear slip-weakening is also possible with a reduction to the dynamic friction angle $\varphi_d$ at a sliding amount equal to $d_c$ \citep{SodM1}.
As to the undisturbed stress regime, the vertical effective stress, $\sigma_v$, is a principal component. The other two horizontal principal components of the effective stress tensor, $\sigma_H$ and $\sigma_h$, are characterized by confinement factors $M_1 = \sigma_h/\sigma_v = 0.64$ and $M_2 = \sigma_H/\sigma_v = 0.83$, with $\sigma_h$ orthogonal to F1, F2, and $\sigma_H$ orthogonal to F4 and F5. 
The value of $\sigma_v$ is computed from the formation density, assuming a hydrostatic pore-pressure distribution, \revB{and corresponds to $\sigma_v \approx 25$~MPa at the average reservoir depth of about 2100~m \citep{SodM1}}.
% In the numerical model, this undisturbed stress state is applied uniformly across all formations. We note that this is a simplification, as the Zechstein Group salt typically exhibits near-lithostatic horizontal stress (e.g., van~Eijs et~al.,~2006).

%------------------
\subsection{Fluid-dynamic model}  % Spostato qui da Results
\begin{figure}%[htbp]
    \centering
    \includegraphics[width=0.9\textwidth]{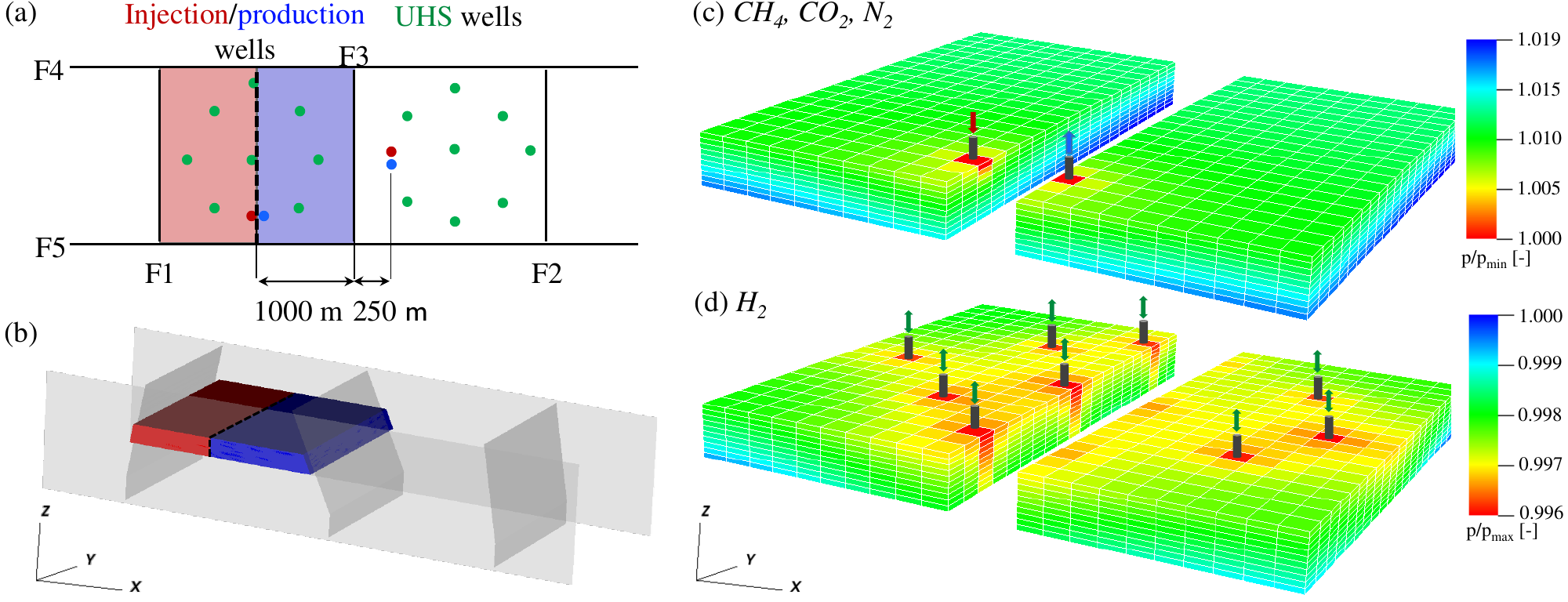}
    \caption{(a) spatial location of the injection and production wells (blue and red for injection/production for CH$_4$,  red for CO$_2$ and N$_2$ injection, and green for UHS cycles). (b) 3D perspective view of block 1, showing the surrounding fault planes and the subdivision into two sub-blocks, left (red) and right (blue), illustrated in the right-hand panels Right: normalized pressure distribution from the flow-dynamic simulation of CH$_4$ (c) and H$_2$ (d) within the right reservoir block. The grid splits vertically along the dashed line shown on the left panel to provide evidence of the pressure distribution along the vertical direction.}
    \label{fig:remap}
\end{figure}

The gas flow dynamic during injection and withdrawal is simulated by different codes according to the selected fluid: OPM Flow \citep{OPM2015, rasmussen_2021_the} for CH$_4$ (UGS), and Eclipse \citep{schlumberger_2014_eclipse} for 
H$_2$ (UHS),  CO$_2$ (CCS), and N$_2$. For all simulations, the same well configuration is used \citep{sodm_web}. 
During the PP stage, two wells located next to the reservoir boundary are open in the first layer (see Fig.~\ref{fig:remap}~a,b,c,d). Then, these wells are converted to injection purposes. Only for UHS, nine storage wells serve for injection and withdrawal. These wells are rather uniformly distributed throughout the reservoir, as there are no structural or petrophysical characteristics to suggest specific well locations. 
In Figure~\ref{fig:remap}~a, the location of the injection and withdrawal wells is shown with respect to the fault system. To avoid interpolation between different computational grids, both OPM and Eclipse finite volume cells coincide with the hexahedral elements of the geomechanical model. 
Each reservoir block is subdivided into 4,000 regular cells with a $20\times20\times10$ partitioning. 
The horizontal and vertical permeability is $k_h = 600$~mD and $k_v = 300$~mD, respectively.
\revB{In the flow simulations, only the reservoir layers are included in the hydraulic domain. The caprock, sideburden and underburden are not part of the flow model, and therefore their pore pressure is kept at the initial hydrostatic value throughout all simulation stages. As a consequence, all pressure variations produced during production and injection remain confined within the reservoir blocks, and the pressure field provided to the geomechanical model reflects changes occurring only inside the reservoir.}

The PP stage is conventionally set to 10 years, with a maximum pore pressure drop of about 20~MPa. 
The targeted final Recovery Factor, i.e., the ratio between the cumulative produced volume of gas and the gas originally in place, is 90\%, in line with theoretical reference behavior for a depletion-driven gas reservoir. 
Then, for CH$_4$ the reservoir pressure recovers to the initial value in two years during the CGI phase and a set of UGS cycles, consisting of 6-month withdrawal and 6-month injection each with a 10-MPa pressure, follows.  
As to the other fluids, the pressure history is targeted to closely approximate the maximum static pressure drop and recovery in reservoir observed for CH$_4$, imposing realistic storage constraints. In addition, no injection/withdrawal cycles are simulated for CO$_2$ and N$_2$ storage.
The CO$_2$ storage scenario considers a maximum injection rate of $4 \cdot 10^6$ Sm$^3$/day per well for 13 years before recovering the initial pressure value. Injection of N$_2$ lasts 6 years with an injection rate of $4 \cdot 10^6$ Sm$^3$/day per well. 
Concerning H$_2$, 10 years of injection/withdrawal activity are simulated after a 8-year refilling phase at a maximum rate of $4 \cdot 10^6$ Sm$^3$/day per well. 
Each UHS cycle consists of a 6-month withdrawal and a 6-month injection period, with a maximum pressure variation of 10~MPa. 
In all simulations, a 1-year loading step (l.s.) is used during PP and CGI/refilling, while during UGS and UHS the loading step is equal to 15~days.

% Fault lubrification in non sealing faults
It is important to recall that in the Netherlands the fluid pressure is not allowed to exceed the initial value $P_i$. %Operational data from Dutch UGS sites indicate that reinjection pressures are typically kept within a narrow safety margin of about 0–1~MPa below $P_i$ (???) \citep{TNO15}.
Pore pressure propagation within the fault zone varies according to the distribution of pressure changes across different reservoir compartments. When differing pressure changes occur on opposite sides of the fault, the resulting pore pressure change within the fault, denoted as $\Delta P_f$, correspond to the average of the pressure changes from both sides.
The fluid pressure in the fault, $P_f$, generates forces that oppose to the action of the pressure variation on $\sigma_n$.

\subsection{Sensitivity analysis}

\begin{table}
\centering
\small
\begin{tabular}{lccc}
Layer & Density [kg/m$^3$] & Young Modulus [GPa] & Poisson ratio [-] \\
\hline
Overburden      & 2200       & 10.0   & 0.25 \\
Zechstein Salt  & 2100       & 40.0   & 0.30 \\
Reservoir (Upper Rotliegend) & 2400   & 11.0  & 0.15 \\
Underburden     & 2600       & 30.0   & 0.20
\end{tabular}
\caption{Formation-dependent geomechanical parameters.}
\label{tab:ref_geomechanics}
\end{table}

A sensitivity analysis has been developed within the range of some geological and geomechanical quantities to investigate the configurations that are more likely to generate ``unexpected'' fault reactivations during the injection and storage of different kinds of fluid. %All the investigated scenarios are then tested by the 3D geomechanical simulator.
The analysis addresses the role of: (i) the geological setting, by varying geometry and %the hydro-geomechanical
properties of faults and reservoir rock, (ii) the poro-elastic stress change with respect to the natural stress regime, and (iii) the space and time distribution of pore pressure gradients in the reservoir compartments. %within the faults bounding/compartmentalizing the reservoir.
%Our goal is to understand which are the forcing factors  that are prone to induced ``unexpected'' (micro)-seismicity during the injection and storage phases for different kind of fluids.
%
The analysis has been carried out in two stages: %in order to completely analyze and assess the influence of numerous factors on the behavior and performance of the system under consideration. In the initial stage,
\begin{enumerate}
\item one parameter at a time is modified with respect to the reference scenario with CH$_4$ storage/production.
This enables a thorough analysis of the individual influence of each parameter. %on the UGS system (i.e., CH$_4$ production and storage).
The considered variables and scenarios are summarized in Table~\ref{tab:sensitivity};
%
% \iffalse    % Commented because is the same thing as the next one
% \begin{table}
% \begin{tabular}{llll}
%  & \textbf{Parameter} & \textbf{Reference} & \textbf{Interval} \\ \hline
% \multirow{2}{*}{\begin{tabular}[c]{@{}l@{}}Reservoir and \\ fault geometry\end{tabular}} & F3 dip & 90 & +65, -65 \\
%  & Compartment offset & 0m & 100m \\
% \multirow{2}{*}{\begin{tabular}[c]{@{}l@{}}Initial stress\\ regime\end{tabular}} & $\theta$ & 0$^\circ$ & 90$^\circ$ \\
%  & M1,M2 & 0.74,0.83 & 0.40, 0.47 \\
% \multirow{4}{*}{\begin{tabular}[c]{@{}l@{}}Mohr-Coulomb \\ criterion\end{tabular}} & c & 2~MPa & 0, 10~MPa \\
%  & $\varphi_s$ & 20$^\circ$ & 20$^\circ$ \\
%  & $\varphi_d$ &  & 10,20$^\circ$ \\
%  & $d_c$ &  & 2mm, 20mm \\
% \begin{tabular}[c]{@{}l@{}}Reservoir \\ stiffness\end{tabular} & E & 11GPa & 8GPa, 20Gpa \\
% \begin{tabular}[c]{@{}l@{}}Caprock \\ formation\end{tabular} & caprock & linear elastic & viscous \\
% Biot coefficient & $\alpha$ & 0.86 & 1.0 \\
% \begin{tabular}[c]{@{}l@{}}UGS pressure \\ change\end{tabular} & $\Delta p_1$, $\Delta p_2$ & 10~MPa & 10 MPa, 20~MPa
% \end{tabular}
% \caption{Model parametrization for the ``reference test case and variability range of the various parameters tested in the sensitivity analysis.}
% \label{tab:sensitivity}
% \end{table}
% \fi
%
\item
%Subsequently, a secondary series of simulations have been run by
simultaneous variation %the reference set-up with
of more than one parameter, as identified %during the previously outlined sensitivity.
from Stage 1.
%The objective %behind this set of simulations
%is twofold: (i) to achieve a configuration that more accurately capure the UGS reservoir conditions in the Netherlands, and (ii) to investigate, from a geomchanical perspective, the feasibility of UGS using different types of gases.
%More specifically, several combinations of cohesion, static friction angle, compartment offset, viscosity, and reservoir stiffness have been tested.
%The selected combination settings (Table~\ref{tab:combinations}) were tested with different injected fluids.
%For the four fluids under consideration,
%We used the combinations of settings summarized in Table \ref{tab:combinations}. %Furthermore, we accounted for geo-
To reflect potential mechanical consequences of fluid–rock chemical interactions, the Young modulus was modified by ±30\% depending on the injected fluid type \citep{harbert_2020}.
\end{enumerate}
Note that the objective of this analysis is to reproduce realistically the actual settings of the Rotliegend formation,
%It should be noted that only possible configurations have been evaluated since the goal of the study is to examine under which conditions this technology is feasible
rather than to investigate ``extreme'' conditions.
%By conducting the sensitivity analysis described above, we were able to explore the main mechanisms that may cause unexpected seismic events, and assessed the safety, in terms of criticality, of injection or storage phases.

%\begin{table}
%\begin{tabular}{ll}
%\hline
%Scenario ID & Parameter/Mechanism \\
%\hline
%1 & reference \\
%% 2 & Biot coefficient $\alpha$ = 1 \\
%3a & Fault F3 with dip angle $\delta=+65^\circ$ \\
%3b & Fault F3 with dip angle $\delta=-65^\circ$ \\
%3c & Block 2 moved down by the offset $o=100$ m \\
%3d & Block 2 moved down by the offset $o=200$ m \\
%4a & Principal stress $\sigma_H$ and $\sigma_h$ rotated by $\theta=90^\circ$ \\
%4b & Ratios $M_1 = \sigma_h/\sigma_v = 0.40$ and  $M_2 = \sigma_H/\sigma_v =0.47$ \\
%5a & Fault cohesion $c= 0.5$~MPa \\
%5b & Fault static friction angle $\varphi_s = 20^\circ$ \\
%5c & Linear slip-weakening with $\varphi_d$ = 10$^\circ$ and $d_c=2$~mm \\
%5d & Linear slip-weakening with $\varphi_d$ = 20$^\circ$ and $d_c=20$~mm \\
%6a & Reservoir Young's modulus $E = 8$ GPa \\
%6b & Reservoir Young's modulus $E = 20$ GPa \\
%7a & Uneven pressure variation during UGS $\Delta P_1 = -10$~MPa, $\Delta P_2 = 0$~MPa \\
%7b & Uneven pressure variation during UGS $\Delta P_1 = -10$~MPa, $\Delta P_2 = -20$~MPa \\
%8 & Pressure variation during UGS $\Delta P_1 = \Delta P_2 = -15$~MPa \\
%9 & Viscous caprock \\ \hline
%\end{tabular}
%\caption{Scenarios addressed in Stage 1 of the sensitivity analysis.}
%\label{tab:sensitivity}
%\end{table}

\begin{table}
\centering
\small
\begin{tabular}{ccl}
\hline
Paper section & Scenario & Parameter/mechanism \\
\hline
\ref{sec:ref} & 1 & reference \\
\hline
\multirow{4}{*}{\ref{sec:geom}}
& 2a & Fault F3 with dip angle $\delta=+65^\circ$ \\
& 2b & Fault F3 with dip angle $\delta=-65^\circ$ \\
& 2c & Block 2 moved down by the offset $o=100$ m \\
& 2d & Block 2 moved down by the offset $o=200$ m \\
\hline
\multirow{6}{*}{\ref{sec:param}}
& 3a & Fault cohesion $c= 0.5$~MPa \\
& 3b & Fault static friction angle $\varphi_s = 20^\circ$ \\
& 3c & Linear slip-weakening with $\varphi_d$ = 10$^\circ$ and $d_c=2$~mm \\
& 3d & Linear slip-weakening with $\varphi_d$ = 20$^\circ$ and $d_c=20$~mm \\
& 4a & Reservoir Young's modulus $E = 8$ GPa \\
& 4b & Reservoir Young's modulus $E = 20$ GPa \\
\hline
\multirow{2}{*}{\ref{sec:press}}
& 5a & UGS uneven $\Delta P$: $\Delta P_1 = -10$~MPa, $\Delta P_2 = 0$~MPa \\
& 5b & UGS uneven $\Delta P$: $\Delta P_1 = -10$~MPa, $\Delta P_2 = -20$~MPa \\
\hline
\multirow{2}{*}{\ref{sec:stress}}
& 6a & Principal stress $\sigma_H$ and $\sigma_h$ rotated by $\theta=90^\circ$ \\
& 6b & Ratios $M_1 = \sigma_h/\sigma_v = 0.40$ and  $M_2 = \sigma_H/\sigma_v =0.47$ \\
\hline
\end{tabular}
\caption{Scenarios addressed in Stage 1 of the sensitivity analysis (UGS).}
\label{tab:sensitivity}
\end{table}

% How to interpret the result
%The model results will provide comprehensive information regarding the 3D displacement and stress fields everywhere of the 3D grid.
The main interest of this study is focused on the fault behavior, and in particular the stress conditions yielding a potential sliding. %of these discontinuity surfaces deserve attention. By employing this modeling approach, it becomes possible to quantitatively evaluate the safety of UGS operations in relation to fault reactivation. The constitutive law describing the behavior of faults is the well-known Coulomb failure criterion \citep{labuz_2012_mohrcoulomb}:
According to the classical Coulomb failure criterion, fault stability is guaranteed if the modulus $\tau$ of the shear stress is smaller than its limiting value $\tau_L$:
\begin{equation}
  \tau < \tau_{L} = c + \mu(\sigma_n - P_f),
  \label{eqn:Coulomb}
\end{equation}
with $\mu$ the friction coefficient, i.e., the tangent of the friction angle, $\sigma_n$ the normal compressive stress, and $P_f$ the fluid pressure in the fault.
%This criterion compares the shear stress modulus acting on the fault surfaces ($\|\boldsymbol{\tau}\|_2$) with respect to the shear stress limit ($\tau_L$), function of the normal stress ($\sigma_n$), the friction angle ($\varphi_s$) and the cohesion ($c$).
Wherever $\tau=\tau_L$, sliding is allowed with a potential energy release, while the condition $\tau>\tau_L$ is forbidden.
Three parameters are used to evaluate quantitatively the fault state: (i) the criticality index $\chi=\tau/\tau_L$, and, where $\tau=\tau_L$, %which represents the ratio between actual tangential stress modulus $\|\boldsymbol{\tau}\|_2$ acting on a fault and the limit tangential stress according with the Coulomb criterion,
(ii) the reactivated area %(or thickness)
$t_{a}$, and (iii) the sliding $d$. Clearly, from Eq.~\eqref{eqn:Coulomb} we have $\chi$ $\in$ [0, 1]. The closer $\chi$ to 1, the more likely is a fault activation. %When the shear stress equals the shear stress limit, the fault is prone to slip.
%Additionally, when a tensile normal stress is present, the fault has the potential to open.
In particular, in the sensitivity scenarios %addressed within the sensitivity analysis,
the maximum criticality factor ($\chi_{\max}$), the maximum sliding ($d_{\max}$), and the fault area ($t_{80}$) where $\chi$ exceeds 0.8 %(i.e., the ration $A_{80}/L$ where $A_{80}$ is the fault area with $\chi > 0.8$ and $L$ is the fault length in the horizontal plane) enable a comprehensive comparison of
are used to lump
the faults' behavior. %Notably, these quantities are closely related. Indeed, the single element can slide only when $\chi = 1$, while $t_{80}$ gives an indication of the largest fault size close to the reactivation ($\chi > 0.8$).
%
%The total number of scenarios and combinations is pretty large. For this reason, in the sequel
%Given the huge number of simulations performed,
%only the most significant simulations will be reported and discussed.

\section{Analyses and results}
\subsection{Fluid-dynamic model}

\revB{The results of the fluid-dynamic model indicate that, during CH$_4$ storage, the pressure change within each UGS cycle remains nearly uniform across the reservoir block during both withdrawal and injection phases. This behavior is consistent with the spatial pressure distributions shown in Figure~\ref{fig:remap}~c,d, where the normalized pressure fields for CH$_4$ and H$_2$ exhibit limited lateral variations within the same block. The vertical cuts displayed in Figure~\ref{fig:remap}~c,d further confirm that the pressure variations are essentially homogeneous across the block thickness.}

\begin{figure}%[htbp]
    \centering
    \includegraphics[width=1.0\textwidth]{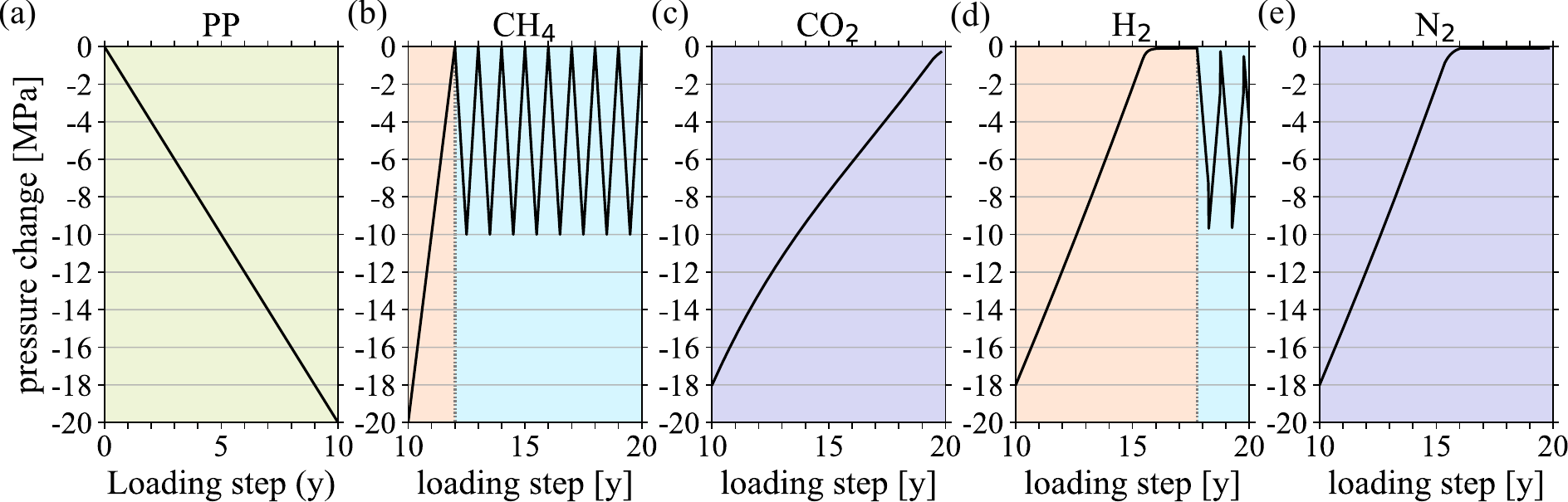}
    \caption{Schematic evolution of pore-pressure variation for the investigated fluids.(a) Primary production (PP): 10-year pressure decline of about 20~MPa. (b) CH$_4$: 2-year cushion-gas injection (CGI) followed by UGS cycles consisting of 6-month withdrawal and 6-month injection (10~MPa swing). (c) CO$_2$: 10-year pressure recovery toward $P_i$. (d) H$_2$: 8-year refilling phase up to $P_i$, followed by 6-month withdrawal and 6-month injection (10~MPa swing) UHS cycles (e) N$_2$: 10-year injection stage reaching $P_i$.}
    %under investigation, i.e., CH$_4$, CO$_2$, H$_2$, N$_2$.}
    \label{fig:pmean}
\end{figure}

\revB{Figure~\ref{fig:pmean} illustrates the imposed pore-pressure evolution for all investigated fluids. 
For CH$_4$, the pressure returns to the initial value following the PP stage and the subsequent CGI phase, after which the prescribed UGS cycles generate the imposed 10-MPa oscillations (Figure~\ref{fig:pmean}~a,b) . 
For CO$_2$ and N$_2$, the pressure increases according to the imposed injection rates until the initial pressure is recovered (Figure~\ref{fig:pmean}~c,e). 
For H$_2$, the imposed refilling and cycling schedule leads to the 10-MPa pressure variations associated with the UHS cycles (Figure~\ref{fig:pmean}~d).}

\subsection{Mechanisms of fault re-activation: reference case (Scenario \#1)}
\label{sec:ref}
The basic fault reactivation mechanisms are studied in Scenario~\#1, which serves as the reference scenario for UGS activity with CH$_4$. \revB{Scenario~\#1 does not represent the full structural complexity of real Rotliegend reservoirs; rather, it is a simplified setting adopted to clarify the fundamental mechanical processes governing fault reactivation during PP, CGI and UGS. The choice of this configuration was made in consultation with \cite{sodm_web}, and aims to reproduce the essential geological, stress and operational conditions of a typical Rotliegend UGS reservoir while maintaining a geometry that allows the isolation of the primary mechanical mechanisms before exploring more complex configurations. For this reason, Scenario~\#1 serves as the "reference" against which all subsequent parametric variations are interpreted.}
A few preliminary outcomes can already be found in previous studies \citep{Tea_etal19, Tea_etal20,SodM1}. For the sake of completeness, they are briefly summarized here as well.

\begin{figure}
    \centering
    \includegraphics[width=0.8\linewidth]{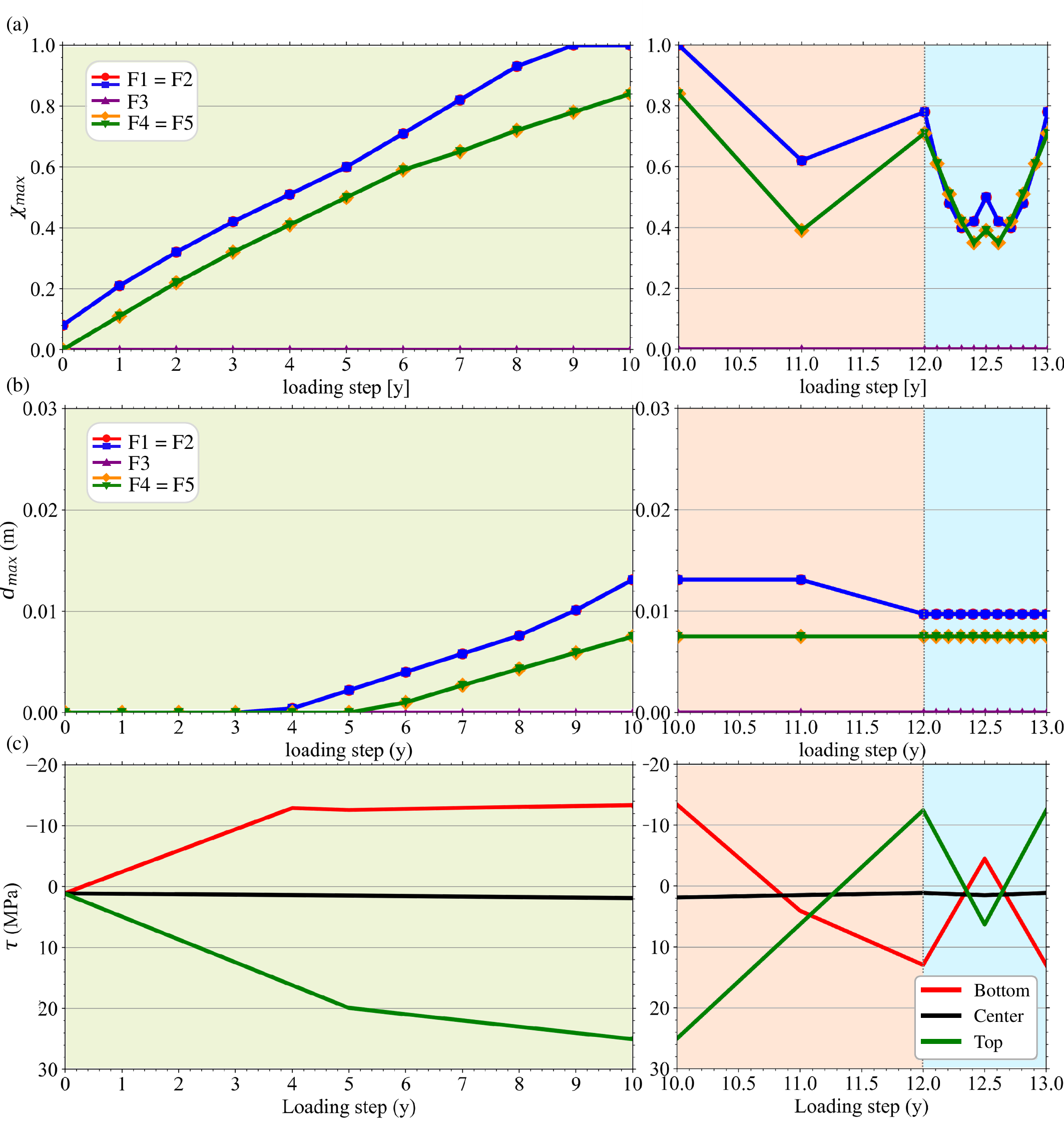}
    \caption{\revB{Scenario \#1 ("Reference"): $\chi_{max}$ (a) and $d_{max}$ (b) over time for all faults. (c) behavior of $\tau$ for three nodes located at the top, bottom, and center of F1. Positive values for the shear stress means that it is upward oriented. These plots summarize the reactivation mechanisms during PP, CGI and UGS in the simplified reference configuration.}}
    \label{fig:CH4-reference}
\end{figure}

Figure~\ref{fig:CH4-reference}~a shows the behavior of the maximum $\chi$ value, $\chi_{max}$, experienced on faults F1 through F5 for 13 years, i.e., PP, CGI and the first UGS cycle.
Faults F1 and F2 can reach a critical condition at the end of the PP stage.
During CGI and UGS, the fault activation risk reduces, with $\chi_{max}$ decreasing to about 0.80. %, which is considered the threshold for fault activation.
During CGI, $\chi_{max}$ initially decreases, i.e., the fault returns more stable, but then increases again, with the shear stress acting
in the opposite direction with respect to that experienced during PP (Figure~\ref{fig:CH4-reference}~c - orange panel).
The time behavior of the maximum sliding  $d_{max}$ computed for each fault is also shown in Figure~\ref{fig:CH4-reference}~b.
These two quantities, i.e., $\chi_{max}$ and $d_{max}$, are closely related to each other, since a single element can slide only when $\chi=1$.
%However, we prefer to propose an averaged version of $\chi$, so as to obtain information on the criticality state of the entire fracture at a given depth. In particular, for the sake of clarity and ease of readability, $\chi$ is represented for each fault as a function of depth only, i.e., for each $z$-value we compute the $\chi$ average for the stripe of fault elements located at the same depth.
Notice that, to obtain information on the criticality state of the entire fracture at a given depth, the provided $\chi_{max}$ represents the value averaged on the stripe of fault elements located at the same depth. For this reason, it is possible to detect sliding conditions even when the $\chi$ value reported on Figure~\ref{fig:CH4-reference}~a is smaller than 1.
Fault F1, F2, F4, and F5 can start sliding between year 4 and 6. The maximum sliding can be of about 1.4~cm.
Fault F3 remains unaffected throughout the loading and injection steps due to the model symmetry. % with its criticality index always equal to 0, meaning no likelihood of reactivation. 
\revB{This behavior, however, directly reflects the simplified structural setting adopted in Scenario~\#1. In this configuration, fault~F3 is perfectly vertical and exhibits no offset, so no shear stress builds up during PP or CGI/UGS. As a consequence, the bounding faults (F1–F2 and F4–F5) display the largest shear variations, since they juxtapose a depleting reservoir block against sideburden regions.}
%In real Rotliegend reservoirs, intra-reservoir faults are rarely perfectly vertical and commonly present significant offset or compartmentalization, conditions under which they may experience the most destabilizing stress changes. The resulting hierarchy of fault criticality observed in Scenario~\#1 is therefore specific to this idealized geometry and motivates the more complex structural configurations investigated in the subsequent scenarios.}

\revA{The behaviour of $\chi$ with respect to depth on all the fracture surfaces at loading step 10 is shown in Figure~\ref{fig:ref_chi_vs_depth} in the Supplementary Materials. These profiles indicate that the most critical conditions develop along the top and bottom of the reservoir. }
% Figure~\ref{fig:ref_chi_vs_depth} shows the behavior of $\chi$ with respect to depth on all the fracture surfaces at loading step 10. The most critical condition develops along the top and bottom of the reservoir.

%\begin{figure}%[htbp]
%    \centering
%    \vcenteredhbox{\includegraphics[width=1\textwidth]{fig_revision_GETE/X_vs_z_upd_rev01.pdf}}
%    \caption{ Scenario \#1 ("Reference"): depth profiles of the criticality index $\chi$ on all faults at the end of PP (loading step 10). High criticality develops at the reservoir top and bottom.}
%    \label{fig:ref_chi_vs_depth}
%\end{figure}

The mechanisms leading to the possible fault reactivation can be better understood in terms of shear stress $\tau$ (\ref{fig:CH4-reference}~c), particularly at loading steps 0, 10 (end of PP), 11, 12 (end of CGI), 12.5 (end of UGS withdrawal), and 13 (end of UGS injection). The largest values are computed at the end of PP, with opposing shear stress orientations at top and bottom of the reservoir. During PP, the direction of the shear stress is always oriented towards the reservoir mid-plane. At loading step 4 and 5 sliding starts (Figure~\ref{fig:CH4-reference}~b) at reservoir top and bottom.
The elements surrounding the activated fault portions increase their shear stress to accommodate the excess to the Coulomb frictional limit not supported by the activated (sliding) elements. At loading step 11, half of the pore pressure change has been recovered. As the reservoir expands due to pressure recovery, the shear stress changes direction and annihilates approximately the shear stress acting on the previously sliding elements. The reservoir continues to recover pressure and re-expand until loading step 12. During the second half of CGI, the shear stress increases again, with a direction opposite to that experienced during PP. \revB{The expansion during CGI can generate again a critical condition on the fault, with a possible localized reactivation at the reservoir top due to the stress redistribution following the sliding. During UGS, at the end of the production phase (loading step 12.5), the shear stress almost equals the condition at loading step 11.} Correspondingly, the shear stress at loading step 13 equals the stress state at loading step 12. \revB{Indeed, during UGS, no new sliding occurs. The faults keep the slip accumulated during PP and CGI, and the porous medium behaves according to a linear elastic constitutive law.} %Although the shear stress reverses sign in this phase, its magnitude remains below failure ($\chi_{max}$ < 1), preventing further reactivation.}
%This complex behavior is also confirmed by the analysis of the stress path experienced by some representative elements, as reported in \cite{SodM1}. 

\revB{These results, together with the stress-path interpretation developed in \cite{SodM1}, provide the basis for understanding the occurrence of “unexpected’’ reactivation during CGI or UGS, even though the reservoir pressure remains below $p_i$. The slip accumulated during primary
production produces a permanent redistribution of stresses along the fault system. As a result, when pressure increases during CGI and UGS, the response does not follow the elastic unloading path of PP. Instead, the system evolves along a different stress trajectory, characterized by a reversal of the shear stress $\tau$ at the top and bottom of the reservoir (see Figure~\ref{fig:CH4-reference}~c). Because of this stress-path reversal,
the fault may approach the failure condition again even at pressures within the previously experienced range. The renewed criticality does not require exceeding the historical maximum loading, but it arises because the system returns toward failure from a different direction in the $(\sigma_n,\tau)$ space, driven by the stress redistribution inherited from PP slip \citep{SodM1}}
\revB{To illustrate this mechanism in simple terms, we consider the evolution of $\chi_{max}$ on the bounding faults F1 and F2 as an example.
During UGS, the imposed pressure variation is of the order of 10~MPa. If such a pressure cycle were applied to a pristine reservoir–fault system, the resulting 
criticality would remain moderate (e.g., $\chi_{max} \approx 0.6$), indicating stable conditions. In the modeled case, however, the same pressure variation acts on a fault system that has already experienced PP and partial pressure recovery. Slip occurring during PP on faults F1 and F2 produces a permanent redistribution of stresses, so that the subsequent CGI and UGS phases do not correspond to a simple elastic unloading–reloading of the initial state. As a result, the same $\Delta P$ applied during UGS leads to significantly higher criticality levels on these faults (e.g., $\chi_{max} \approx 0.8$), increasing the likelihood of approaching failure even though the 
stress magnitudes remain lower than those reached during PP.}

\subsection{Sensitivity analysis: Stage 1}
%In scenario \#1, UGS operations do not appear to stress the faults. However, a $\chi_{max}$ close to the threshold value even when pressure conditions are not typically associated with fault activity suggests that faults might reactivate unexpectedly even under ``normal'' pressure conditions.

Starting from the reference configuration of scenario \#1, we explore how variations in individual geometric or mechanical parameters can affect \revA{the model outcomes, including fault reactivation potential and the associated stress and slip patterns.} We conducted an extensive investigation, of which we report here the most significant results, referring to the scenarios summarized in Table \ref{tab:sensitivity} for the UGS analysis with CH$_4$. Other simulations were carried out, but are not discussed here for the sake of brevity.
%into multiple scenarios to better understand the factors influencing faults reactivation in the subsurface system.
The outcomes %of these scenarios, specifically for CH$_4$,
are cross-compared to identify the settings that make the subsurface system more prone to fault reactivation. For the analysis, we use the parameter $t_{80}$, i.e., the areal extent where $\chi\geq0.80$. %as it provides a reliable measure for our investigation and its evolution over time is comparable to that of the criticality index $\chi$.
Such areal extent is scaled by 2,000~m, which is the characteristic size of each reservoir block, so that $t_{80}$ has the size of a length.
For the sake of simplicity, we focus on faults F2 and F3, unless differently specified. %The criticality in Faults F4 and F5 is generally lower due to factors such as the initial stress regime (excluding scenario 4a), the vertical orientation of the faults, and the constant depth (without any offset) in the Upper Rotliegend along the direction orthogonal to F4 and F5.  Conversely, fault F2 represents the discontinuities bounding the reservoir, while F3 experiences the highest stress level under realistic conditions.

\subsubsection{Reservoir and fault geometry}
\label{sec:geom}
\begin{figure}[htbp]
    \centering
    \includegraphics[width=1\textwidth]{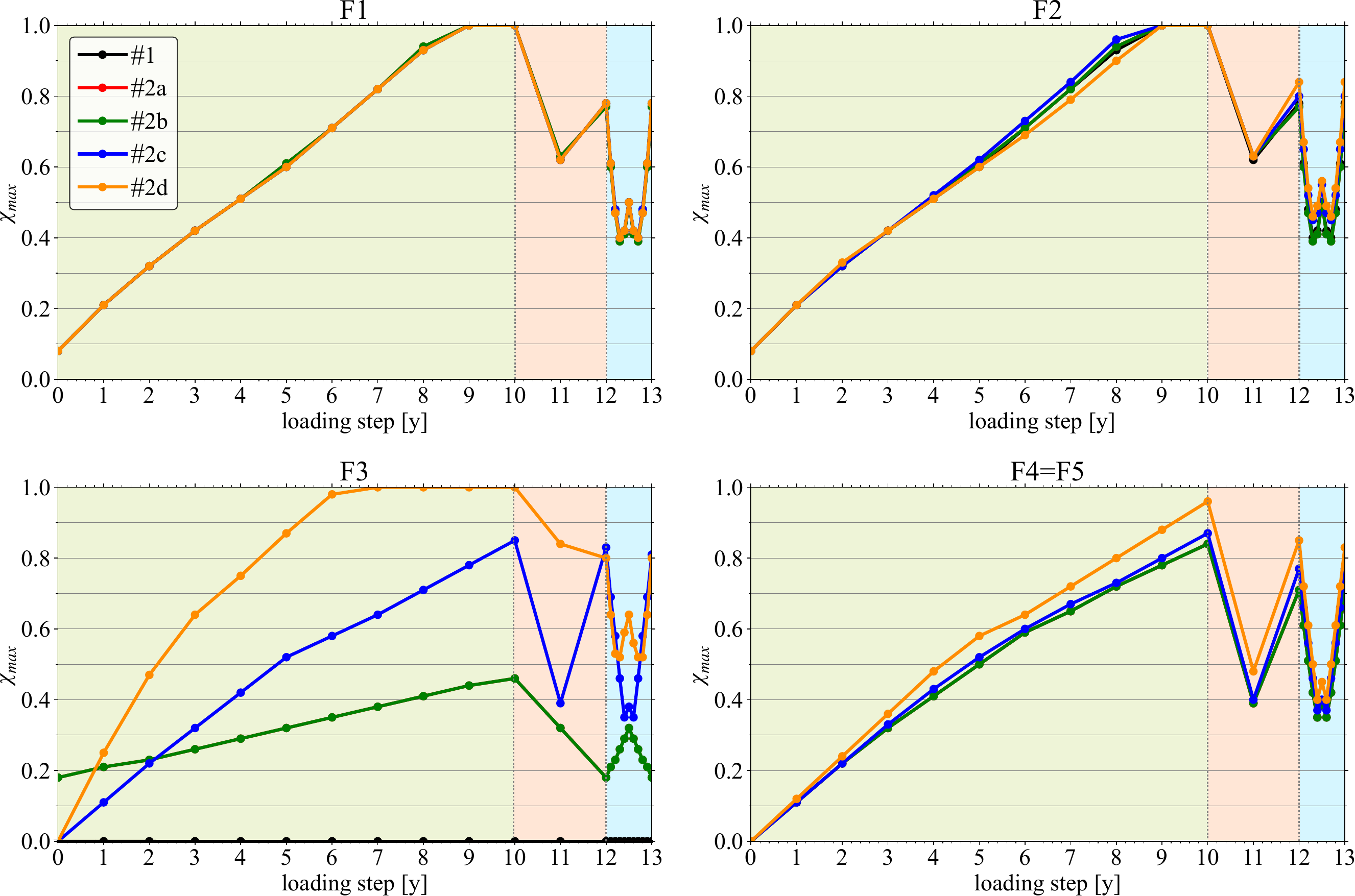}
    \caption{Scenarios \#2: $\chi_{max}$ over time  for each fault. Note that F4 and F5 exhibit an identical behavior, due to symmetric conditions. \revB{Scenario 2 includes four geometrical variants: (\#2a)  F3 dip $+65^\circ$, (\#2b) F3 dip $-65^\circ$, (\#2c) Block 2 lowered by 100 m, (\#2d) Block 2 lowered by 200 m.}}
    \label{fig:ref_chi}
\end{figure}

The impact of fault geometry is investigated by varying the dip angle of the central fault F3 (\#2a~-~\#2b) and increasing the offset between reservoir compartments of half (\#2c) and entire (\#2d) reservoir thickness. Figure~\ref{fig:ref_chi} shows $\chi_{max}$ over time for the reference and the \#2a~-~\#2d scenarios. Interestingly, while the geometry variation %, i.e., both F3 dip and compartment offset,
has no effect on fault F1, %because the variation does not change the loading conditions on these faults,
a large offset between the reservoir compartments (scenario \#2d) enhances the potential for instability in faults F4 and F5.
The instability of fault F2 is slightly increased by a variation of the fault geometry, however this becomes significant only when a large offset between compartments is introduced (\#2d).
% Add scenario 4 initial stress regime and #9 viscosity to the X picture

\begin{figure}%[htbp]
    \centering
    \includegraphics[width=0.45\textwidth]{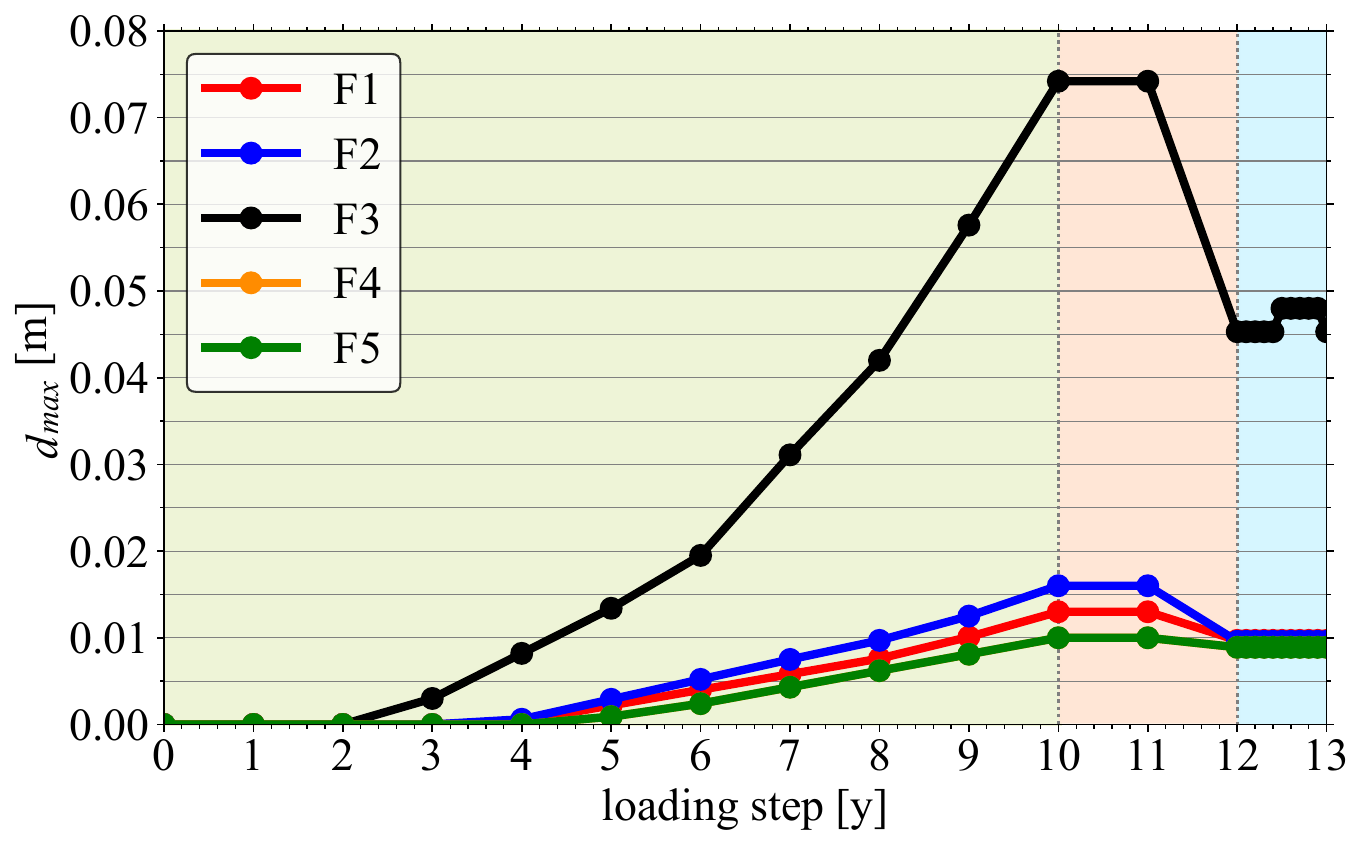}
    \caption{Scenario \#2d: maximum sliding $d_{max}$. \revB{Scenario \#2 includes (\#2a) F3 dip $+65^\circ$, (\#2b) F3 dip $-65^\circ$, (\#2c) Block~2 lowered by 100~m, (\#2d) Block~2 lowered by 200~m.}}
    \label{fig:CH4-sliding}
\end{figure}

In the latter configuration, unlike the reference scenario, fault F3 is significantly affected by the asymmetry of the system. Indeed, it becomes active and experiences the highest stress levels at the top and bottom of each compartment. %($z = -1800$~m for block 1 top, and from $-1960$~m to $-2040$~m, indicating simultaneous effects at the bottom of block 1 and the top of block 2, and $-2200$~m for the bottom of block 2).
%The behavior of $\chi$ in depth is shown in Figure~\ref{fig:ref_chi_vs_depth} at the end of the PP stage.
Differently from the reference scenario, $\chi_{max}$ on fault F3 can now reach the limiting value. %(i.e. $\chi_{max} = 1$)
Furthermore, during CGI and UGS, $\chi_{max}$ reaches approximately 0.8. Although this does not correspond to reactivation, it indicates that the fault is close to a critically stressed state.

Furthermore, $\chi_{max}$ achieves the alert threshold of 0.8 also during CGI and UGS, suggesting that UGS operations in this condition might be critical. %stress the faults in terms of activation.
Consistently, fault F3 exhibits the highest sliding value $d_{max}=7.4$~cm at the end of the PP phase (Figure~\ref{fig:CH4-sliding}). %This fact deserves particular attention as it resulted to be the maximum value of sliding obtained in this study.
Some instability can also occur at the end of the UGS cycle. %, and sliding increases again.
The other faults can activate around the loading step 4, %resulting at first in sliding of the interface elements
again at the reservoir top first, and at the reservoir bottom later.
The active elements within the embedded fault discretization at loading step 6, 10 (end of PP), and 12.5 (middle of a UGS cycle) are shown in Figure~\ref{fig:CH4-3_3D_active} in the Supplementary Materials.

%\begin{figure}%[htbp]
%    \centering
%   \includegraphics[width=1\textwidth]{figs/CH4/3-3d_activation}
%    \caption{Scenario \#2d: active interface elements (IE) on the fault system at loading step 6 (left), loading step 10 (center), and loading step 12.5 (right). Only the elements at the reservoir top, and then bottom, are activated.}
%    \label{fig:CH4-3_3D_active}
%\end{figure}

\subsubsection{Geo-mechanical parameters of faults and reservoir}
\label{sec:param}
The investigation continues further into the effect of Coulomb rupture criterion and reservoir stiffness. %and viscosity.
Concerning the Coulomb parameters, the variation of both the cohesion (\#3a) and friction angle (\#3b) exerts a significant influence on fault stability, %(excluding F3 because of symmetry conditions),
mainly during the PP and CGI phases. A reduced friction angle, in particular, weakens the resistance to fault reactivation.

\begin{figure}%[htbp]
    \centering
    \includegraphics[width=0.32\textwidth]{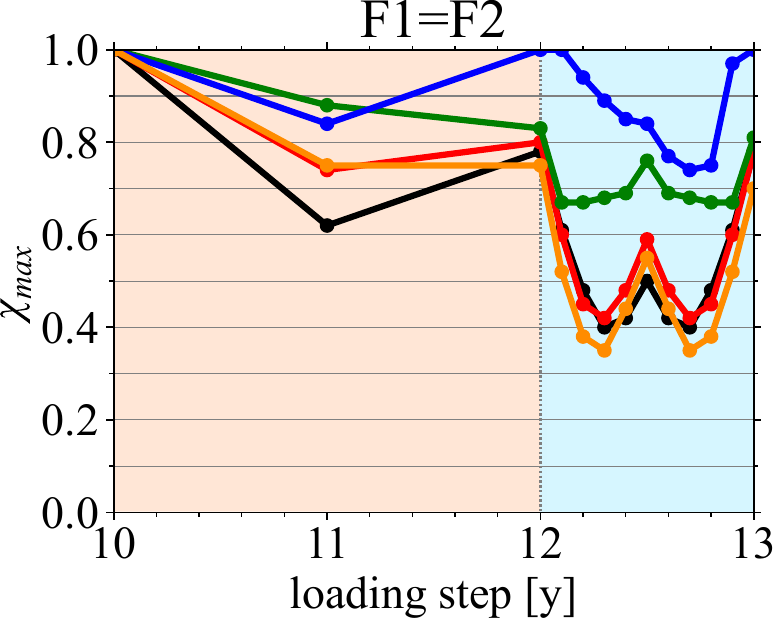}
    \includegraphics[width=0.32\textwidth]{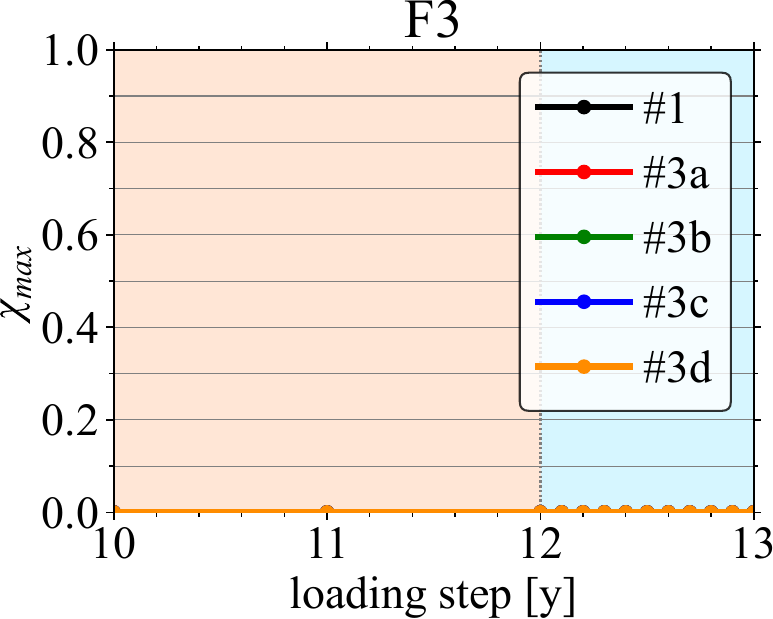}
    \includegraphics[width=0.32\textwidth]{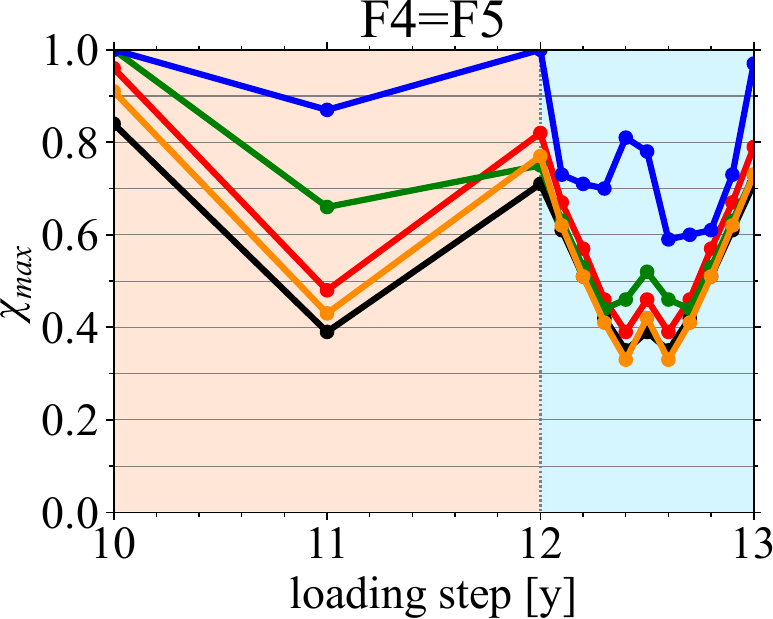}
    \caption{Scenarios \#3: $\chi_{max}$ over time for the CGI and UGS phases.}  \revB{Scenario \#3 includes four fault–strength variants: (\#3a) cohesion $c = 0.5$ MPa,
(\#3b) friction angle $\varphi_s = 20^\circ$, (\#3c) slip–weakening with $\varphi_d = 10^\circ$ and $d_c = 2$ mm, (\#3d) slip–weakening with $\varphi_d = 20^\circ$ and $d_c = 20$ mm. }
    \label{fig:CH4-5_chi_zoom}
\end{figure}

The mechanism of fault weakening is investigated by using slip-weakening constitutive law for the fracture behavior \citep{SodM1}:
\begin{equation}
  \varphi = \begin{cases}
      \varphi_s + \dfrac{\varphi_d - \varphi_s}{d_c} d & \text{for } d < d_c \\
      \varphi_d & \text{for } d \geq d_c
  \end{cases}
\end{equation}
The two parameters defining the new constitutive law are $ \varphi_d$ and $d_c$, i.e., the dynamic friction angle and the slip weakening distance, respectively. Slip weakening occurs after sliding begins by reducing the value of the friction angle. The most critical results in terms of fault stability are obtained for scenario \#3c, characterized by a reduction of the friction angle from $\varphi_s = 30^\circ$ to $\varphi_d = 10^\circ$ in a slip distance of $d_c = 2$~mm.
%The outcomes of this analysis, already shown by \cite{SodM1}, %and we briefly report them here for completeness.
%
The behavior of $\chi_{max}$ for scenarios \#3 during the CGI and UGS stages is reported in Figure~\ref{fig:CH4-5_chi_zoom}, showing that
%A comparison within the CGI and UGS phases in terms of criticality index varying the friction angles, cohesion and slip-weakening parameters shows that the new constitutive law
slip-weakening may
cause %F1, F2, F4 and F5 to slide also at the end of the CGI and UGS phases,
a reactivation also at loading steps 11 and 13,
but not at loading step 12.5, i.e., in the middle of a UGS cycle. %at the end of the 6-month UGS-P production phase (see the zoom in the right panel of Figure~\ref{fig:CH4-5_chi_zoom}).

\begin{figure}%[htbp]
    \centering
    \includegraphics[width=0.32\textwidth]{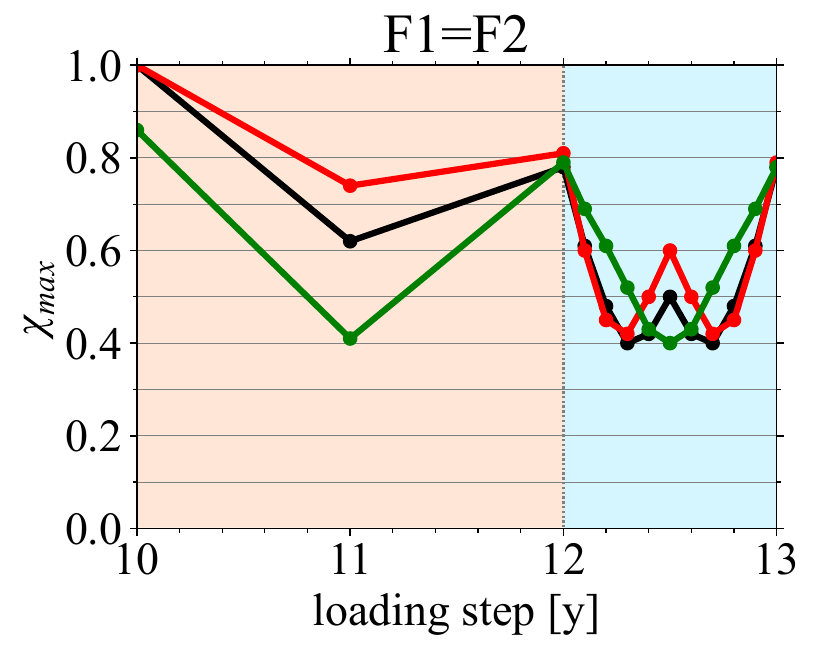}
    \includegraphics[width=0.32\textwidth]{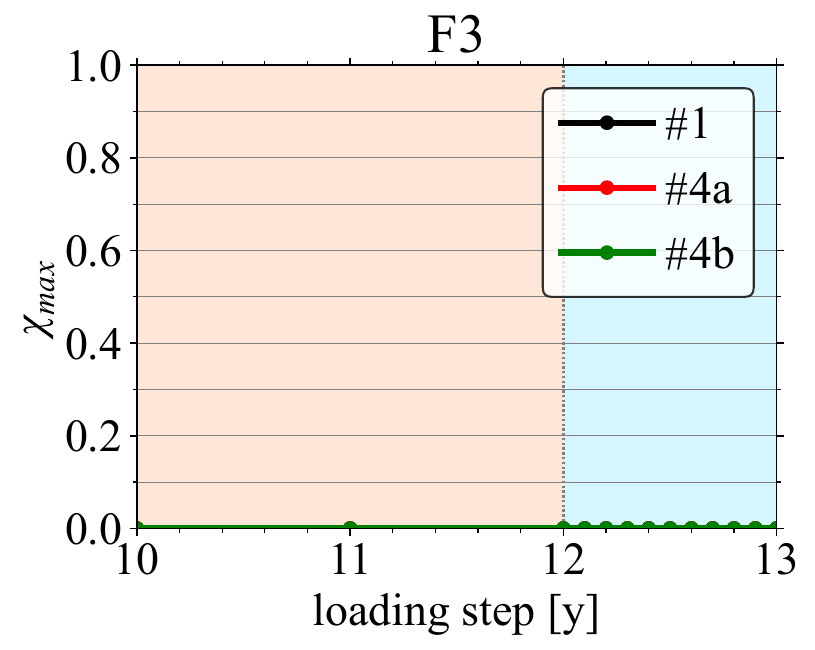}
    \includegraphics[width=0.32\textwidth]{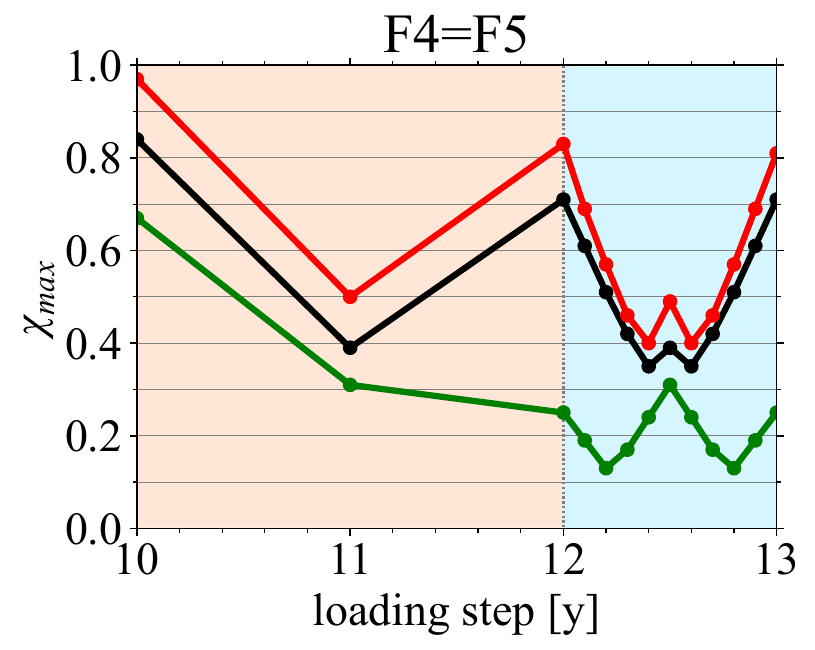}
    \caption{Scenarios \#4: $\chi_{max}$ over time for the CGI and UGS phases. \revB{Scenario \#4 includes (\#4a) Young's modulus $E = 8$ GPa,
(\#4b) Young's modulus $E = 20$ GPa.}} 
    \label{fig:CH4-6_chi_zoom}
\end{figure}

%Figure~\ref{fig:CH4-sliding} provides the behavior of the maximum sliding along the faults for the proposed scenario. It can be seen that the current sliding in scenario \#5c is more than twice what obtained using a static friction angle, as in the reference case.
%The stress path and the yield bound are quite complex due to weakening. Due to the very small friction angle ($\varphi = 10^\circ$), a large part of UGS is characterized by a stress state that develops either on the yield surface or very close to it.
%Because of the reduced friction angle, the stress path of an element located at the top of F2 shows that the yield surface is reached more easily than in scenario \#1 during PP, at the end of CGI and UGS-S.
%As observed for the reference scenario, the elastic phases develop with an almost constant normal stress because of the selected ratios between the reservoir and overburden stiffness and between the pressure change in the reservoir and within the fault. Indeed, the ``fortuitous'' combination of $E_{\text{Zechstein}} = 2E_{\text{reservoir}}$ and $\Delta P_{\text{reservoir}} = 2\Delta P_f$ contributes to the maintenance of an almost constant stress tensor within this stage.

% Reservoir stiffness
The relationship between reservoir stiffness and fault behavior has been further investigated in Scenarios \#4. A significant contrast between the reservoir and the surrounding caprock, sideburden, and underburden stiffness may play an important role to induce a critical stress state on faults.
In scenario \#4a, the reservoir is softer than the surroundings,
and this generally increases the likelihood of a fault reactivation (see Figure \ref{fig:CH4-6_chi_zoom}).
Conversely, in scenario \#4b, the fault reactivation likelihood is notably reduced. During the UGS stage, $\chi_{max}$ remains below 0.7, regardless of the reservoir stiffness value (see Figure \ref{fig:CH4-6_chi_zoom}).

\subsubsection{Uneven pore pressure excursion in the reservoir compartments}
\label{sec:press}
\begin{figure}%[htbp]
    \centering
    \includegraphics[width=0.32\textwidth]{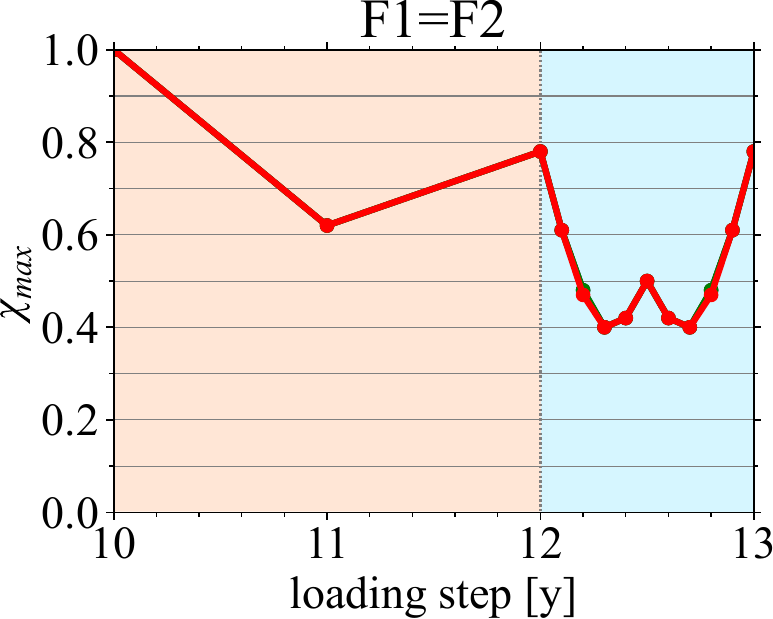}
    \includegraphics[width=0.32\textwidth]{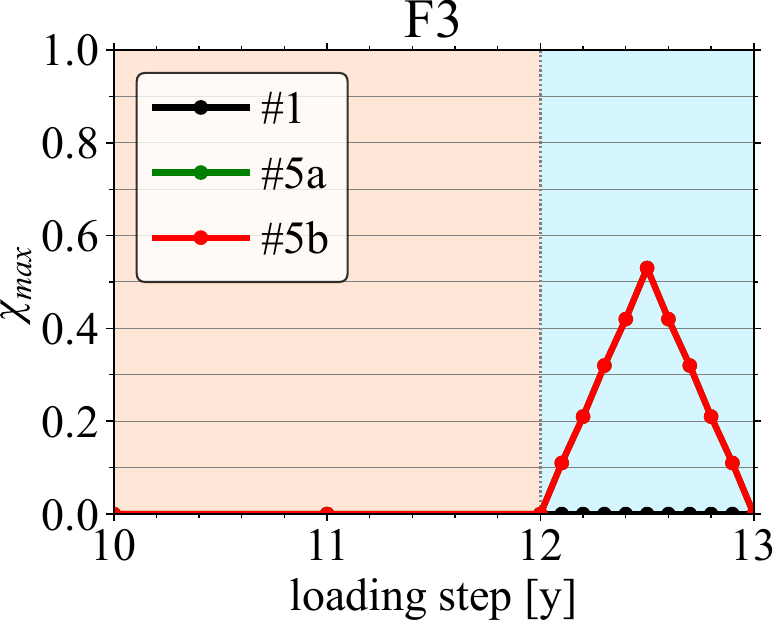}
    \includegraphics[width=0.32\textwidth]{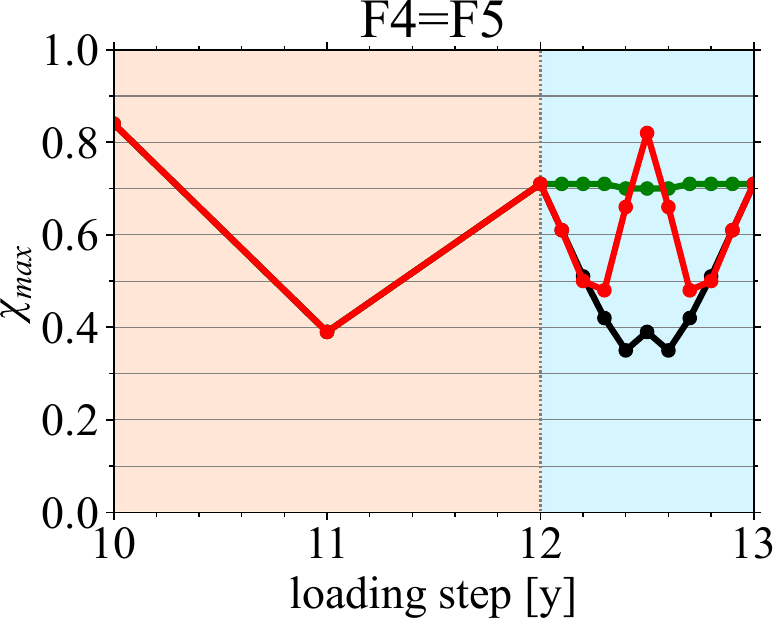}
    \caption{Scenarios \#5: $\chi_{max}$ over time for the CGI and UGS phases. \revB{Scenario \#5 includes two uneven–pressure variants: (\#5a) $\Delta P_1 = -10$ MPa and $\Delta P_2 = 0$ MPa, (\#5b) $\Delta P_1 = -10$ MPa and $\Delta P_2 = -20$ MPa.}}
    \label{fig:CH4-7_chi_zoom}
\end{figure}

\begin{figure}%[htbp]
    \centering
    \includegraphics[width=1.0\textwidth]{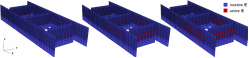}
    \caption{Active interface elements (IE) on the fault system at loading step 12.5 for scenario \#1 (left), \#5a (center) and \#5b (right).}
    \label{fig:7-3D_active}
\end{figure}

The effect of a different pore pressure change in the two reservoir compartments is investigated in Scenarios \#5 (Table \ref{tab:sensitivity}). %The simulations are identical to the reference until the end of the CGI phase; then, during UGS, in block 1 the pore pressure change $\Delta P_1$ is kept equal to reference test case, i.e., $\Delta P_1 = -10$~MPa, whereas in block 2 two different sub-scenarios are analyzed for the end of the production phase: $\Delta P_2 = 0$~MPa for the entire UGS in scenario \#7a and $\Delta P_2 = -20$~MPa in \#7b. Given that the pressure variation behaves as in the reference scenario until the end of CGI, the behavior of the criticality index remains unaffected for all the faults. The same holds for fault F1 during UGS, as the pressure affecting block 1 remains unaltered.
%Conversely,
During UGS, fault F3 and, secondarily, faults F4 and F5,
may exhibit larger $\chi$ values at loading step 12.5 (Figure~\ref{fig:CH4-7_chi_zoom}). 
When $\Delta P_2 = 0$~MPa, faults F4 and F5 keeps the same critical degree experienced at the end of CGI over the UGS phase.
Conversely, with $\Delta P_2 = -20$ MPa the most critical condition is achieved in the middle of the UGS cycle, where $\chi_{max}$ is larger than the threshold value 0.80. %In scenario 7b, UGS production initially unloads F2 respect to the end of CGI and then loads it again with its activation at loading step 12.5 (Figure~\ref{fig:7-3D_active}), when the maximum UGS depletion is reached. The consequent new configuration is then initially unloaded during the first part of UGS injection, but then loaded again when the pressure increase reaches larger values, with a new reactivation at loading steps 12.9 and 13.
Figure~\ref{fig:7-3D_active} shows the distribution of the potentially active elements ($\chi \geq 0.80$) in Scenarios \#5, again within the embedded fault discretization at loading step 12.5.

\begin{figure}
    \centering
    \includegraphics[width=1.0\linewidth]{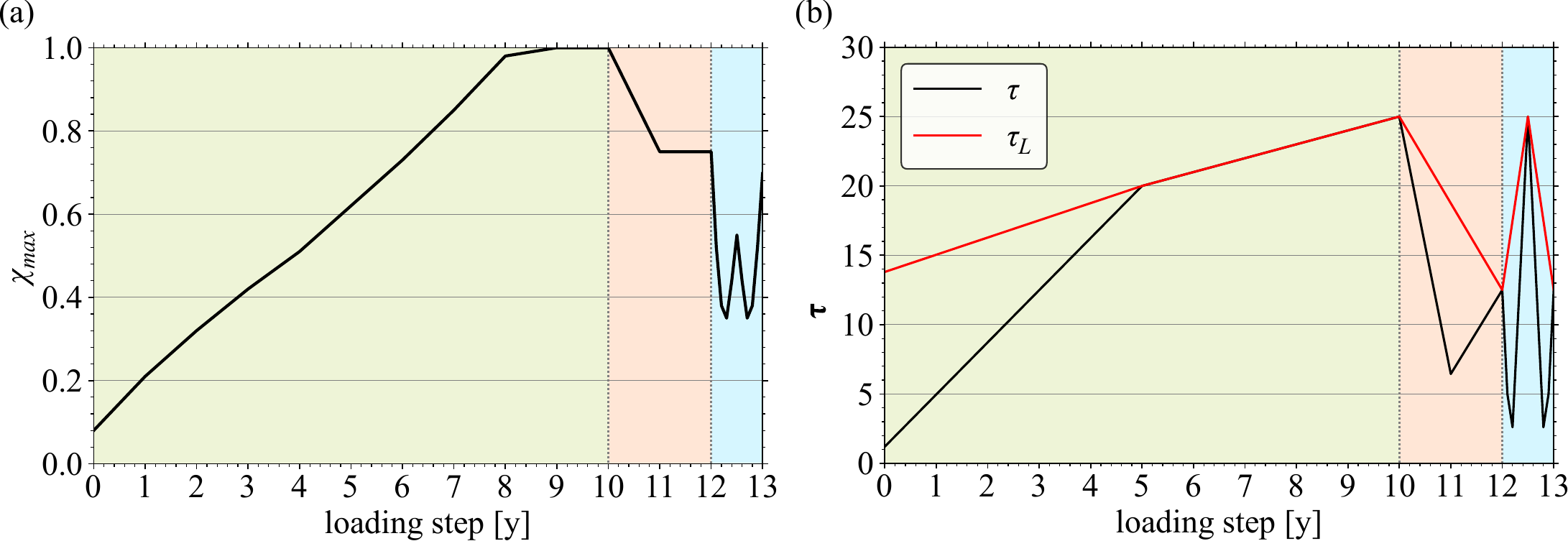}
    \caption{Scenario \#5b: (a) $\chi_{max}$ over time on fault F1 , (b) temporal evolution of the shear stresses $\tau$ and $\tau_L$  on fault F4, showing how 
the shear loading evolves during PP, CGI and UGS.}
    \label{fig:s7-tau}
\end{figure}

%The stress behavior of node 2665 in scenario 7b during PP and CGIis equal to that developed in the reference scenario 1d  but the condition is much more critical then in the reference case around the end of the UGS production (between l.s. \#12.4 - 12.5) and injection (between l.s. \#12.9 - 13) phases when the stress state develops on the yield bound because of the large (20~MPa) pressure changes.
For scenario \#5b, Figure~\ref{fig:s7-tau}~(a) shows the time behavior of $\chi_{max}$ and (b) the shear stress modulus $\tau$ at the reservoir bottom on fault F4 . The shear stress increases during the PP stage, reaching the limit shear stress $\tau_L$ at loading step 5. Then, $\tau$ equates $\tau_L$ at the end of CGI and at the end of UGS production and injection cycle. On fault F2, $\chi_{max}$ approaches 1 at the end of the UGS injection phase too, but remains well below the criticality value at the end of CGI phase. 

\subsubsection{Initial stress regime}
\label{sec:stress}
%On top of the above critical factors stressing the faults, not only during PP but also during CGI and UGS, i.e., when pressure conditions are considered ``normal''  and not typically associated with fault activity, it is worth making a short digression on the potential effects of a variation in the initial stress regime.
The orientation of the principal components of the stress tensor in undisturbed conditions is another key factor for predicting fault stability during reservoir management. This is analyzed in Scenarios \#6. In the reference scenario, the in-situ effective stress field is defined by M$_1$ = $\sigma_h$/$\sigma_v$  = 0.64 and M$_2$ = $\sigma_H$/$\sigma_v$ = 0.83. This configuration implies a stress hierarchy of $\sigma_v$ $>$ $\sigma_H$ $>$ $\sigma_h$, which corresponds to a normal fault regime, consistent with field observations reported in the Tab.~\ref{tab:regime}, with $\sigma_H$ applied orthogonal to the F4-F5 faults and $\sigma_h$ applied orthogonal to the F1-F2 faults.
{To investigate the influence of stress orientation on fault reactivation, we simulate in Scenario \#6a a 90$^{\circ}$ rotation of the principal horizontal stress directions. This test serves two main purposes: (i) assessing the sensitivity of the fault system to changes in the orientation of the stress tensor, and (ii) \revA{accounting for the uncertainty in the relative alignment between the regional stress field and the fault network, so that $\sigma_H$ becomes orthogonal to faults F1–F2 instead.}

%Our observations
The modeling outcome revealed that a rotation of the maximum horizontal stress by 90$^\circ$, corresponding to a fault system rotated in the opposite direction with respect to the regional stress field (i.e., $\sigma_H$ applied orthogonal to F1-F2), does not provide a significant variation with respect to the reference case. %is almost negligible,
By distinction,
an initial stress regime close to normally consolidated conditions (scenario \#6b)
%a reduction of the normal stress on the fault system, as seen in scenario \#4b,
significantly increases the fault reactivation likelihood, %leading to early fault reactivation
as can be seen from the $t_{80}$ values in Figure~\ref{fig:s4-t80}.
Although the occurrence of such a scenario is unlikely, it is important to consider its implications. With the undisturbed stress regime of scenario \#5b, %new initial stress conditions have proven to be significantly more critical compared to the reference condition. When reducing the horizontal principal components to half of the most likely values (i.e., $M_1 = 0.40$ and $M_2 = 0.47$),
the limit $\chi_{max} = 1$ is reached at an earlier point of the PP stage on F1 and F2 faults (loading step 6) and F4 and F5 (loading step 8), resulting in a larger fault area that is prone to potential reactivation (Figure~\ref{fig:s4-t80}). %and significant sliding between the fault surfaces.
The earlier activation during PP causes a reactivation of F1 and F2 faults, %are re-activated also
at the end of the CGI and UGS stages as well. %and at the end of the UGS injection phases.
%The largest fault size close to the reactivation amounts to 175 m and 96 m for faults F1 and F2, and F4 and F5, respectively, over the last period of PP. As the CGI and UGS phases progress, fault areas in a critically stressed state tend to expand further, leading to re-activation of F1 and F2. Finally, we report that F3 is not influenced by any of this variations.

\begin{figure}
    \centering
    \includegraphics[height=0.32\linewidth]{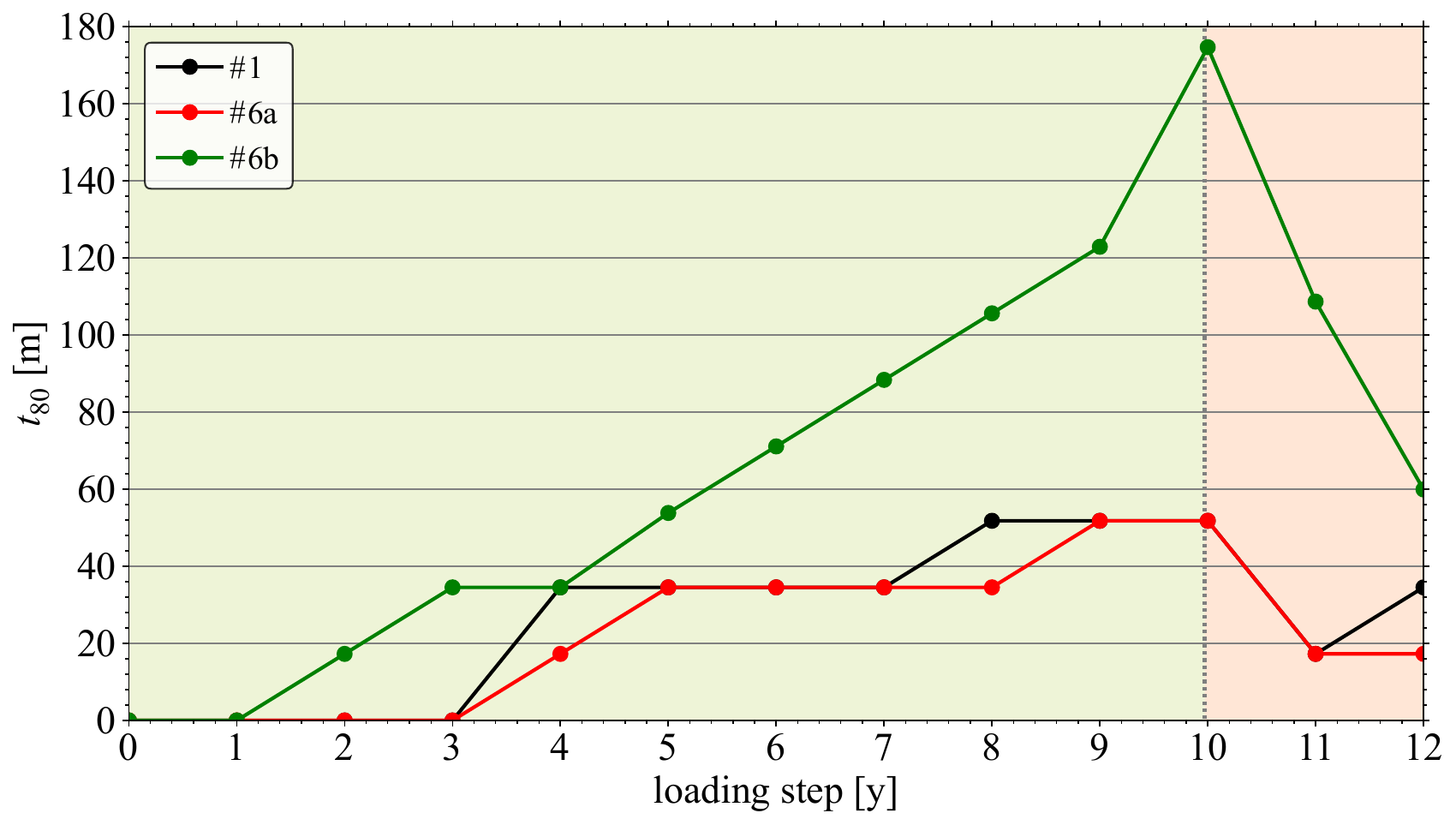}
    \includegraphics[height=0.32\linewidth]{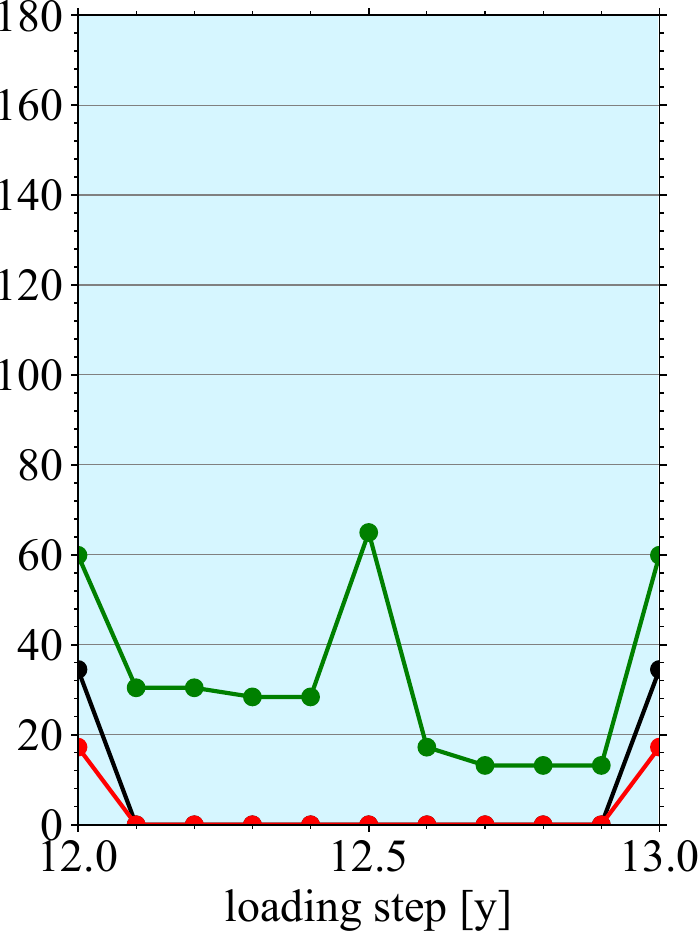}
    \caption{Scenarios \#6: $t_{80}$ over time on fault F2. A zoom for the UGS stage is
      provided on the right panel of the figure. \revB{Scenario \#6 includes two in-situ stress–regime variants: (\#6a) stress rotation, (\#6b) modified horizontal stress ratios ($M_1 = 0.40$, $M_2 = 0.47$).}}
    \label{fig:s4-t80}
\end{figure}

\subsubsection{Ranking factors favoring fault reactivation}
The objective of this section is to define criteria to rank the potential of inducing ``unexpected'' seismic events
for the different investigated scenarios.
%based on various mechanisms, geological settings, geomechanical, and production parameters.
The analysis focuses on fault F3, which is the most stressed in the realistic condition of a compartment offset, and fault F2, which represents the behavior at the reservoir boundaries. 

\revB{The ranking process takes into account the following quantities}
\begin{enumerate}
    \item $\chi_{max}$ during the UGS stage;
    \item the maximum value of average sliding $d_{avg}$;
    \item the loading step of first activation.
\end{enumerate}
\revB{and follows a hierarchical comparison, giving priority to the $\chi_{max}$ indicator and using the others only when needed to discriminate between cases.} The results are presented in Table~\ref{table:ranking}. 

\begin{table}%[htp]
\centering
\small
\begin{tabular}{c|cccc|c|cccc}
\hline
 & \multicolumn{4}{c|}{Fault F2} & & \multicolumn{4}{c}{Fault F3} \\ \hline
 \multirow{2}{*}{Scenario} & $\chi_{max}$ & Max & Act. & $\chi_{max}$ & \multirow{2}{*}{Scenario} & $\chi_{max}$ & Max & Act. & $\chi_{max}$ \\
 & UGS & $d_{avg}$ [m] & year & PP & & UGS & $d_{avg}$ [m] & year & PP \\
 \hline
3c & 1.00 & 0.026 & 7 & 1.00 & 2d & 0.84 & 0.033 & 7 & 1.00 \\
5b & 0.97 & 0.010 & 9 & 1.00 & 2c & 0.81 & 0.007 & - & 0.85 \\
6b & 0.96 & 0.018 & 6 & 1.00 & 5a & 0.53 & 0.000 & - & 0.53 \\
2d & 0.84 & 0.008 & 9 & 1.00 & 5b & 0.53 & 0.000 & - & 0.53 \\
3b & 0.81 & 0.010 & 7 & 1.00 & 2a/b & 0.32 & 0.000 & - & 0.46 \\
%2c & 0.80 & 0.007 & 9 & 1.00 & 1 & 0.00 & 0.000 & - & 0.00 \\
%4a & 0.79 & 0.010 & 8 & 1.00 & 6a & 0.00 & 0.000 & - & 0.00 \\
%3a & 0.78 & 0.008 & 8 & 1.00 & 6b & 0.00 & 0.000 & - & 0.00 \\
%1 & 0.78 & 0.007 & 8 &  1.00 & 3a & 0.00 & 0.000 & - & 0.00 \\
%5a & 0.78 & 0.007 & 9 & 1.00 & 3b & 0.00 & 0.000 & - & 0.00 \\
%2a/b & 0.77 & 0.007 & 8 & 1.00 & 3c & 0.00 & 0.000 & - & 0.00 \\
%6a & 0.74 & 0.007 & 10 & 1.00 & 3d & 0.00 & 0.000 & - & 0.00 \\
%3d & 0.70 & 0.009 & 8 & 1.00 & 4a & 0.00 & 0.000 & - & 0.00 \\
%4b & 0.78 & 0.005 & - & 0.78 & 4b & 0.00 & 0.000 & - & 0.00 \\
2c & 0.80 & 0.007 & 9 & 1.00 & 1 & \multicolumn{4}{c}{\textit{no activation due to symmetry}} \\
4a & 0.79 & 0.010 & 8 & 1.00 & 6a & \multicolumn{4}{c}{\textit{no activation due to symmetry}} \\
3a & 0.78 & 0.008 & 8 & 1.00 & 6b & \multicolumn{4}{c}{\textit{no activation due to symmetry}} \\
1 & 0.78 & 0.007 & 8 &  1.00 & 3a & \multicolumn{4}{c}{\textit{no activation due to symmetry}} \\
5a & 0.78 & 0.007 & 9 & 1.00 & 3b & \multicolumn{4}{c}{\textit{no activation due to symmetry}} \\
2a/b & 0.77 & 0.007 & 8 & 1.00 & 3c & \multicolumn{4}{c}{\textit{no activation due to symmetry}} \\
6a & 0.74 & 0.007 & 10 & 1.00 & 3d & \multicolumn{4}{c}{\textit{no activation due to symmetry}} \\\
3d & 0.70 & 0.009 & 8 & 1.00 & 4a & \multicolumn{4}{c}{\textit{no activation due to symmetry}} \\
4b & 0.78 & 0.005 & - & 0.78 & 4b & \multicolumn{4}{c}{\textit{no activation due to symmetry}} \\
\end{tabular}
\caption{
%Fault F3 simulated scenarios are ranked based on their potential to induce fault sliding during CGI and UGS phases. The maximum value of the criticality index ($\chi_{max}$) is considered over the entire reservoir life, including the PP period. The ranking employs bold-italic, bold, and italic fonts to denote a critical, almost critical, and safe condition, respectively.}
Ranking of the simulated scenarios according to the largest $\chi_{max}$ in the UGS stage, the largest sliding, and the earliest first activation step.}
\label{table:ranking}
\end{table}

As expected, fault activation during CGI and UGS cycles exhibits distinct critical factors for boundary and central faults. For fault F2, stability is mainly influenced by the initial stress regime of the system, geomechanical properties, such as reservoir stiffness, and fault characteristics, such as cohesion, static friction angle, and the presence of slip weakening.
In contrast, for fault F3, geometric parameters %characterizing the fault/reservoir system
play a major role. %in this case.
The stability of fault F3 is strongly influenced by a compartment offset, non-vertical fault planes, and differential pressure changes in the two compartments.

% Commented out by Andrea
%In our simulations, the presence of a viscous caprock (scenario \#9) %yields unexpected results compared to previous modeling studies (e.g., \cite{wassing_2017}). Tables 2 and 3 reveal that the existence of a salt viscous formation sealing the top of the reservoir is not
%does not appear to be
%a crucial factor in favoring fault reactivation. This discrepancy with previous findings, e.g., \cite{wassing_2017}, %which focused on production reservoirs,
%can be attributed to the short-term (seasonal) pressure fluctuations in a UGS reservoir, which limit long-term viscous effect.
%%It is worth mentioning that the mesh size used in the 3D modeling study (20 m along the vertical direction) is notably larger (one order of magnitude) than the meshes employed in 2D investigations specifically dedicated to studying caprock behavior. This relatively large mesh size could have smoothed out the effects of viscosity.
%
%This ranking may be a tool to prioritize attention on most impacting parameters, ensuring they are accurately considered to avoid potential hazards. %The proposed approach of ranking scenarios can be utilized in other studies with similar configurations to assess the degree of criticality and identify the most comparable scenarios for further analysis.
%It is worth noting that while this section only presents the ranking for methane, the same general order of criticality applies to other fluids as they exhibit similar behaviors.

\subsection{Sensitivity analysis: Stage 2}
%The upcoming step in this investigation involves gaining valuable insight on how the storage of different gases for various purpose impacts fault instability. To achieve this,
In the second stage of the sensitivity analysis,
we combine the most influential factors to produce a potential fault reactivation, as identified in Stage 1, and test the new scenarios for the storage of different gases, namely CH$_4$, CO$_2$, H$_2$, and N$_2$.
The scenarios discussed here are summarized in Table \ref{tab:combinations}.
We recall that the primary objective of this analysis is to explain realistic
configurations likely to be encountered in the Rotliegend formation, rather than to
explore ``extreme'' conditions  by combining the most unfavorable parameter values identified in Stage 1. Therefore, taking into account the most influential
factors from Stage 1, we chose to focus on geologically plausible combinations  presented in Table
\ref{tab:combinations}. They combine a compartment offset, non-vertical fault planes, variations in the fault characteristics, i.e., cohesion and static friction angle, and reservoir stiffness. Two dedicated combinations analyze the effect of different geomechanical properties,
i.e., decreasing (W) and increasing (H) of  reservoir stiffness most likely associated to chemo-mechanical effects after the PP phase.

%anticipated by Pietro (lines with 5 %%%%%)
%%%%%\begin{table}
%%%%%\centering
%%%%%\begin{tabular}{cl}
%%%%%\hline
%%%%%Combination & Parameter/mechanism \\
%%%%%\hline
%%%%%C1 & $c = 0$~MPa, %F3 fault dip angle
%%%%%$\delta=+65^\circ$, %and compartment offset of 100~m
%%%%%$o=100$ m\\
%%%%%C2 & $\varphi_s = 20^\circ$, %F3 fault dip angle
%%%%%$\delta=+65^\circ$, %viscous caprock, %and compartment\\
   %& offset of 100~m \\
%%%%%   $o=100$ m\\
%C3 & $c = 0$~MPa, %F3 fault dip angle
%$\delta=+65^\circ$, %viscous caprock, %and compartment\\
%%   & offset of 100~m \\
%$o=100$ m\\
%C4 & $c = 0$~MPa, %F3 fault dip angle
%$\delta=+65^\circ$, $o=100$ m, %compartment offset of 100~m and\\
%%   & reservoir Young modulus of
%$E=20$~GPa\\
%C5 & $\varphi_s = 20^\circ$, %F3 fault dip angle
%$\delta=+65^\circ$, %viscous caprock, %compartment\\
%%   & offset of 100~m and reservoir Young modulus of
%$o=100$ m, $E=20$~GPa\\
%%%%%W  &  $E$ decreased by 30\% during injection \\
%%%%%H  &  $E$ increased by 30\% during injection \\
%%%%%\hline
%%%%%\end{tabular}
%%%%%\caption{Combination of settings addressed in Stage 2 of the sensitivity analysis for the fluids under consideration.}
%%%%%\label{tab:combinations}
%%%%%\end{table}

\begin{table}
\centering
\small
\begin{tabular}{cl}
\hline
Combination & Parameter/mechanism \\
\hline
C1 & $c = 0$~MPa, %F3 fault dip angle
$\delta=+65^\circ$, %and compartment offset of 100~m
$o=100$ m\\
C2 & $\varphi_s = 20^\circ$, %F3 fault dip angle
$\delta=+65^\circ$, %viscous caprock, %and compartment\\
   %& offset of 100~m \\
   $o=100$ m\\
%C3 & $c = 0$~MPa, %F3 fault dip angle
%$\delta=+65^\circ$, %viscous caprock, %and compartment\\
%   & offset of 100~m \\
%$o=100$ m\\
C3 & $c = 0$~MPa, %F3 fault dip angle
$\delta=+65^\circ$, $o=100$ m, %compartment offset of 100~m and\\
%   & reservoir Young modulus of
$E=20$~GPa\\
C4 & $\varphi_s = 20^\circ$, %F3 fault dip angle
$\delta=+65^\circ$, %viscous caprock, %compartment\\
%   & offset of 100~m and reservoir Young modulus of
$o=100$ m, $E=20$~GPa\\
W  &  $E$ decreased by 30\% after PP \\
H  &  $E$ increased by 30\% after PP \\
\hline
\end{tabular}
\caption{Combination of settings addressed in Stage 2 of the sensitivity analysis for the fluids under consideration.}
\label{tab:combinations}
\end{table}

%Several configurations (see Table~\ref{tab:combinations}) are investigated to determine which combination poses the greatest risk of fault failure, regardless of aseismic slip. Reduced Coulomb parameters, system geometry, stiffness variations, and the presence of a viscous caprock are among the parameters considered.
The potential hazard of fault reactivation is evaluated for each combination by the parameter $t_{80}$, whose behavior in time is shown in Figure~\ref{fig:fluids_t80}.
Regardless of the fluid, Figure \ref{fig:fluids_t80} shows that the combination C2 provides the conditions yielding the larger $t_{80}$ values. 
Therefore, all subsequent considerations will focus on this specific scenario.
Scenario C2 is characterized by a compartment offset and a reduced friction angle. Conditions close to reactivation occur during both CGI and UGS/UHS, or long-term storage, stages. Notice that $t_{80}$ after the end of PP is larger if the faults already experienced a significant reactivation (i.e., with a large $t_{80}$) during PP. %; ii) the reservoir pressure at the end of injection closely reaches the initial value $P_i$ (Figure~\ref{fig:sketch_p}); iii)
%the reservoir pressure remain close to the $P_i$ throughout the UHS cycles.
%If a fault is not reactivated during the PP, $t_{80}$ is kept at zero in the ST and CGI stages and in the UHS cycles.
%
As compared to other critical combined configurations, scenario C2 is more critical primarily because of the reduction of the friction angle. %and the inclusion of viscous effects. However, the overall influence of viscosity on fault instability is not found to be significant. Instead, the decrease in friction angle plays a major role in comparison to the decrease in cohesion.
%This observation is further supported by analyzing scenario C3, where the only difference with C2 is the reduction in cohesion instead of the friction angle.
%In cases where the reservoir stiffness increases, the criticality %(as measured by $t_{80}$)
%generally decreases (see scenario C4). %This finding is validated by combination C4, which differs from C1 only in the increased reservoir stiffness.
%Since a larger reservoir stiffness %enhances fault stability by
%reduces the reservoir deformation %compaction and contraction of the reservoir.
%and
%As a result,
%the critical condition for $\chi_{max}$ is never reached (see scenario \#4b), this combination is not considered. %leading to smaller fault displacements in this configuration.
%The same trend is observed in combination C5, which shares the same reduction in friction angle %and viscosity
%as C2, but has an increased reservoir stiffness. %The higher reservoir stiffness further improves fault stability, particularly on faults F4 and F5. Consistent outcomes have been achieved regardless of the fluid type (Figure~\ref{fig:fluids_t80}).

\begin{figure}
    \centering
    \includegraphics[width=1\linewidth]{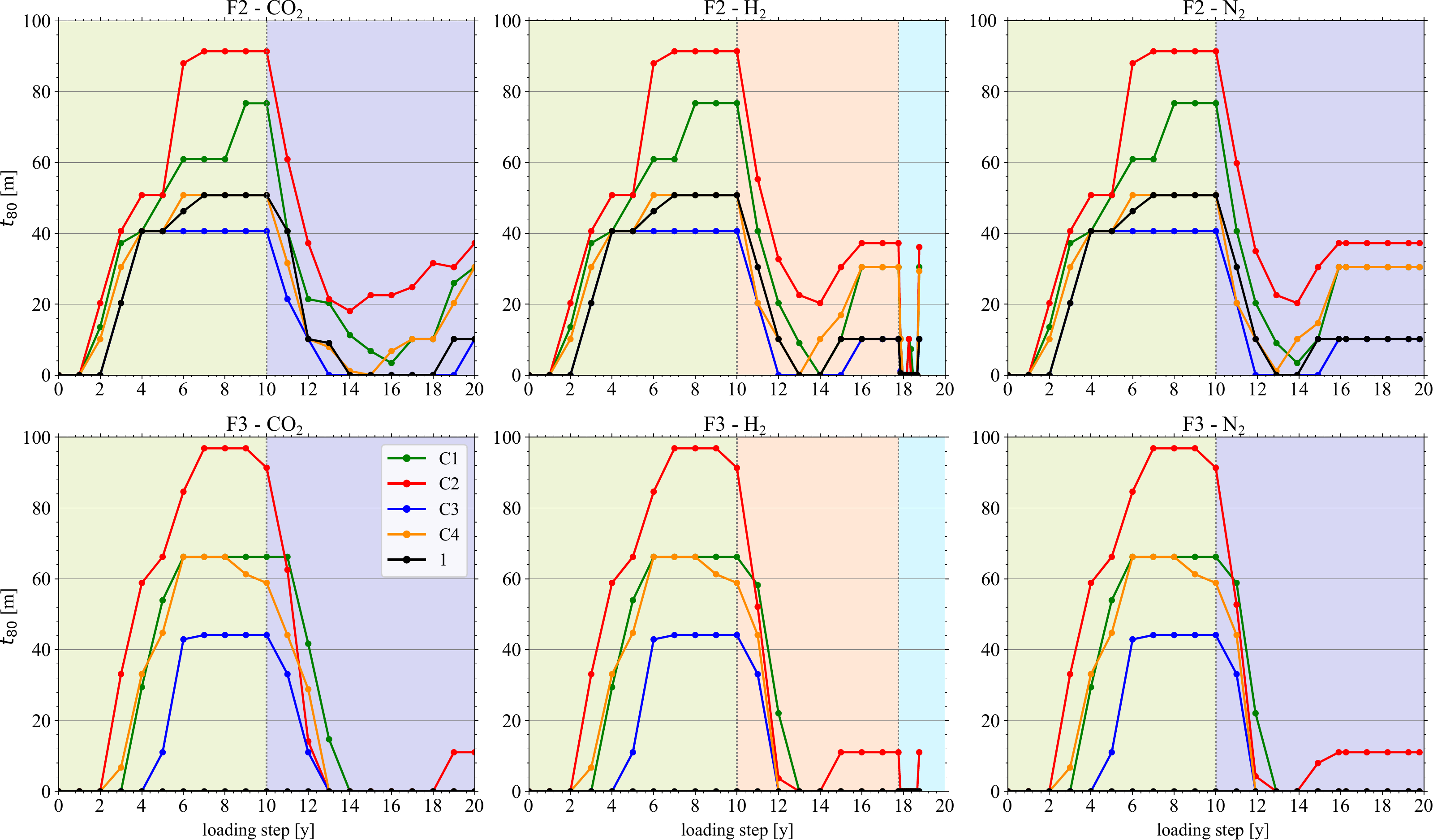}
    \caption{Sensitivity analysis, Stage 2: $t_{80}$ over time for fault F2 and F3 and the different fluids under investigation. \revA{The legend refers to the combinations of settings addressed in Stage 2 of the sensitivity analysis and reported in Table~\ref{tab:combinations }}}
    \label{fig:fluids_t80}
\end{figure}

Due to its significance,
the outcomes of the critical combination C2 have been separately described for CH$_4$ and the other gases. This distinction is made based on the differences in pressure history (see Figure \ref{fig:pmean}). %It is important to remind that a two-year CGI has been used for methane, while for other gases, it spans about 8-10 years, depending on the type of fluid, to ensure it is as physically feasible as possible.

\subsubsection{UGS critical scenario}
\revB{The combination of factors in configuration C2 causes a potential reactivation, mainly during PP. During UGS, faults F2 and F3 reach near-critical values of $\chi_{max}$ at the end of both the withdrawal and storage phases.} Figure~\ref{fig:CH4-C2-chi}, which shows $\chi_{max}$ in time for scenario C2 as compared to the reference case (scenario \#1), highlights that critical conditions for fault reactivation develop earlier during PP. $\chi_{max}$ is also greater than, or close to, 0.8 during CGI and UGS. %, with $\chi$ reaching 0.79 for fault F2-F4-F5 at l.s. 13 and 0.73 for fault F1. Faults F1, F2, and F3 slightly slide also at the end of the UGS.
Because of the loss of symmetry, F2 is more stressed than F1 and  $\chi_{max}$ achieves values larger than 0.8 on fault F3 as well.
%increases also on fault F3, with values up to 0.8 at the end of PP, and close to the threshold during CGI and UGS injection phases.
Figure~\ref{fig:CH4-C2-3D} provides a 3D view of $\chi$ distribution in space at a few significant loading steps while Figure~\ref{fig:CH4-C2-slide} shows the time behaviour of the maximum sliding $d_max$.

\begin{figure}%[htbp]
     \centering
     \includegraphics[width=0.41\linewidth]{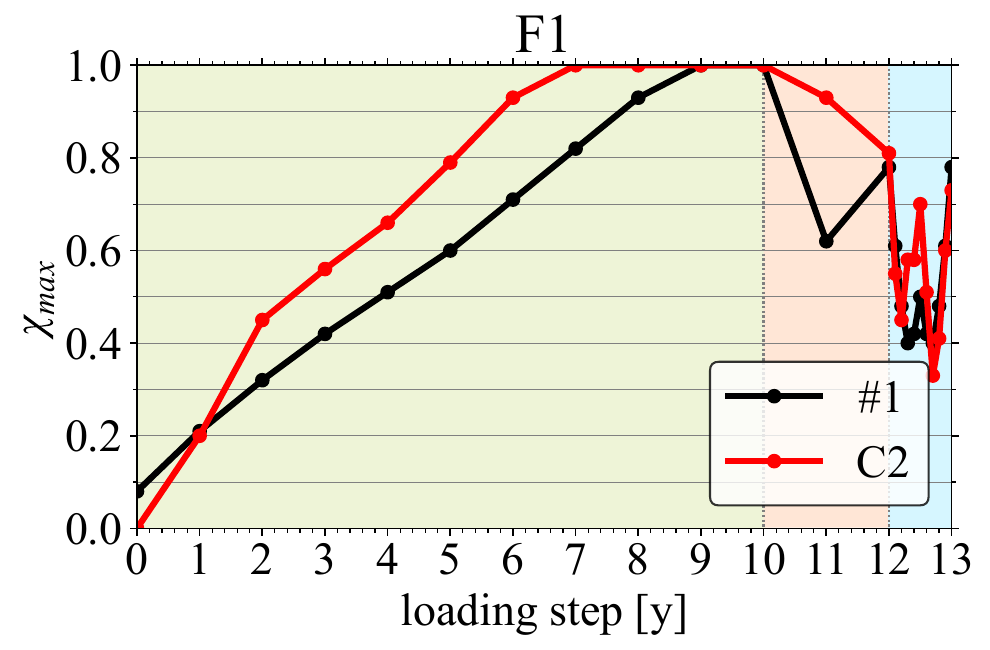}
     \includegraphics[width=0.4\linewidth]{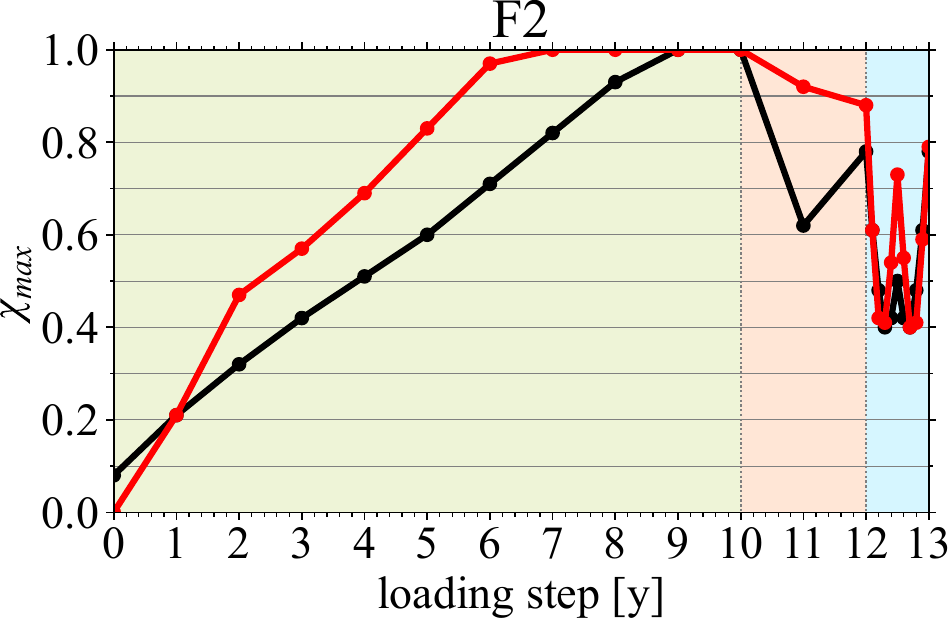}
     \includegraphics[width=0.4\linewidth]{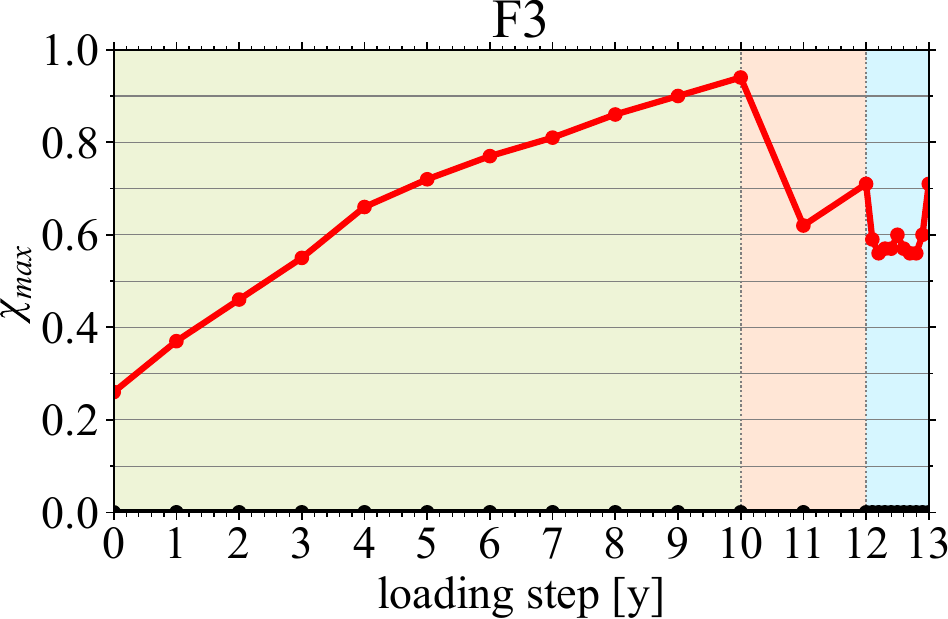}
     \includegraphics[width=0.4\linewidth]{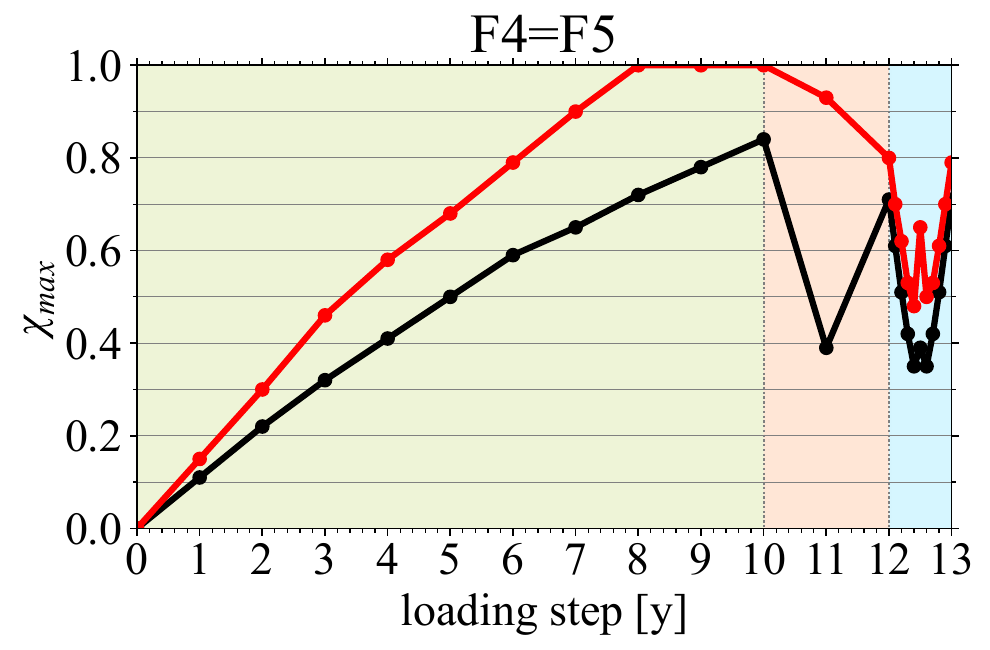}
     \caption{Scenario C2, UGS: %(combination C2):
     $\chi_{max}$ over time for all faults compared to the reference case (scenario \#1). \revB{ Scenario C2 corresponds to the low-friction configuration 
with $\varphi_s = 20^\circ$, fault dip $\delta = +65^\circ$, and Block~2 offset of 100~m.}} 
     \label{fig:CH4-C2-chi}
 \end{figure}

\begin{figure}%[htbp]
     \centering
     \includegraphics[width=1.0\linewidth]{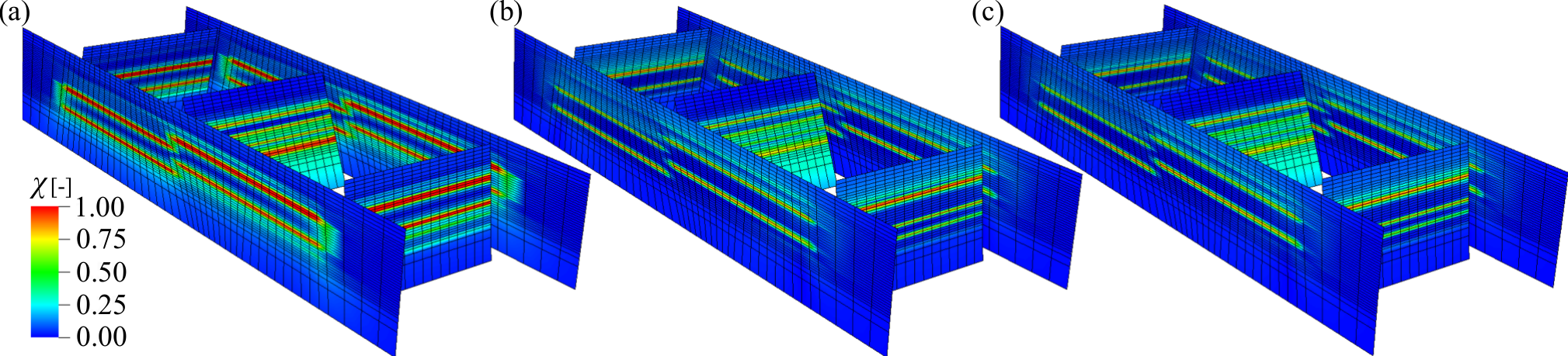}
     \caption{Scenario C2: %(combination C2):
     $\chi$ distribution in space on the fault system at the end of PP (a), CGI (b) and UGS (c).} 
     \label{fig:CH4-C2-3D}
 \end{figure}

The behavior of $d_{max}$ for all the faults is shown in Figure~\ref{fig:CH4-C2-slide}. During the CGI phase, all the faults slide, primarily during the period between loading step 11 and 12, in the opposite direction to what occurred during PP. Additionally, faults F1, F2, and F3 slightly slide again at the end of the UGS injection phase.

\begin{figure}%[htbp]
    \centering
    \includegraphics[width=0.45\linewidth]{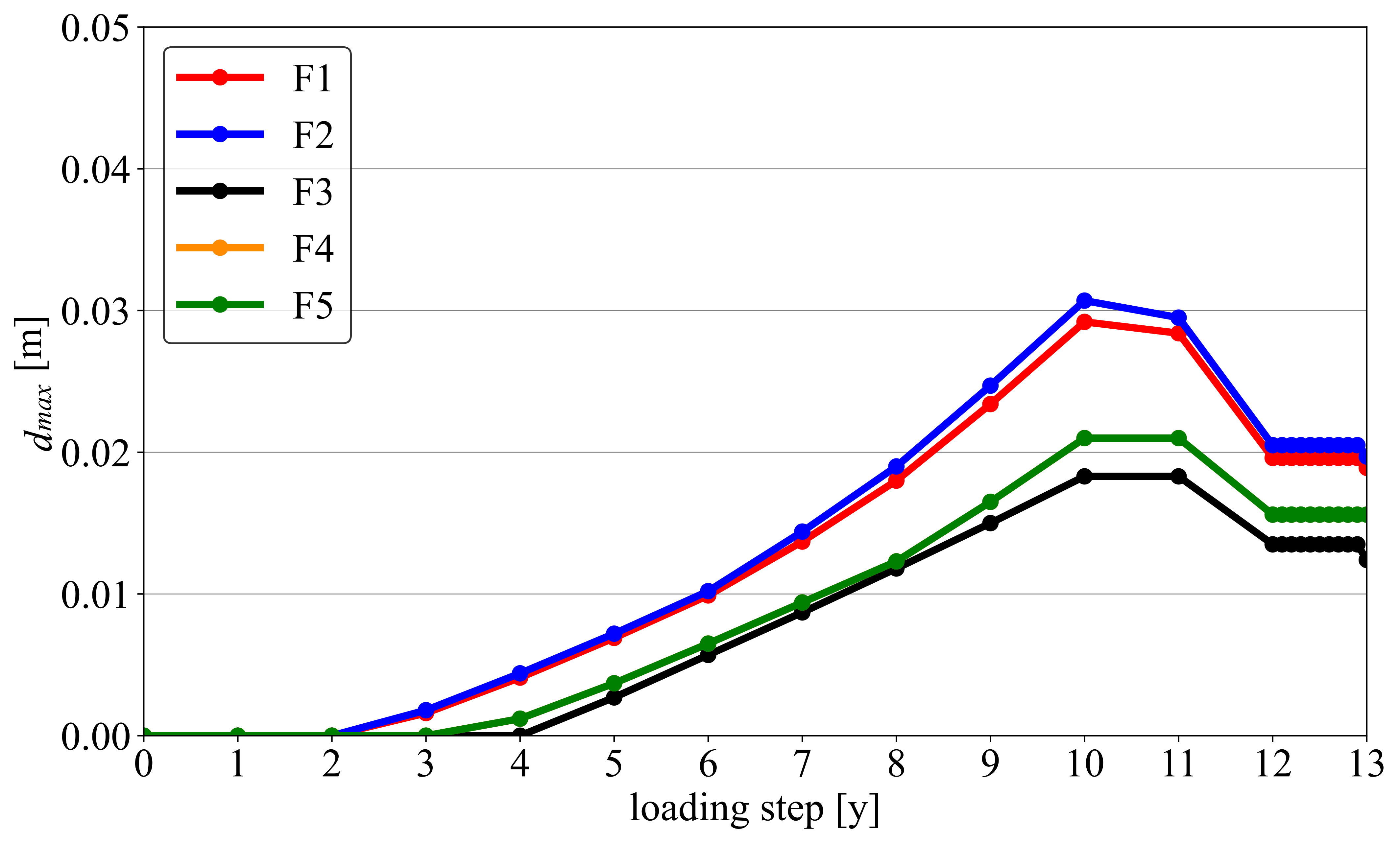}
    \caption{Scenario C2, UGS: maximum sliding $d_{max}$ over time for all the faults.}
    \label{fig:CH4-C2-slide}
\end{figure}

\subsubsection{CO$_2$, H$_2$, N$_2$ storage dynamics and fault activation}
The goal of this section is to compare how the temporary or permanent storage of different gases may influence the fault stability, as measured by $\chi_{max}$ in the critical scenario C2. Since the pore pressure variation during PP is the same for each fluid, $\chi$ has the same behavior until loading step 10. %, approaching a critical condition of  $\chi_{max} = 1$ already at loading step 4 for faults F1, F2, F4 and F5.
%During this stage, the reservoir compaction increases the shear stress along the fault, while its contraction unloads the normal stress. Instead, fault F3 experiences an increase of $\chi_{max}$ due to the increase in shear stresses, particularly at the end of this phase.
After that, %since different gases have different pressure histories,
the impact on fault stability must be considered based on whether temporary or permanent storage scenarios occur.

For the permanent storage scenarios, i.e., CO$_2$ and N$_2$, %the faults are stabilized
at the beginning of the storage $\chi_{max}$ decreases due to the unloading of the shear stress %along the fault
caused by reservoir expansion. However, as the storage continues, the shear stress is loaded with the opposite sign, reaching a relative maximum at the end of the injection stage. This point corresponds to the full pore pressure recovery %return of the pore-pressure in the reservoir to its initial value, before the beginning of the PP
(Figure~\ref{fig:C2_gases}).
%The effect of the viscous caprock %over- and side-burden
%is slightly appreciable only at the end of gas injection but overall it does not appear to have a major impact on the fault stability. Note that the oscillation in $\chi_{max}$ observed on F1 and F2 %. This peculiarity
%is due to a change in the elements experiencing the most critical condition.
%The same considerations hold for the CGI phase with H$_2$ (l.s. 10 to 18).
As far as it concerns the injection/withdrawal cycles for H$_2$, $\chi_{max}$ on F1, F2 and F4-F5 reaches its maximum at the end of the injection stage, corresponding to the highest value of pressure in the reservoir. For F3, $\chi_{max}$ exhibits an additional peak at the end of the UHS withdrawal stage. %, when the pore pressure reaches the lowest value, due to the unloading of the normal stress on the fault by the contraction of the two reservoir blocks.

\begin{figure}%[ht]
    \centering
    \includegraphics[width=0.9\linewidth]{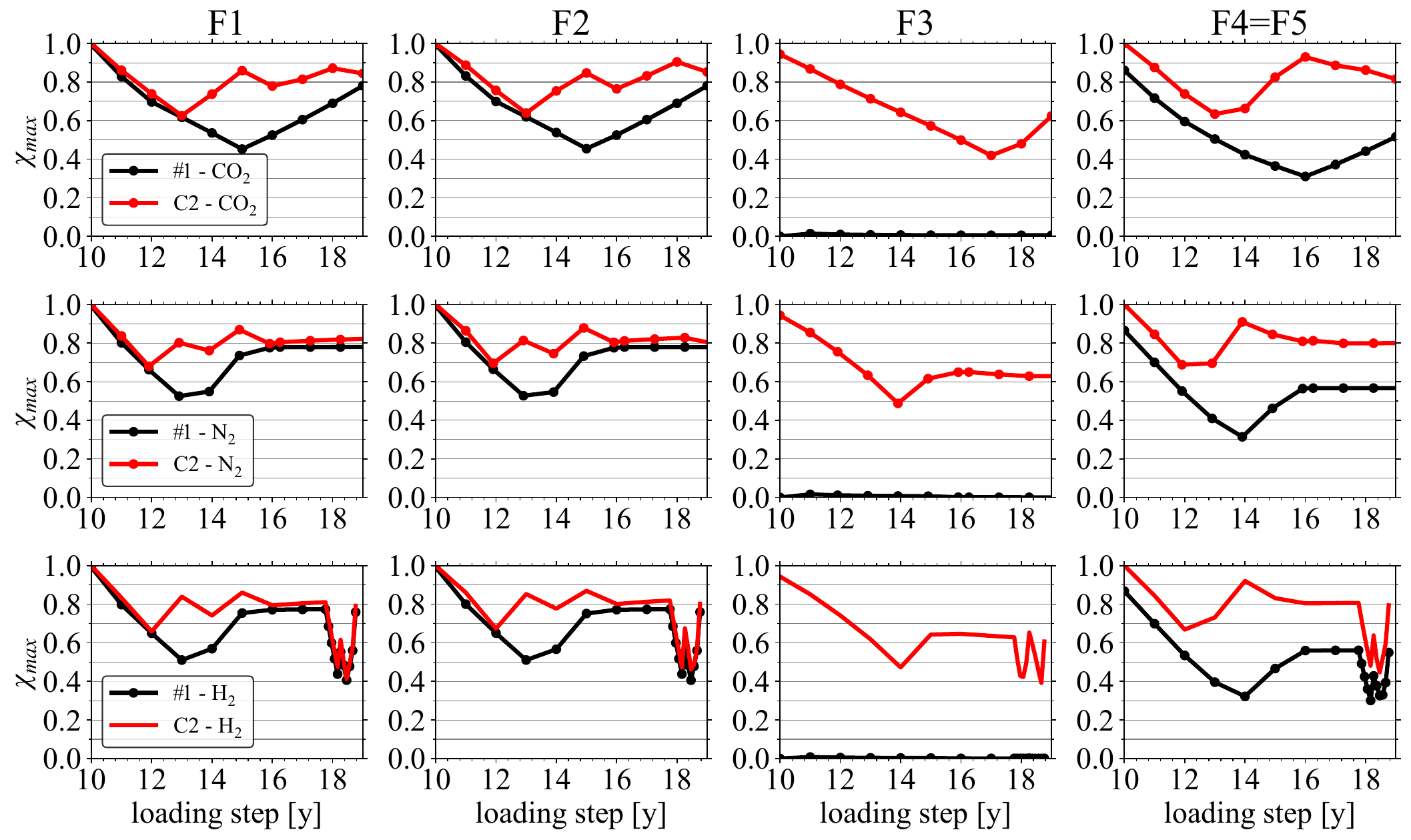}
    \caption{Scenario C2: $\chi_{max}$ over time for CO$_2$, N$_2$ and H$_2$ after the end of PP.} %with respect to loading steps for each fault from l.s. 10 (end of PP) to l.s. 19 (end of ST/UHS). The behavior during PP (l.s. 1 to 10) is exactly the same for different fluids since the pressure history is unaltered. Due to symmetry, F4 and F5 behave identically.}
    \label{fig:C2_gases}
\end{figure}

% \iffalse
% \begin{figure}
%     \centering
%     \includegraphics[width=1.00\linewidth]{figs/plot/chi_CO2-N2}
%     \caption{Scenario 11 - CO$_2$-N$_2$: $\chi_{max}$ vs loading steps for each fault. Note that due to symmetry F4 and F5 behave identically}
%     \label{fig:C2_CO$_2$}
% \end{figure}
% \begin{figure}
%     \centering
%     \includegraphics[width=1.00\linewidth]{figs/plot/chi_H2}
%     \caption{Scenario 11 - H$_2$: $\chi_{max}$ vs loading steps for each fault. Note that due to symmetry F4 and F5 behave identically}
%     \label{fig:C2_H$_2$}
% \end{figure}
% \fi

\subsubsection{Geochemical effects on fault stability}
%sout{Geochemical effects are also considered. Based on the literature review summarized in Section \ref{sec:fluids}, we have incorporated the possibility to account for geochemical influences by introducing a hardening/weakening behavior of reservoir rocks.}
Due to the uncertainty observed in the literature review (Section \ref{sec:fluids}), geochemical effects are considered in a simplified, indirect manner. We account for potential fluid-dependent influences on reservoir rock stiffness by modifying the Young modulus.
Specifically, after PP the Young modulus is increased or decreased by 30\%, depending on the fluid involved. For CO$_2$, both weakening and hardening scenarios have been observed, depending on the actual rock mineralogical composition. Conversely, H$_2$ might induce a weakening effect, while N$_2$ shows a relatively neutral geochemical impact on the reservoir rock.

Figure~\ref{fig:W-H} shows the behavior of $\chi_{max}$ for the weakening/hardening scenarios (W/H). %on fault F2.
\revB{Fault F3 remains stable throughout the simulation due to its vertical dip and the absence of offset, which limit the development of shear stress on the fault plane.}
During CO$_2$ injection, with a hardening behavior, $\chi_{max}$ first decreases and then increases with respect to the reference case. The opposite occurs for the weakening scenario. %shows an increase, particularly during the initial half of the storage stage. This implies that the fault becomes more resistant to deformation and less prone to slip or movement. On the contrary, when accounting for a weakening behavior, fault stability decreases, especially in the first half of the storage stage.
Similar considerations hold for H$_2$ in the weakening scenario, with an even smaller difference with respect to the reference case.
%Regarding UHS, the overall stability of the faults affected by weakening behavior is slightly reduced during the first half of the CGI stage. Notably, during the UHS cycles, increased instability is observed primarily toward the end of the UHS-P stage, coinciding with the lowest pressure values in the reservoir.
Based on these results, we can preliminarily conclude that geochemical effects on the reservoir rock have a minor impact on fault stability. %weakening or hardening, resulting in a small impact on the stability of the fault, as indicated by a slight increase in the value of $t_{80}$.

\begin{figure}%[htbp]
    \centering
    \includegraphics[width=0.9\textwidth]{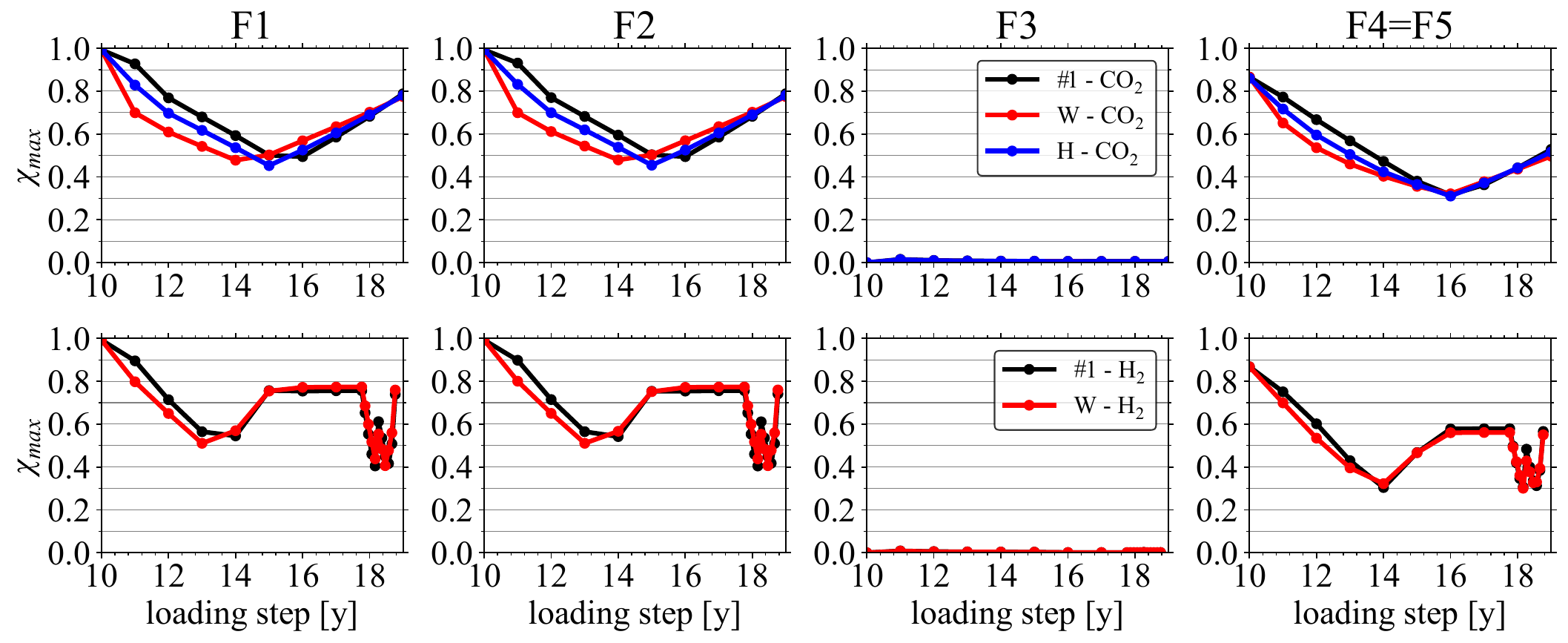}
    \caption{Scenarios W and H: $\chi_{max}$ over time for all faults after PP, shown for CO$_2$ (upper row) and H$_2$ (lower row). \revB{Scenario W corresponds to reservoir weakening 
($E$ decreased by 30\% after PP), whereas Scenario H represents reservoir hardening 
($E$ increased by 30\% after PP). %$\chi_{max}$ on Fault F3 is null due to symmetry.
} }
    \label{fig:W-H}
\end{figure}

\section{Discussion}
% Add brief summary of results (Rev1)

\subsection{Comparison with previous studies}

A few papers addressing the topics of our work have recently been published.
Works by \cite{wassing_2017} and \cite{Hau_etal18} are of particular interest because they focus on Rotliegend reservoirs in the Netherlands and northern Germany, providing insights into fault reactivation and fault rupture during PP.
Although the authors: i) do not consider induced seismicity during UGS stages, ii) use a simplified 2D setting, and iii) employ a different modeling approach concerning fault activation, these studies carry out a parametric analysis on the same geometric/geomechanical features investigated in this paper.

The findings of our work mostly align with the outcomes by %Wassing et al. (2017)
\cite{wassing_2017} regarding fault failure mechanisms and slip initiation. Our study confirms that fault slip initiates at the top of the reservoir rocks where the initial \textit{shear capacity utilization} SCU~\citep{wassing_2017}, corresponding to our safety factor $\chi$, reaches the maximum values. Both studies agree that faults without offset experience reactivation at later stages of depletion, and fault rupture does not extend upwards into the caprock. Scenarios \#2c, \#2d, and C2, for example, confirm that in the case of a compartment offset, an earlier fault activation is favored and remains bounded within the reservoir depth range. Regarding the influence of the contrast in elastic properties distribution, our study agrees with %Wassing et al. (2017)
\cite{wassing_2017} only partially. They observed a secondary peak in shear stress at the bottom of the footwall reservoir block and a localized decrease in shear stresses and SCU at the bottom of the hanging wall reservoir block. Qualitatively, the results presented above (see, for example, Figure~\ref{fig:C2_gases}) agree. However, it must be considered that \cite{wassing_2017} assigned an uniform density and elastic properties to all rocks while our study, in particular scenarios \#4a and \#4b, demonstrates that the contrast in elastic properties may significantly impact on stress concentration and, therefore, fault failure (Figure~\ref{fig:CH4-6_chi_zoom}).

% Commented out by Andrea
%%Differences in modeling approaches should be considered when comparing the results.
%It is also important to acknowledge some differences in the outcomes, which can be attributed to disparities in modeling approaches, initialization, and timing. While \cite{wassing_2017} suggested that the viscoelastic creep behavior of the caprock promotes fault reactivation and early fault rupture, our results indicate a lower degree of fault reactivation facilitated by viscous creep.
%%These variations can be attributed to differences in initialization and timing between the two modeling approaches.
%Instead, our study suggests that the viscoelastic caprock promotes a more critical fault condition during later underground gas storage cycles.

Comparing our results with the findings presented by %Haug et al. (2018)
\cite{Hau_etal18} regarding factors influencing fault reactivation and criticality,
we agree that the initial stress regime (scenarios \#6a and \#6b, Figure ~\ref{fig:s4-t80}) plays significant roles in fault behavior. Obviously, this relation is a well knows in the literature \citep{Fou_etal18,Mun_etal15,walsh_2016}. Indeed, fault criticality largely increases as the horizontal components of the natural stress regime decrease relative to the vertical stress. Furthermore, our outcomes confirm that the stiffness contrast between the reservoir and surrounding rocks (scenarios \#4a and \#4b) governs the stress redistribution and the degree of stress rotation during the reservoir development, impacting fault reactivation.
According to %Haug et al. (2018)
\cite{Hau_etal18} the depletion of thicker reservoir horizons results in a stronger fault-loading compared to the depletion of thinner reservoir horizons, as a thick reservoir undergoes a relatively larger strain for a same change of pore fluid pressure. Although our sensitivity analysis did not directly consider reservoir thickness as a parameter, we observed consistent outcomes when a larger strain is attributed to a larger pressure decrease (scenario \#5b, Figure~\ref{fig:s7-tau}).

However, there are a few aspects where our project disagrees with findings from %Haug et al. (2018)
\cite{Hau_etal18}. While they suggested that the fault throw should be half of the reservoir thickness to obtain the maximum shear stress ratio values, we found that the most critical condition occurred when the fault offset is equal to the entire reservoir thickness. This discrepancy can be attributed to the different  modeling setup, including the horizontal-to-vertical ratio of the natural stress components and the reservoir depth.
%
%The project largely supports the findings of \cite{Hau_etal18}, but there are some areas of disagreement.
These disparities emphasize the importance of the modeling approach and setup parameters, highlighting the need for further research to fully understand fault behavior in different geological contexts.

% To add:
% - Effect of elasto-plasticity or visco-plasticity of the reservoir on the stress path. See e.g. the articles of Pijnenburg et al. 2019, 2020 on Rotliegend sandstones.
%- Horizontal stress variations across lithologies and implications for the initial stress on bounding faults which are typically juxtaposed against clays of salt.
%- Refer to works of van den Bogert, Buijze, and others for exhaustive analysis of the effect of offset, dip, and other parameters.
%- Discuss the mechanisms and possible alternatives. Thermal effects should be discussed, as well as locally elevated pressures during injection, preferential fluid flow pathways, etc. 

\subsection{Definition of safe operational bandwidths}
After the end of PP, when the natural fluid pressure $P_i$ is reduced to $P_{min,PP}$, the reservoir experiences a relatively fast pressure recovery to $P_{max,CGI/CCS}$ during CGI or gas storage (CO$_2$ or N$_2$). Afterward, UGS or UHS are characterized by a cyclic pressure fluctuation between $P_{min,UGS/UHG}$ and $P_{max,UGS/UHS}$ at the end of the production and injection phases, respectively. In the usual practice and in alignment with the legislation of some countries, e.g., the Netherlands \citep{minEZK_grijpskerk2022,tno_ondergrondse2018}, $P_{max,CGI/CCS} \simeq P_{max,UGS/UHS} \simeq P_i$ and $P_{min,UGS/UHS} \geq P_{min,PP}$ (Figure~\ref{fig:P_bandwidth}-left). In this framework, guidelines for the definition of ``safe operation bandwidths'' in gas storage, i.e., operations with a reduced risk of fault reactivation, must identify proper values for the aforementioned pressure bounds. Nevertheless, it has to be recalled that a fault reactivation could always occur aseismically \citep{cappa_2019_stabilization}.

The interpretation of the modeling results allows outlining some key guidelines.
\revA{When a reservoir is depleted during PP, the
faults adjust to a new, permanently deformed configuration. If the same faults slipped
during PP, they become more sensitive to later pressure changes, because the stress
state no longer returns to its original configuration when pressure is increased again.
This means that even pressure variations that remain within the range already
experienced during PP can bring a fault back close to reactivation. The safe operating
bandwidth therefore depends not only on the absolute pressure values, but also on how
much the pressure changes with respect to the amounts that previously caused (or did
not cause) slip. A physically consistent operational rule is therefore to limit the pressure excursions
during CGI/CCS and UGS/UHS so that they do not exceed the pressure changes that
were sustained during PP without inducing slip. The following operational guidelines summarize how this principle
translates into practical limits for pressure cycling in depleted gas reservoirs.}

%\revB{Once a fault is activated during primary
%production, the subsequent pressure cycles during CGI/CCS and UGS/UHS reload the
%fault along a different stress path due to the permanent stress redistribution induced by
%slip. As a consequence, a pressure variation applied in the opposite direction may bring
%the system back toward failure even if the absolute pressure remains below $P_i$.  
%A physically consistent operational rule is therefore to limit the pressure excursions
%during CGI/CCS and UGS/UHS so that they do not exceed the pressure changes that
%were sustained during PP without inducing slip. On this basis, the following guidelines
%directly follow from the mechanisms identified in the sensitivity analysis.}
%%
They are necessarily qualitative because of the theoretical/general framework of the modeling application and the quasi-static nature of the implemented model, which properly simulates the possible inception of fault slip but not the seismic evolution. The numerous scenarios investigated within the study have clearly revealed that fault failure is more likely to happen during CGI/CCS and UGS/UHS in depleted reservoirs when:
\begin{enumerate}
    \item fault reactivation is occurred during PP. We refer to the pressure of seismic occurrence as $P_{seis,PP}$. As an example, in scenario \#4b faults are not active during PP and remain far from critical conditions during CGI and UGS as well (Figure~\ref{fig:CH4-6_chi_zoom} and Table~\ref{table:ranking});
    \item the reservoir pressure approaches $P_{max,CGI/ST}$, $P_{max,UGS/UHS}$, or $P_{min,UGS/UHS}$ (e.g., Figures~\ref{fig:CH4-C2-chi} and~\ref{fig:C2_gases}).
\end{enumerate}

The resulting guidelines, which link $P_{max,CGI/CCS}$, $P_{max,UGS/UHS}$, and $P_{min,UGS/UHS}$ to $P_{min,PP}$ and $P_{seis,PP}$, can be stated as follows. The outcome of scenario C2 is particularly illustrative in this regard (Figures~\ref{fig:CH4-C2-chi} and~\ref{fig:C2_gases}).

\begin{figure}%[tbp]
    \centering
    \includegraphics[width=1\linewidth]{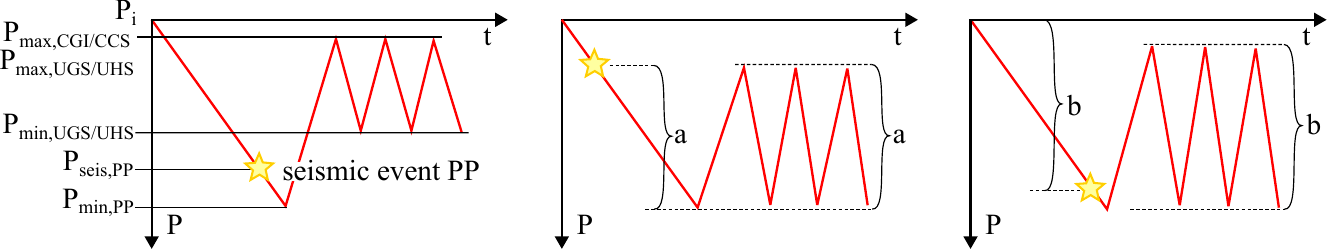}
    \caption{Left: pressure values to be accounted for in the definition of operational safety guidelines. Middle and right: safe pressure variation during CGI, CCS, and UGS/UHS phases in the case of a seismic event occurring during PP at pressure close to $P_i$ and $P_{min,PP}$, respectively.}
    \label{fig:P_bandwidth}
\end{figure}

\begin{itemize}
    \item If a fault reactivation occurs during PP, the pressure change $\Delta P$ spanned during CGI/CCS and UGS/UHS phases should be kept smaller than $|P_i - P_{min,PP}|$. Indeed, a number of investigated scenarios %(see, for example, the outcomes of scenarios \#2, Figure~\ref{fig:ref_chi}) 
    reveals that fault activation during PP leads to a stress redistribution and a new deformed ``balanced'' configuration that is newly loaded, in the opposite direction, when the pressure variation changes its sign. \revB{This behavior is clearly observed in Scenario~\#2
    (Figure~\ref{fig:ref_chi}), Scenario~C2 (Figures~\ref{fig:CH4-C2-chi} and ~\ref{fig:C2_gases}), and also in Scenarios~\#3 and \#4a, where the fault approaches criticality during CGI despite $p < P_i$.}
    A reasonable rule is to keep $\Delta P$ smaller than maximum between $|P_{seis,PP} - P_{min,PP}|$ (range highlighted with \textit{a} in the middle panel of Figure~\ref{fig:P_bandwidth}) and $|P_i - P_{seis,PP}|$ (range highlighted with \textit{b} in the right panel of Figure~\ref{fig:P_bandwidth}), i.e., the maximum pressure difference experienced by the reservoir during PP without the occurrence of a seisimic event. The classification into these two classes is determined by the value of $P_{seis,PP}$ in relation to $P_i$ and $P_{min,PP}$, as illustrated in Figure~\ref{fig:P_bandwidth}, particularly in the middle and right panels.
    %Generally, case a) or case b) prevails depending on the value of $P_{seis,PP}$ with respect to $P_i$ and $P_{min,PP}$ (Figure~\ref{fig:P_bandwidth}-Middle and Right).
    \item In case a, represented in the middle panel of Figure~\ref{fig:P_bandwidth}, $P_{max,CGI/CCS}$ and $P_{max,UGS/UHS}$ should be kept below $P_i$. A reasonable rule of thumb could be to keep $P_{max,CGI/CCS}$ and $P_{max,UGS/UHS}$ smaller than $P_{seis,PP}$.
    \item In case b, right panel of Figure~\ref{fig:P_bandwidth}, $P_{min,UGS/UHS}$ should be kept above $P_{min,PP}$. A reasonable rule of thumb could be to keep $P_{min,UGS/UHS}$ larger than $P_{seis,PP}$.
    \item If during PP activation occurs on a fault that separates two reservoir compartments, during CGI/CCS or UGS/UHS the pressure difference between adjacent blocks should be kept safely less than $P_i - P_{seis,PP}$ (scenario \#5b, Figure~\ref{fig:s7-tau}).
    \item If no activation occurs during PP, $P_{max,CGI/CCS}$ can equate $P_i$ with no particular risk of unexpected events during CGI and UGS. Moreover, $\Delta P$ during UGS/UHS can safely span the whole pressure change between $P_i$ and $P_{min,PP}$. %In fact, the system lies in reloading conditions and behaves practically in an elastic way within the pressure range already experienced
    In fact, the system operates under reloading conditions and exhibits predominantly elastic behavior within the pressure range experienced previously 
    (scenario \#4b, Figure~\ref{fig:CH4-6_chi_zoom}).
\end{itemize}
A graphical representation of the guidelines is provided in Figure~\ref{fig:flowsheet}.

\begin{figure}[H]
    \centering
    \includegraphics[width=1\linewidth]{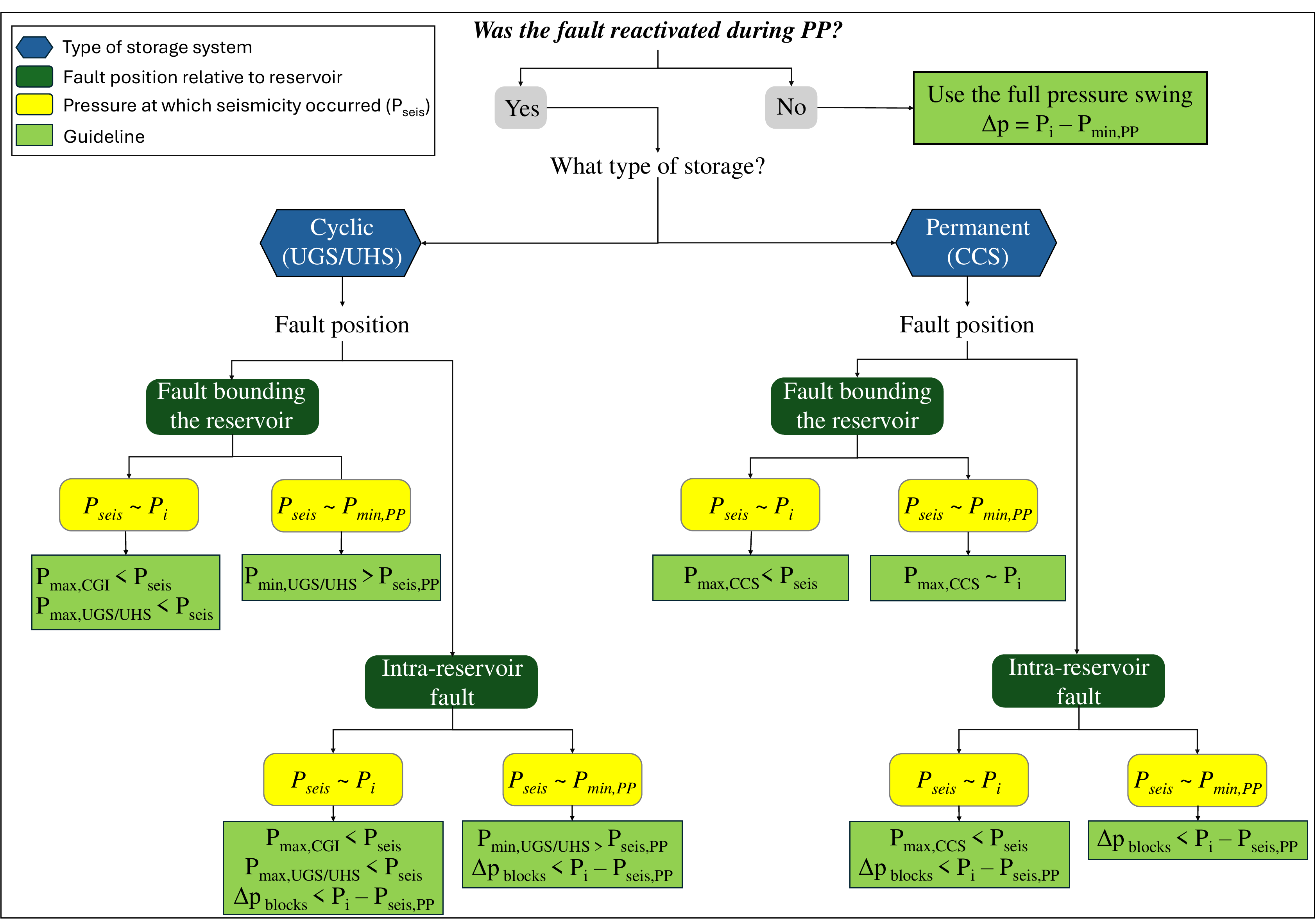}
    \caption{Graphycal representation of the operational guidelines for a compartmentalized, previously produced gas reservoir. Aa key factor is represented by the pressure at which fault reactivation occurred during PP.}
    \label{fig:flowsheet}
\end{figure}

\subsection{Limitations of the current modeling approach}
The modeling framework adopted in this study involves several simplifications that are necessary to address the physical problem at the scale of real-world gas reservoirs. The principal limitations are discussed in this section.

Although the model setup is based on representative geological configurations, both the fault-reservoir geometry and the spatial distribution of physical properties are simplified to allow the problem to be computationally manageable at the basin scale. Geometries are schematized and material properties are homogeneously assigned per domain block. As a result, the model does not capture fine-scale heterogeneities, such as local discontinuities, grain-scale damage zones, or small-scale structural complexity. It is crucial to highlight that fault reactivation is largely site-sensitive, depending on the geometry of the fault/reservoir (e.g., the presence of sloped faults, dislocation of the reservoir compartments), differential pore pressure between adjacent reservoir compartments and within each reservoir block, the geomechanical properties of the faults and the reservoir, caprock, under- and over-burden. Therefore, specific investigations are of paramount importance to characterize the actual reservoir setting and quantify more specific bounds.

For CH$_4$ and H$_2$, the loading history follows a synthetic temporal evolution, designed to replicate typical seasonal injection–withdrawal amplitudes observed in Dutch UGS operations. Unlike actual field conditions, where pressure amplitudes tend to increase gradually over time, our model imposes high-amplitude cycles from the beginning of UGS/UHS ($\sim$50\% of the depletion range) (Figure~\ref{fig:geol-seismic}, bottom-right panel). This idealization is consistent with long-term UGS behavior (see, e.g., %Norg and Grijpskerk 
\cite{TNO15}) and allows us to evaluate the fault response under representative upper-bound stress conditions.

The simulations adopt a one-way hydro-mechanical coupling strategy, where the pore pressure $p$ is provided by a flow model and its contribution for the geomechanical model is regarded as an external source of strength. This modeling choice is recognized as an approximation: it may overestimate pressure buildup near faults and anticipate fault reactivation compared to fully coupled simulations, especially under long-term or cyclic loading. Previous studies have shown that the divergence between one-way and two-way coupling can increase significantly over time \citep{prevost2013}. However, the use of one-way coupling in this study is a practical compromise and it is fully warranted on the space and time scale of interest \citep{bui2019,gambolati2000}. It enables the execution of a large number of 3D, elasto-plastic simulations with a fault system and realistic boundary conditions, while still capturing the relative sensitivity of fault stability to pressure changes. 

The model adopts a quasi-static formulation and does not include inertial effects or dynamic rupture propagation \citep{SodM1}. As a result, it cannot distinguish between aseismic and unstable seismic fault slip, nor simulate the nucleation and growth of dynamic instabilities. Fault activation is evaluated based on a static failure criterion (in our case, Coulomb's formulation), and is used here solely as an indicator of the onset of slip. The model is therefore intended to identify critical loading conditions for reactivation, not to reproduce the full physics of earthquake ruptures.

Building on the previous point, the operational guidelines derived from our simulations are formulated as conservative thresholds that aim to reduce the likelihood of any fault reactivation, regardless of whether slip would manifest aseismically or seismically. %As suggested by field and laboratory studies \cite{}, aseismic slip can occur without immediate hazard, and seismic nucleation typically requires additional conditions such as sufficient fault patch length and stress drop. Accordingly, 
%The criteria derived from fault reactivation are defined conservatively, aiming to provide upper-bound thresholds suitable for guiding safety-focused operational decisions. 
It is worth noting that these guidelines are not universally applicable: they are restricted to a specific class of reservoirs, namely previously depleted gas fields with well-defined geometries and pressure histories, as commonly found in the Dutch subsurface. Extrapolation to other types of storage sites would require additional modeling and site-specific calibration.

Moreover, the need of more specific analyses holds in relation to the geochemical effects on fault mechanical properties caused by the substitution of formation fluids with CO$_2$, H$_2$, or N$_2$. In the modeling approach used in this work, these potential effects on faults have only been superficially addressed, due to the lack of laboratory testing. %which poses a significant challenge in gathering quantitative information to incorporate them into the modeling analysis.
Finally, we underline that the above rules are aimed at avoiding, or at least limiting, the probability of fault activation. For example, according to \cite{wassing_2017}, ``after the onset of fault reactivation, a further pore pressure decline of 3.7~MPa (or 1.6~MPa with no fault offset) and a critically stressed length of approximately 30~m is needed for the nucleation of a seismic event''. Therefore, since fault reactivation can develop with aseismic slip, the recommendations listed above can be considered conservative.

%To add:
% -  substantiate the use of X max - address the effect of cell size

% - no pressure change in seal and underburden - describe assumptions properly. Pressure diffusion in seal and basement are important for the stress changes , as these would smooth out the stress singularities

\revB{Another important limitation concerns the treatment of pore pressure in the seal and underburden. In the current modeling framework, these units are not included in the flow simulation and are assumed to remain at their initial hydrostatic pressure throughout all production and storage stages. 
However, pressure diffusion into the seal and underburden would act to smooth pressure gradients at the reservoir boundaries, thereby reducing stress localization and mitigating the development of stress singularities at fault–layer intersections.}
% Effects of locally higher injection pressures during CGI, UGS + % Effects of cooling
\revB{Moreover, spatial pressure heterogeneities and locally elevated injection pressures 
that may develop near wells during CGI and UGS are not represented. Thermo-mechanical 
effects associated with cooling during gas injection are also neglected.}

% - Bounding faults are the most critical in the current setup, however, this is typically not what is observed or expected from other modeling studies. Bounding faults often have large offsets and are juxtaposed against the Zechstein salt or clay-rich overburden. These rock are prone to creep, likely causing a stable stress on bounding faults. In addition, frictional sliding along clay-filled faults is likely stable. 

% Role of elasto-plastic and visco-plastic deformation of the reservoir during PP, hardening effects, and effect of this on the stress path are essential to consider. Consider the works of Pijnenburg et al., Hol et al. who show a different stress path parameter during loading and unloading. 
\revB{In the present modeling framework, reservoir deformation is represented using a simplified mechanical description that does not explicitly account for elasto-plastic or visco-plastic deformation, strain hardening, or time-dependent effects. Laboratory studies on Slochteren sandstone show that inelastic deformation contributes significantly to reservoir compaction already at small strains and leads to stress-path hysteresis, with different stress evolution during loading and unloading phases \citep{pijnenburg2019inelastic}.
As a consequence, the stress paths modeled here during pressure depletion and subsequent pressurization are largely elastic and reversible, and do not capture the irreversible stress-path shifts documented experimentally. Including inelastic and time-dependent reservoir behavior would likely modify the magnitude and timing of stress transfer to faults, particularly by reducing elastic stress recovery during pressurization and by shifting the conditions at which fault criticality is reached.}

% For CO2, H2, N2 changes in Young’s modulus are considered, but it is Poisson’s ratio that drives the poro-elastic stress changes. Do CO2, H2, N2 change Poisson’s ratio? Please discuss.

\subsection{Conclusions and future prospects}
Depleted gas reservoirs can serve not only as storage facilities for natural gas (i.e., CH$_4$) but also for CO$_2$, N$_2$, and H$_2$. The injection and/or withdrawal of these fluids induce changes in the hydraulic and mechanical state of the reservoir. Consequently, deformations and variations in the stress regime develop beneath the subsurface, potentially leading to (seismic/aseismic) reactivation of faults near the reservoir. While the majority of human-induced seismic activities can be associated to injection or withdrawal of fluids at pressures above the original formation pressure causing significant pressure decline \citep{ellsworth_2013_injectioninduced}, which trigger shear stress along faults to reach their limit strength, a few recorded events do not fit this explanation. These seismic events, somehow ``unexpected'', develop in reservoirs where:
\begin{itemize}
    \item seismic events already occurred during the primary production phase;
    \item in correspondence to a pressure value already experienced by the reservoir during PP.
\end{itemize}

%\sout{The main aim of this work has been to understand whether, and under which circumstances, faults can be ``unexpectedly'' reactivated during USS activities. A one-way quasi-static coupled strategy is adopted to deal with interaction between fluid and the continuous-discontinuous porous formation. The study has been carried out using the typical geological setting of reservoirs located in the Rotliegend formation, in the Netherlands, where systems of almost orthogonal faults split the formation into various compartments. Four underground storage plans, i.e., cyclic storage of CH$_4$ and H$_2$ and permanent storage of CO$_2$ and N$_2$, have been investigated. This objective is achieved by analyzing which are the main factors controlling the reactivation of faults under permanent or cyclic storage conditions. The consequences of using different fluids, which affect the pressure evolution over time and the Rotliegend geomechanical properties, have been evaluated.}

The simulation of various realistic scenarios have allowed to define the critical factors influencing fault activation during CGI/CCS and UGS/UHS cycles. The stability of the faults bounding the reservoir compartments is mainly jeopardized by an initial stress regime with small horizontal principal components, low friction angle, large pressure change because of injection/withdrawal, and significant contrast between the reservoir and the over-, side-, and underburden stiffness. Concerning faults located between producing blocks, the drivers mainly influencing fault instability are the geometrical setting of the fault/reservoir system, i.e., the offset between reservoir compartments and the fault dip angle, together with the different pressure change in adjacent compartments.

Notice that the pressure recovery and drop addressed in the simulations have been defined based on regulations in the Netherlands, which do not allow the pressure to rise above the initial natural value. Due to this constrain, the mechanisms causing fault reactivation remain similar, irrespective of the fluid considered.
The possible mechanical weakening and hardening of the reservoir associated to non-natural pore fluids interaction with the solid grains does not impact significantly the outcomes.

The modeling outcomes have enabled the formulation of general guidelines to define safe operational bandwidths for USS sites, specifically the pressure range over which ``unexpected'' seismic event can be excluded. These minimum and maximum pressure thresholds are closely linked to the pressures at which seismic events occurred during primary production. The occurrence of a seismic event during PP provides valuable insights for delineating pressure bandwidths within which fault reactivation is highly unlikely during UGS, UHS, CGI, and CCS activities.

These conservative recommendations can serve as a preliminary methodology for reducing the potential risks associated to seismic activity in similar contexts. However, more specific and in-depth evaluations must follow, taking into account the peculiarities of each individual case study.

\section{Supplementary Materials}

\begin{figure}%[htbp]
    \centering
    \vcenteredhbox{\includegraphics[width=1\textwidth]{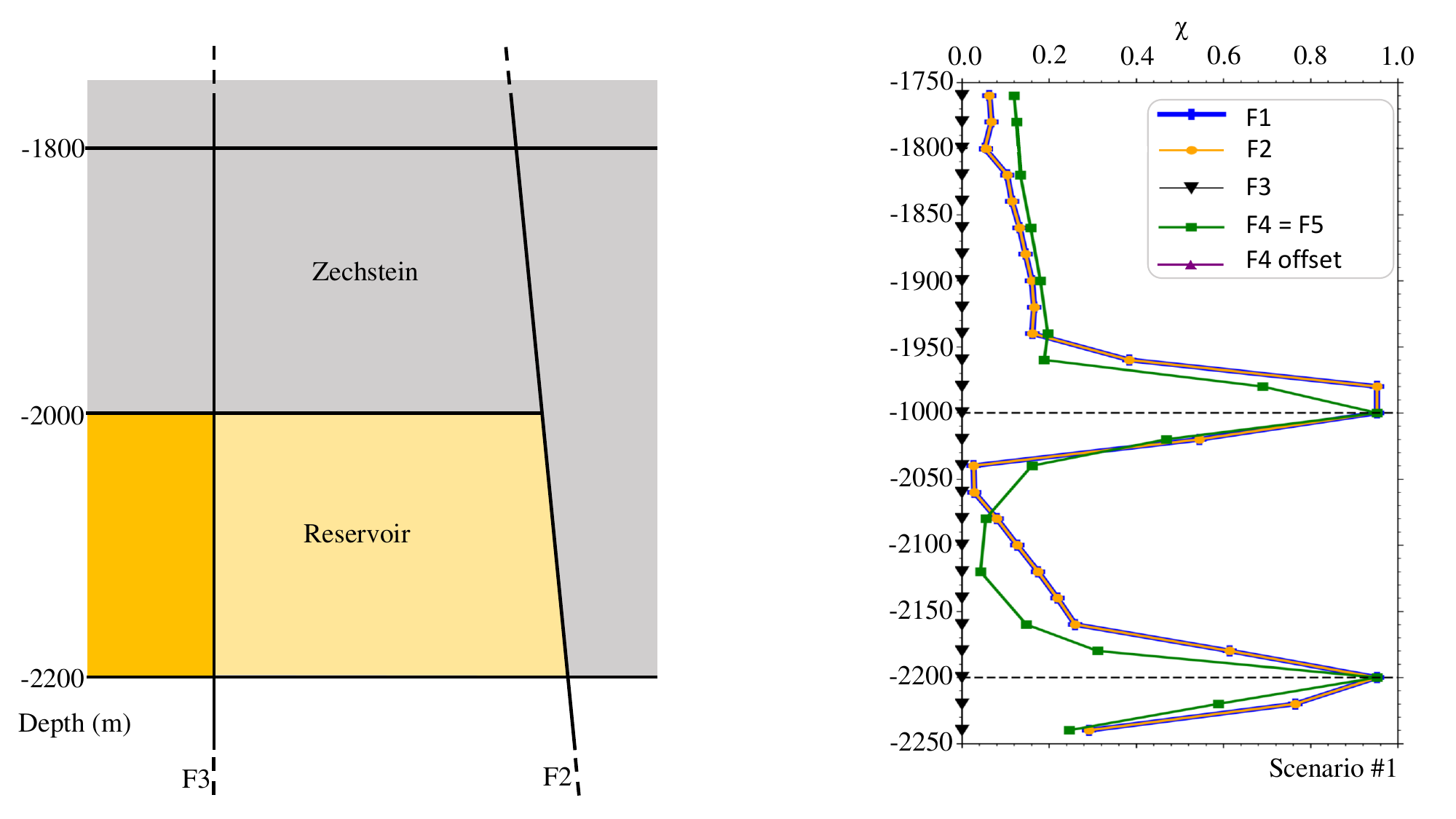}}
    \caption{ Scenario \#1 ("Reference"): depth profiles of the criticality index $\chi$ on all faults at the end of PP (loading step 10). High criticality develops at the reservoir top and bottom.}
    \label{fig:ref_chi_vs_depth}
\end{figure}

\begin{figure}%[htbp]
    \centering
    \includegraphics[width=1\textwidth]{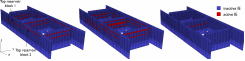}
    \caption{Scenario \#2d: active interface elements (IE) on the fault system at loading step 6 (left), loading step 10 (center), and loading step 12.5 (right). Only the elements at the reservoir top, and then bottom, are activated.}
    \label{fig:CH4-3_3D_active}
\end{figure}
%\input{flowsheet-guidelines}

%% The Appendices part is started with the command \appendix;
%% appendix sections are then done as normal sections
%% \appendix

%% \section{}
%% \label{}

% \section*{CRediT authorship contribution statement}
% \textbf{Selena Baldan}: Investigation, Formal analysis, Writing - original draft.
% \textbf{Massimiliano Ferronato}: Conceptualization, Methodology, Writing - review and editing, Supervision.
% \textbf{Andrea Franceschini}: Methodology, Software, Writing - review and editing.
% \textbf{Carlo Janna}: Software.
% \textbf{Claudia Zoccarato}: Conceptualization, Methodology, Formal analysis, Investigation.
% \textbf{Matteo Frigo}: Software, Investigation.
% \textbf{Giovanni Isotton}: Conceptualization, Methodology, Formal analysis, Software, Funding acquisition.
% \textbf{Cristiano Collettini}: Formal analysis, Writing - review and editing.
% \textbf{Chiara Deangeli}: Methodology, Formal analysis, Writing - review and editing.
% \textbf{Vera Rocca}: Methodology, Formal analysis, Writing - review and editing.
% \textbf{Francesca Verga}: Methodology, Formal analysis, Writing - review and editing.
% \textbf{Pietro Teatini}: Conceptualization, Methodology, Formal analysis, Writing - review and editing, Supervision, Funding acquisition.

\section*{Declaration of competing interest}
The authors declare that they have no known competing financial interests or personal
relationships that could have appeared to influence the work reported in this paper.

\section*{Acknowledgements}
This research was supported by the State Supervision of Mines (SodM), Ministry of Economic Affairs (The Netherlands), projects KEM01 ``Safe Operational Bandwidth of Gas Storage Reservoirs'' and KEM39 ``Study within the Mining Effects Knowledge Program (KEM-39) on the cyclic storage of gases in the Netherlands''. %Portions of this work were performed within the 2020 INdAM-GNCS project ``Optimization and advanced linear algebra for PDE-governed problems''.
M.F., A.F., and C.J. from the University of Padova are members of the Gruppo Nazionale Calcolo Scientifico -- Istituto Nazionale di Alta Matematica (GNCS-INdAM).
Computational resources were provided by University of Padova Strategic Research Infrastructure Grant 2017: ``CAPRI: Calcolo ad Alte Prestazioni per la Ricerca e l'Innovazione''.

%% If you have bibdatabase file and want bibtex to generate the
%% bibitems, please use
%%
%%  \bibliographystyle{elsarticle-num}
%%  \bibliography{<your bibdatabase>}
%\Urlmuskip=0mu plus 1mu
%\bibliographystyle{elsarticle-num-names}
%\bibliographystyle{elsarticle-harv}
\bibliographystyle{elsarticle-harv}
\bibliography{biblio}

@article{borello_2024_underground,
  author = { Borello, Eloisa Salina and Bocchini, Sergio and Chiodoni, Angelica and Coti, Christian and Fontana, Marco and Panini, Filippo  and Peter, Costanzo and  Pirri, Candido Fabrizio and Tawil, Michel and Mantegazzi, Andrea and Marzano, Francesco and Pozzovivo, Vincenzo and Verga, Francesca and Viberti, Dario},
  month = {01},
  pages = {394},
  publisher = {Multidisciplinary Digital Publishing Institute},
  title = {Underground Hydrogen Storage Safety: Experimental Study of Hydrogen Diffusion through Caprocks},
  doi = {10.3390/en17020394},
  volume = {17},
  year = {2024},
  journal = {Energies}
}

@article{Bui_etal17,
  title        = {
    {Fault reactivation mechanisms and dynamic rupture modelling of depletion-induced
    seismic events in a Rotliegend gas reservoir}
  },
  author       = {
    Buijze, Loes and van den Bogert, Peter A.J. and Wassing, Brecht B.T. and Orlic,
    Bogdan and ten Veen, Johan
  },
  year         = 2017,
  journal      = {Netherlands Journal of Geosciences},
  volume       = 96,
  number       = 5,
  pages        = {s131-s148},
  doi          = {10.1017/njg.2017.27}
}

@article{bui2019,
  title={Nucleation and arrest of dynamic rupture induced by reservoir depletion},
  author={Buijze, L and van den Bogert, Peter AJ and Wassing, BBT and Orlic, B},
  journal={Journal of Geophysical Research: Solid Earth},
  volume={124},
  number={4},
  pages={3620--3645},
  year={2019},
  publisher={Wiley Online Library}
}

@article{Can_etal19,
  title        = {
    Depletion-Induced Seismicity at the {G}roningen Gas Field: {C}oulomb Rate-and-State
    Models Including Differential Compaction Effect
  },
  author       = {
    Candela, Thibault and Osinga, Sander and Ampuero, Jean-Paul and Wassing, Brecht and
    Pluymaekers, Maarten and Fokker, Peter A. and van Wees, Jan-Diederik and de Waal,
    Hans A. and Muntendam-Bos, Annemarie G.
  },
  year         = 2019,
  journal      = {Journal of Geophysical Research: Solid Earth},
  volume       = 124,
  number       = 7,
  pages        = {7081--7104},
  doi          = {10.1029/2018JB016670}
}

@article{cappa_2019_stabilization,
  author = {Cappa, F. and Scuderi, M. M. and Collettini, C. and Guglielmi, Y. and Avouac, J.P.},
  month = {03},
  title = {Stabilization of fault slip by fluid injection in the laboratory and in situ},
  doi = {10.1126/sciadv.aau4065},
  volume = {5},
  year = {2019},
  journal = {Science Advances}
}

@article{Ces_etal14,
  title        = {
    The 2013 {S}eptember-{O}ctober seismic sequence offshore {S}pain: a case of seismicity
    triggered by gas injection?
  },
  author       = {
    Cesca, Simone and Grigoli, Francesco and Heimann, Sebastian and Gonzalez, Alvaro and
    Buforn, Elisa and Maghsoudi, Samira and Blanch, Estefania and Dahm, Torsten
  },
  year         = 2014,
  journal      = {Geophysical Journal International},
  volume       = 198,
  number       = 2,
  pages        = {941--953},
  doi          = {10.1093/gji/ggu172}
}

@article{collettini_2008_fault,
  author = { Collettini, C. and Cardellini, C. and Chiodini, G. and  De Paola, N. and Holdsworth, R E and Smith, S.A.F},
  month = {01},
  pages = {175-194},
  publisher = {Geological Society of London},
  title = {Fault weakening due to {CO}$_2$ degassing in the {N}orthern {A}pennines: short- and long-term processes},
  doi = {10.1144/sp299.11},
  urldate = {2023-09-13},
  volume = {299},
  year = {2008},
  journal = {Geological Society, London, Special Publications}
}

@article{Def_etal95,
  title        = {
    {Microseismic surveying and repeated VSPs for monitoring an underground gas storage
    reservoir using permanent geophones}
  },
  author       = {Deflandre, J. P. and Laurent, J. and Michon, D. and Blondin, E.},
  year         = 2018,
  journal      = {First Break},
  volume       = 13,
  number       = 4,
  pages        = {129--138},
  doi          = {10.3997/1365-2397.1995008}
}

@article{ellsworth_2013_injectioninduced,
  author = {Ellsworth, W.L.},
  month = {07},
  pages = {1225942-1225942},
  title = {Injection-Induced Earthquakes},
  doi = {10.1126/science.1225942},
  volume = {341},
  year = {2013},
  journal = {Science}
}

@article{Fou_etal18,
  title        = {{Global review of human-induced earthquakes}},
  author       = {
    Gillian R. Foulger and Miles P. Wilson and Jon G. Gluyas and Bruce R. Julian and
    Richard J. Davies
  },
  year         = 2018,
  journal      = {Earth-Science Reviews},
  volume       = 178,
  pages        = {438--514},
  doi          = {10.1016/j.earscirev.2017.07.008}
}

@article{gambolati2000,
  title={Importance of poroelastic coupling in dynamically active aquifers of the Po river basin, Italy},
  author={Gambolati, Giuseppe and Teatini, Pietro and Ba{\'u}, Domenico and Ferronato, Massimiliano},
  journal={Water Resources Research},
  volume={36},
  number={9},
  pages={2443--2459},
  year={2000},
  publisher={Wiley Online Library}
}

@article{HANGX2013,
title = {The effect of CO2 on the mechanical properties of the Captain Sandstone: Geological storage of CO$_2$ at the Goldeneye field (UK)},
journal = {International Journal of Greenhouse Gas Control},
volume = {19},
pages = {609-619},
year = {2013},
issn = {1750-5836},
doi = {https://doi.org/10.1016/j.ijggc.2012.12.016},
url = {https://www.sciencedirect.com/science/article/pii/S1750583612003283},
author = {Suzanne Hangx and Arjan {van der Linden} and Fons Marcelis and Andreas Bauer},
}

@article{harbert_2020,
  author = {Harbert, William and Goodman, Angela and Spaulding, Richard and Haljasmaa, Igor  and Crandall, Dustin and Sanguinito, Sean and Kutchko, Barbara and Tkach, Mary and Fuchs, Samantha and Werth, Charles J and Tsotsis, Theodore and Dalton, Laura and Jessen, Kristian and Shi, Zhuofan and Frailey, Scott},
  month = {09},
  pages = {103109},
  publisher = {Elsevier BV},
  title = {{CO}$_2$ induced changes in {M}ount {S}imon sandstone: {U}nderstanding links to post {CO}$_2$ injection monitoring, seismicity, and reservoir integrity},
  doi = {10.1016/j.ijggc.2020.103109},
  volume = {100},
  year = {2020},
  journal = {International Journal of Greenhouse Gas Control}
}

@article{Hau_etal18,
  title        = {
    {Assessment of geological factors potentially affecting production-induced seismicity
    in North German gas fields}
  },
  author       = {C. Haug and J.-A. N\"uchter and A. Henk},
  year         = 2018,
  journal      = {Geomechanics for Energy and the Environment},
  volume       = 16,
  pages        = {15--31},
  doi          = {10.1016/j.gete.2018.04.002}
}

@article{Het_etal00,
  title        = {
    Production-Induced Compaction of a Sandstone Reservoir: The Strong Influence of
    Stress Path
  },
  author       = {
    Hettema, M. H. H. and Schutjens, P. M. T. M. and Verboom, B. J. M. and Gussinklo, H.
    J.
  },
  year         = 2000,
  journal      = {SPE Reservoir Evaluation \& Engineering},
  volume       = 3,
  pages        = {342--347},
  doi          = {10.2118/65410-PA}
}

@article{Hun_etal17,
author = {Hunfeld, L. B. and Niemeijer, A. R. and Spiers, C. J.},
title = {Frictional Properties of Simulated Fault Gouges from the Seismogenic {G}roningen Gas Field Under In Situ {P--T}-Chemical Conditions},
journal = {Journal of Geophysical Research: Solid Earth},
volume = {122},
number = {11},
pages = {8969--8989},
doi = {10.1002/2017JB014876},
year = {2017}
}

@article{Hun_etal20,
author = {Hunfeld, Luuk B. and Chen, Jianye and Hol, Sander and Niemeijer, Andr\'e R. and Spiers, Christopher J.},
title = {Healing Behavior of Simulated Fault Gouges From the {G}roningen Gas Field and Implications for Induced Fault Reactivation},
journal = {Journal of Geophysical Research: Solid Earth},
volume = {125},
number = {7},
pages = {e2019JB018790},
doi = {10.1029/2019JB018790},
year = {2020}
}

@article{hun2021,
  title={{Seismic slip-pulse experiments simulate induced earthquake rupture in the Groningen gas field}},
  author={Hunfeld, Luuk B and Chen, Jianye and Niemeijer, Andr{\'e} R and Ma, Shengli and Spiers, Christopher J},
  journal={Geophysical Research Letters},
  volume={48},
  number={11},
  pages={e2021GL092417},
  year={2021},
  publisher={Wiley Online Library}
}

@article{KerWei18,
  title        = {Induced Seismicity},
  author       = {Keranen, Katie M. and Weingarten, Matthew},
  year         = 2018,
  journal      = {Annual Review of Earth and Planetary Sciences},
  volume       = 46,
  number       = 1,
  pages        = {149--174},
  doi          = {10.1146/annurev-earth-082517-010054}
}

@article{Iso_etal19,
  title        = {
    {Robust numerical implementation of a 3D rate-dependent model for reservoir
    geomechanical simulations}
  },
  author       = {
    Isotton, Giovanni and Teatini, Pietro and Ferronato, Massimiliano and Janna, Carlo
    and Spiezia, Nicol\`o and Mantica, Stefano and Volonte, Giorgio
  },
  year         = 2019,
  journal      = {International Journal for Numerical and Analytical Methods in Geomechanics},
  volume       = 43,
  number       = 18,
  pages        = {2752--2771},
  doi          = {10.1002/nag.3000}
}

@article{Jia_etal21,
author = {Jiang, Guoyan and Liu, Lin and Barbour, Andrew J. and Lu, Renqi and Yang, Hongfeng},
title = {Physics-Based Evaluation of the Maximum Magnitude of Potential Earthquakes Induced by the {H}utubi ({C}hina) Underground Gas Storage},
journal = {Journal of Geophysical Research: Solid Earth},
volume = {126},
number = {4},
pages = {e2020JB021379},
doi = {10.1029/2020JB021379},
year = {2021}
}

@article{Liu_etal23,
title = {Evaluation of fault stability and seismic potential for {H}utubi underground gas storage due to seasonal injection and extraction},
journal = {Underground Space},
volume = {13},
pages = {74--85},
year = {2023},
doi = {10.1016/j.undsp.2023.03.006},
author = {Cexuan Liu and Fengshou Zhang and Quan Wang and Bin Wang and Qi Zhang and Bin Xu},
}

@article{Mic_etal23,
author = {Johannes Miocic  and Niklas Heinemann  and Katriona Edlmann  and Jonathan Scafidi  and Fatemeh Molaei  and Juan Alcalde },
title = {Underground hydrogen storage: a review},
journal = {Geological Society, London, Special Publications},
volume = {528},
number = {1},
pages = {73--86},
year = {2023},
doi = {10.1144/SP528-2022-88},
}

@article{Mun_etal15,
  title        = {
    {A guideline for assessing seismic risk induced by gas extraction in the Netherlands}
  },
  author       = {A. G. Muntendam-Bos  and  J. P. A. Roest  and  J. A. de Waal},
  year         = 2015,
  journal      = {The Leading Edge},
  volume       = 34,
  number       = 6,
  pages        = {672--677},
  doi          = {10.1190/tle34060672.1}
}

@article{Seg_etal94,
  title        = {
    {Poroelastic stressing and induced seismicity near the Lacq gas field, southwestern
    France}
  },
  author       = {Segall, Paul and Grasso, Jean-Robert and Mossop, Antony},
  year         = 1994,
  journal      = {Journal of Geophysical Research: Solid Earth},
  volume       = 99,
  number       = {B8},
  pages        = {15423--15438},
  doi          = {10.1029/94JB00989}
}

@article{shi2019,
  title={{Impact of Brine/CO$_2$ exposure on the transport and mechanical properties of the Mt Simon sandstone}},
  author={Shi, Z and Sun, L and Haljasmaa, I and Harbert, W and Sanguinito, S and Tkach, M and Goodman, A and Tsotsis, TT and Jessen, K},
  journal={Journal of Petroleum Science and Engineering},
  volume={177},
  pages={295--305},
  year={2019},
  publisher={Elsevier}
}

@article{tar2020,
  title={Influence of CO2 injection on the poromechanical response of Berea sandstone},
  author={Tarokh, Ali and Makhnenko, Roman Y and Kim, Kiseok and Zhu, Xuan and Popovics, John S and Segvic, Branimir and Sweet, Dustin E},
  journal={International Journal of Greenhouse Gas Control},
  volume={95},
  pages={102959},
  year={2020},
  publisher={Elsevier}
}

@article{Tea_etal20,
  title        = {{About geomechanical safety for UGS activities in faulted reservoirs}},
  author       = {
    Teatini, P. and Zoccarato, C. and Ferronato, M. and Franceschini, A. and Frigo, M.
    and Janna, C. and Isotton, G.
  },
  year         = 2020,
  journal      = {Proceedings of the International Association of Hydrological Sciences},
  volume       = 382,
  pages        = {539--545},
  doi          = {10.5194/piahs-382-539-2020}
}

@incollection{Hul10,
  title        = {
    {Geological factors effecting compartmentalization of Rotliegend gas fields in the
    Netherlands}
  },
  author       = {Van Hulten, F. F. N.},
  year         = 2010,
  booktitle    = {{Reservoir Compartmentalization}},
  publisher    = {Geological Society of London},
  pages        = {301--315},
  doi          = {10.1144/SP347.17}
}

@article{Wee_etal14,
  title        = {
    {Geomechanics response and induced seismicity during gas field depletion in the
    Netherlands}
  },
  author       = {
    Van Wees, J. and Buijze, Loes and Thienen-Visser, Karin and Nepveu, Manuel and
    Wassing, Brecht and Orlic, Bogdan and Fokker, Peter
  },
  year         = 2014,
  month        = 10,
  journal      = {Geothermics},
  volume       = 52,
  pages        = {206--219},
  doi          = {10.1016/j.geothermics.2014.05.004}
}

@article{walsh_2016,
  author = {Walsh, F. Rall and Zoback, Mark D.},
  month = {10},
  pages = {991-994},
  title = {Probabilistic assessment of potential fault slip related to injection-induced earthquakes: Application to north-central {O}klahoma, {USA}},
  doi = {10.1130/g38275.1},
  volume = {44},
  year = {2016},
  journal = {Geology}
}

@article{verga_2018,
  author = {Verga, Francesca},
  month = {05},
  pages = {1245},
  title = {What's Conventional and What's Special in a Reservoir Study for Underground Gas Storage},
  doi = {10.3390/en11051245},
  volume = {11},
  year = {2018},
  journal = {Energies}
}

@article{Vil_etal21,
author = {Vilarrasa, V\'ictor and De Simone, Silvia and Carrera, Jesus and Villase\~nor, Antonio},
title = {Unraveling the Causes of the Seismicity Induced by Underground Gas Storage at {C}astor, {S}pain},
journal = {Geophysical Research Letters},
volume = {48},
number = {7},
pages = {e2020GL092038},
doi = {10.1029/2020GL092038},
year = {2021}
}

@article{Zho_etal19,
  title        = {
    Seismological Investigations of Induced Earthquakes Near the {H}utubi Underground Gas
    Storage Facility
  },
  author       = {Zhou, Pengcheng and Yang, Hongfeng and Wang, Baoshan and Zhuang, Jiancang},
  year         = 2017,
  journal      = {Journal of Geophysical Research: Solid Earth},
  volume       = 124,
  number       = 8,
  pages        = {8753--8770},
  doi          = {10.1029/2019JB017360}
}

@techreport{Bai_etal16,
  title        = {
    {Induced Seismicity in the Bergermeer Field: Hyopcenter Relocation and
    Interpretation}
  },
  author       = {Baisch, Stefan and Koch, Christopher and Voros, Robert and Rothert, Elmar},
  year         = 2016,
  institution  = {Taqa Energy B.V. TAQA003}
}

@techreport{Gau03,
  title        = {
    {Carboniferous-Rotliegend Total Petroleum System Description and Assessment Results
    Summary}
  },
  author       = {Donald L. Gautier},
  year         = 2003,
  institution  = {U.S. Geological Survey Bulletin 2211}
}

@techreport{groenenberg_2020_largescale,
  title        = {
    {Large-Scale Energy Storage in Salt Caverns and Depleted Fields (LSES) - Project
    Findings}
  },
  author       = {
    Groenenberg, R. M. and Koornnef, J. M. and Sijm, J. P. M. and Janssen, G. J. M. and
    Morales Espa\~{n}a, G. A. and van Stralen, J. and Hernandez-Serna, R. and Smekens, K.
    E. L. and Juez-Larre, J. and Goncalves Machado, C. and Wasch, L. J. and Dijkstra, H.
    E. and Wassing, B. B. L. and Orlic, B. and Brunner, L. G. and van der Valk, K. and
    van Unen, M. and Hajonides van der Meulen, T. C. and Kranenburg-Bruinsma, K. J. and
    Winters, E. and Puts, H. and van Popering-Verkerk, J. and Duijn, M.
  },
  year         = 2020,
  institution  = {TNO R12006}
}

@article{JAG07,
author = {de Jager, Jan and Geluk, M.C.},
year = {2007},
month = {01},
pages = {237-260},
title = {Petroleum geology},
journal = {Geology of the Netherlands. Royal Dutch Academy of Arts and Sciences, Amsterdam}
}

@techreport{MIT09,
  title        = {{Technical Review of Bergermeer Seismicity Study TNO Report 2008-U-R1071/B}},
  author       = {Hager, B. H. and Toksoz, M. N.},
  year         = 2009,
  institution  = {Massachusetts Institute of Technology}
}

@techreport{NAM16,
  title        = {{Norg UGS fault reactivation study and implications for seismic threat}},
  author       = {{Nederlandse Aardolie Maatschappij BV}},
  year         = 2016,
  institution  = {NAM EP201610208045}
}

@techreport{TNO15,
  title        = {
    {Injection-Related Induced Seismicity and its Relevance to Nitrogen Injection:
    Description of Dutch Field Cases}
  },
  author       = {TNO},
  year         = 2015,
  institution  = {TNO R10906}
}

@incollection{deJ07,
  title        = {{Geological development}},
  author       = {De Jager, J.},
  year         = 2007,
  booktitle    = {{Geology of the Netherlands}},
  publisher    = {Royal Netherlands Academy of Arts and Sciences, Amsterdam},
  pages        = {5--26},
  editors      = {Wong, T. E., Batjes, D. A. J. and De Jager, J.}
}

@article{de2017,
  title={Geology of the Groningen field--an overview},
  author={De Jager, Jan and Visser, Clemens},
  journal={Netherlands Journal of Geosciences},
  volume={96},
  number={5},
  pages={s3--s15},
  year={2017},
  publisher={Cambridge University Press}
}

@inproceedings{Kra_etal13,
  title        = {
    Microseismic Monitoring and Subseismic Fault Detection in an Underground Gas
    Storage
  },
  author       = {Kraaijpoel, D.A. and Nieuwland, D.A. and Dost, B.},
  year         = 2013,
  booktitle    = {{Proc. 4th EAGE Passive Seismic Workshop, Amsterdam, Netherlands}},
  publisher    = {European Association of Geoscientists \& Engineers},
  pages        = {1--3},
  doi          = {10.3997/2214-4609.20142354}
}

@inproceedings{Orl_etal13,
  title        = {
    {Field scale geomechanical modeling for prediction of fault stability during
    underground gas storage operations in a depleted gas field in the Netherlands}
  },
  author       = {Orlic, B. and Wassing, B. B. T. and Geel, C. R.},
  year         = 2013,
  booktitle    = {{47th US Rock Mechanics / Geomechanics Symposium}},
  publisher    = {Paper \#ARMA 13-300, American Rock Mechanics Association},
  pages        = {1--11}
}

@inproceedings{Tea_etal19,
  title        = {
    {Gas storage in compartmentalized reservoirs: a numerical investigation on possible
    ``unexpected'' fault activation}
  },
  author       = {
    Teatini, P. and Ferronato, M. and Franceschini, A. and Frigo, M. and Janna, C. and
    Zoccarato, C.
  },
  year         = 2019,
  booktitle    = {{53rd US Rock Mechanics / Geomechanics Symposium}},
  publisher    = {Paper \#ARMA 19-1991, American Rock Mechanics Association},
  pages        = {1--9}
}

@phdthesis{Uta17,
  title        = {
    {Recent Intraplate Earthquakes in Northwest Germany - Glacial Isostatic Adjustment
    and/or a Consequence of Hydrocarbon Production}
  },
  author       = {Uta, Philipp},
  year         = 2017,
  doi          = {10.15488/9088},
  school       = {Leibniz University Hannover}
}

@misc{SodM1,
  title        = {
    {Unexpected fault activation in underground gas storage. Part I: Mathematical model
    and mechanisms}
  },
  author       = {
    Franceschini, A. and Baldan, S. and Ferronato, M. and Janna, C. and Zoccarato, C. and
    Frigo, M. and Isotton, G. and Teatini, P.
  },
  year         = 2024,
  month        = {08},
  doi          = {10.48550/arXiv.2308.02198}
}

@misc{OPM2015,
  title        = {{Open Porous Media}},
  author       = {{The Open Porous Media Initiative}},
  year         = 2023,
  note         = {Accessed: 2024-09-03},
  howpublished = {\url{https://opm-project.org}}
}

@article{rasmussen_2021_the,
  title        = {{The Open Porous Media Flow reservoir simulator}},
  author       = {
    Rasmussen, A. and Sandve, T. and Bao, K. and Lauser, A. and Hove, J. and Skaflestad,
    B. and Kl\"ofkorn, R. and Blatt, M. and Rustad, A. and S\ae vareid, O. and Lie, K.
    and Thune, A.
  },
  year         = 2021,
  month        = {01},
  journal      = {Computers \& Mathematics with Applications},
  volume       = 81,
  pages        = {159--185},
  doi          = {10.1016/j.camwa.2020.05.014},
  urldate      = {2021-11-02}
}

@misc{schlumberger_2014_eclipse,
  title        = {{ECLIPSE Industry Reference Reservoir Simulator}},
  author       = {Schlumberger},
  year         = 2014,
  url          = {https://www.software.slb.com/products/eclipse},
  organization = {www.software.slb.com}
}

@article{alshafi_2023_a,
  title        = {
    {A review on underground gas storage systems: Natural gas, hydrogen and carbon
    sequestration}
  },
  author       = {Al-Shafi, Manal and Massarweh, Osama and Abushaikha, Ahmad S. and Bicer, Yusuf},
  year         = 2023,
  month        = 12,
  journal      = {Energy Reports},
  volume       = 9,
  pages        = {6251--6266},
  doi          = {10.1016/j.egyr.2023.05.236}
}

@book{netzero2050,
  title        = {{Going climate-neutral by 2050 - A strategic long-term vision for a
                   prosperous, modern, competitive and climate-neutral EU economy}},
  author       = {{European Commission and Directorate-General for Climate Action}},
  year         = 2019,
  publisher    = {Publications Office},
  doi          = {10.2834/02074}
}

@misc{nlog,
  title        = {{Natural resources and geothermal energy in the Netherlands - Annual review 2020}},
  author       = {{Ministry of Economic Affairs and Climate Policy}},
  year         = 2021,
  url          = {https://www.nlog.nl/sites/default/files/2021-06/pre-publ_ch4_yearbook2020_en.pdf},
  urldate      = {2023-08-18}
}

@article{mun_etal22,
  title        = {{An overview of induced seismicity in the Netherlands}},
  author       = {
    Muntendam-Bos, Annemarie G. and Hoedeman, Gerco and Polychronopoulou, Katerina and
    Draganov, Deyan and Weemstra, Cornelis and van der Zee, Wouter and Bakker, Richard R.
    and Roest, Hans
  },
  year         = 2022,
  journal      = {Netherlands Journal of Geosciences},
  publisher    = {Cambridge University Press},
  volume       = 101,
  pages        = {e1},
  doi          = {10.1017/njg.2021.14}
}

@article{rohmer_2016_mechanochemical,
  title        = {{Mechano-chemical interactions in sedimentary rocks in the context of
                   CO$_2$ storage: Weak acid, weak effects?}},
  author       = {Rohmer, J. and Pluymakers, A. and Renard, F.},
  year         = 2016,
  month        = {06},
  journal      = {Earth-Science Reviews},
  volume       = 157,
  pages        = {86--110},
  doi          = {10.1016/j.earscirev.2016.03.009}
}

@article{peter_2022_a,
  title        = {
    A Review of the Studies on {CO}$_2$-Brine-Rock Interaction in Geological Storage
    Process
  },
  author       = {Peter, Ameh and Yang, Dongmin and Eshiet, Kenneth Imo-Imo Israel and Sheng, Yong},
  year         = 2022,
  month        = {04},
  journal      = {Geosciences},
  volume       = 12,
  pages        = 168,
  doi          = {10.3390/geosciences12040168}
}

@article{rimmel_2010_evolution,
  title        = {
    Evolution of the Petrophysical and Mineralogical Properties of Two Reservoir Rocks
    Under Thermodynamic Conditions Relevant for {CO}$_2$ Geological Storage at 3 km Depth
  },
  author       = {Rimmel\'e, G. and Barlet-Gou\'edard, V. and Renard, F.},
  year         = 2010,
  month        = {07},
  journal      = {Oil \& Gas Science and Technology - Revue de l'Institut Fran\c{c}ais du P\'etrole},
  volume       = 65,
  pages        = {565--580},
  doi          = {10.2516/ogst/2009071}
}

@article{bolourinejad_2015_chemical,
  title        = {
    {Chemical effects of sulfur dioxide co-injection with carbon dioxide on the reservoir
    and caprock mineralogy and permeability in depleted gas fields}
  },
  author       = {Bolourinejad, Panteha and Herber, Rien},
  year         = 2015,
  month        = {08},
  journal      = {Applied Geochemistry},
  volume       = 59,
  pages        = {11--22},
  doi          = {10.1016/j.apgeochem.2015.03.003}
}

@phdthesis{bolourinejad2015effects,
  author       = {Bolourinejad, Panteha},
  title        = {Effects of impurities on subsurface CO$_2$ storage in gas fields in the northeast Netherlands},
  school       = {University of Groningen},
  year         = {2015},
  url          = {https://hdl.handle.net/11370/ef1d6dd6-4b9d-4e4e-bc97-d257deb6724d},
  note         = {Chapter 2, p.~28}
}

@article{mikhaltsevitch_2014_measurements,
  title        = {
    {Measurements of the elastic and anelastic properties of sandstone flooded with
    supercritical CO$_2$}
  },
  author       = {Mikhaltsevitch, Vassily and Lebedev, Maxim and Gurevich, Boris},
  year         = 2014,
  month        = {09},
  journal      = {Geophysical Prospecting},
  volume       = 62,
  pages        = {1266--1277},
  doi          = {10.1111/1365-2478.12181}
}

@article{hu_2016_a,
  title        = {
    {A modified true triaxial apparatus for measuring mechanical properties of sandstone
    coupled with CO$_2$-H$_2$ biphase fluid}
  },
  author       = {
    Hu, Shaobin and Li, Xiaochun and Bai, Bing and Shi, Lu and Liu, Mingze and Wu,
    Haiqing
  },
  year         = 2016,
  month        = {09},
  journal      = {Greenhouse Gases: Science and Technology},
  volume       = 7,
  pages        = {78--91},
  doi          = {10.1002/ghg.1637}
}

@article{kim_2022_short,
  title        = {Short- and Long-Term Responses of Reservoir Rock Induced by {CO}$_2$ Injection},
  author       = {Kim, Kiseok and Makhnenko, Roman Y},
  year         = 2022,
  month        = {08},
  journal      = {Rock Mechanics and Rock Engineering},
  publisher    = {Springer Science+Business Media},
  volume       = 55,
  pages        = {6605--6625},
  doi          = {10.1007/s00603-022-03032-1}
}

@article{dnicolasespinoza_2018_co2,
  title        = {
    {CO$_2$ charged brines changed rock strength and stiffness at Crystal Geyser, Utah:
    Implications for leaking subsurface CO$_2$ storage reservoirs}
  },
  author       = {
    D. Nicolas Espinoza and Jung, Hojung and Major, Jonathan D and Sun, Zhuang and Ramos,
    Matthew J and Eichhubl, Peter and Balhoff, Matthew T and Choens, R C and Dewers,
    Thomas A
  },
  year         = 2018,
  month        = {06},
  journal      = {International Journal of Greenhouse Gas Control},
  publisher    = {Elsevier BV},
  volume       = 73,
  pages        = {16--28},
  doi          = {10.1016/j.ijggc.2018.03.017}
}

@article{fuchs_2019_geochemical,
  title        = {
    {Geochemical and geomechanical alteration of siliciclastic reservoir rock by
    supercritical CO$_2$-saturated brine formed during geological carbon sequestration}
  },
  author       = {
    Fuchs, Samantha J. and Espinoza, D. Nicholas and Lopano, Christina L. and Akono,
    Ange-Therese and Werth, Charles J.
  },
  year         = 2019,
  month        = {09},
  journal      = {International Journal of Greenhouse Gas Control},
  volume       = 88,
  pages        = {251--260},
  doi          = {10.1016/j.ijggc.2019.06.014}
}

@article{samuelson_2012_fault,
  title        = {
    {Fault friction and slip stability not affected by CO$_2$ storage: Evidence from
    short-term laboratory experiments on North Sea reservoir sandstones and caprocks}
  },
  author       = {Samuelson, Jon and Spiers, Christopher J.},
  year         = 2012,
  month        = 11,
  journal      = {International Journal of Greenhouse Gas Control},
  volume       = 11,
  pages        = {S78-S90},
  doi          = {10.1016/j.ijggc.2012.09.018}
}

@article{heinemann_2021_enabling,
  title        = {{Enabling large-scale hydrogen storage in porous media - the scientific challenges}},
  author       = {
    Heinemann, Niklas and Alcalde, Juan and M. Miocic, Johannes and T. Hangx, Suzanne J.
    and Kallmeyer, Jens and Ostertag-Henning, Christian and Hassanpouryouzband, Aliakbar
    and M. Thaysen, Eike and J. Strobel, Gion and Schmidt-Hattenberger, Cornelia and
    Edlmann, Katriona and Wilkinson, Mark and Bentham, Michelle and Haszeldine, R. Stuart
    and Carbonell, Ramon and Rudloff, Alexander
  },
  year         = 2021,
  journal      = {Energy \& Environmental Science},
  volume       = 14,
  pages        = {853--864},
  doi          = {10.1039/D0EE03536J}
}

@incollection{wassing_2017,
    author = {Wassing, B. B. T. and Buijze, L. and Ter Heege, J. H. and Orlic, B. and Osinga, S.},
    title = {The Impact of Viscoelastic Caprock on Fault Reactivation and Fault Rupture in Producing Gas Fields},
    booktitle = {U.S. Rock Mechanics/Geomechanics Symposium},
publisher = {Paper \#ARMA-2017-0355, American Rock Mechanics Association},
    year = {2017},
    month = {06}
}

@article{alyaseri_2023_experimental,
  author = {Al-Yaseri, A. and Amao, A. and Fatah, A.},
  month = {08},
  pages = {128272-128272},
  publisher = {Elsevier BV},
  title = {Experimental investigation of shale/hydrogen geochemical interactions},
  doi = {10.1016/j.fuel.2023.128272},
  volume = {346},
  year = {2023},
  journal = {Fuel}
}

@article{vasile_2024_innovative,
  author = {Vasile, Nicol\`o Santi and Bellini, Ruggero and Bassani, Ilaria and Vizzarro, Arianna and Azim, Annalisa Abdel  and Coti, Christian and Barbieri, Donatella and Scapolo, Matteo  and Viberti, Dario and Verga, Francesca and Pirri, Fabrizio and Menin, Barbara},
  month = {01},
  pages = {41-50},
  publisher = {Elsevier BV},
  title = {Innovative high pressure/high temperature, multi-sensing bioreactors system for microbial risk assessment in underground hydrogen storage},
  doi = {10.1016/j.ijhydene.2023.10.245},
  volume = {51},
  year = {2024},
  journal = {International Journal of Hydrogen Energy}
}

@article{shoushtari2023utilization,
  title     = {{Utilization of CO$_2$ and N$_2$ as cushion gas in underground gas storage process: A review}},
  author    = {Shoushtari, Sharif and Namdar, Hamed and Jafari, Arezou},
  journal   = {Journal of Energy Storage},
  volume    = {67},
  pages     = {107596},
  year      = {2023},
  publisher = {Elsevier},
  doi       = {10.1016/j.est.2023.107596}
}

@misc{PorthosCO2,
  title        = {{CO$_2$ reduction through storage under the North Sea}},
  author       = {{Porthos CO$_2$ Transport and Storage C.V.}},
  year         = 2024,
  note         = {Accessed: 2024-09-03},
  howpublished = {\url{https://www.porthosco2.nl/en}}
}

@misc{minEZK_grijpskerk2022,
  title        = {{Instemmingsbesluit Gasopslag Grijpskerk}},
  author       = {{Ministerie van Economische Zaken en Klimaat}},
  year         = {2022},
  month        = {July},
  url          = {https://www.rvo.nl/sites/default/files/2022/02/Instemmingsbesluit-Gasopslag-Grijpskerk.pdf},
  institution  = {Ministerie van Economische Zaken en Klimaat},
  language     = {Dutch}}

@techreport{tno_ondergrondse2018,
  title        = {{Ondergrondse Opslag in Nederland – Technische Verkenning}},
  author       = {{TNO}},
  year         = {2018},
  institution  = {TNO – Nederlandse Organisatie voor Toegepast Natuurwetenschappelijk Onderzoek},
  address      = {Utrecht, The Netherlands},
  language     = {Dutch}
}

@article{pijn19,
  title={Inelastic deformation of the Slochteren sandstone: Stress-strain relations and implications for induced seismicity in the {G}roningen gas field},
  author={Pijnenburg, RPJ and Verberne, BA and Hangx, SJT and Spiers, CJ},
  journal={Journal of Geophysical Research: Solid Earth},
  volume={124},
  number={5},
  pages={5254--5282},
  year={2019},
  publisher={Wiley Online Library}
}

@article{hol2018,
  title={Rock physical controls on production-induced compaction in the {G}roningen Field},
  author={Hol, Sander and van der Linden, Arjan and Bierman, Stijn and Marcelis, Fons and Makurat, Axel},
  journal={Scientific reports},
  volume={8},
  number={1},
  pages={7156},
  year={2018},
  publisher={Nature Publishing Group UK London}
}

@article{Zhao24,
    author = {Zhao, Xiaoxi and Jha, Birendra},
    title = {Role of Plasticity in Induced Seismicity Risk Mitigation: A Case of the {G}roningen Gas Field},
    journal = {SPE Journal},
    volume = {29},
    number = {12},
    pages = {7060-7073},
    year = {2024},
    month = {12},
    issn = {1086-055X},
    doi = {10.2118/223611-PA},
    url = {https://doi.org/10.2118/223611-PA},
    eprint = {https://onepetro.org/SJ/article-pdf/29/12/7060/4233366/spe-223611-pa.pdf},
}

@article{manj2023,
  title={{Role of CO2 in geomechanical alteration of Morrow Sandstone across micro--meso scales}},
  author={Manjunath, GL and Akono, Ange-Therese and Haljasmaa, I and Jha, Birendra},
  journal={International Journal of Rock Mechanics and Mining Sciences},
  volume={163},
  pages={105311},
  year={2023},
  publisher={Elsevier}
}

@inproceedings{prevost2013,
  title={One-way versus two-way coupling in reservoir-geomechanical models},
  author={Prevost, Jean H},
  booktitle={Poromechanics V: Proceedings of the Fifth Biot Conference on Poromechanics},
  pages={517--526},
  year={2013}
}

@article{van2015social,
  title={Social impacts of earthquakes caused by gas extraction in the Province of Groningen, The Netherlands},
  author={Van der Voort, Nick and Vanclay, Frank},
  journal={Environmental Impact Assessment Review},
  volume={50},
  pages={1--15},
  year={2015},
  publisher={Elsevier}
}

@article{pijnenburg2019inelastic,
  title={Inelastic deformation of the Slochteren sandstone: Stress-strain relations and implications for induced seismicity in the Groningen gas field},
  author={Pijnenburg, RPJ and Verberne, BA and Hangx, SJT and Spiers, CJ},
  journal={Journal of Geophysical Research: Solid Earth},
  volume={124},
  number={5},
  pages={5254--5282},
  year={2019},
  publisher={Wiley Online Library}
}

@book{geluk2005stratigraphy,
  title={Stratigraphy and tectonics of Permo-Triassic basins in the Netherlands and surrounding areas},
  author={Geluk, Marinus Cornelis},
  year={2005},
  publisher={Utrecht University}
}

@book{ziegler1990,
  title={Geological atlas of western and central Europe},
  author={Ziegler, Peter A and others},
  volume={52},
  year={1990},
  publisher={Shell Internationale Petroleum Maatschappij BV The Hague}
}

@misc{sodm_web,
  author = {{Staatstoezicht op de Mijnen (SodM)}},
  title = {About Staatstoezicht op de Mijnen (SodM)},
  url = {https://www.sodm.nl/over-ons/about-sodm-english},
  organization = {Sodm.nl}
}

%superscritp
%\Urlmuskip=0mu plus 1mu
%\bibliographystyle{elsarticle-num-names}
%\bibliography{biblio}

%% else use the following coding to input the bibitems directly in the
%% TeX file.

%%\begin{thebibliography}{00}
%%
%%%% \bibitem{label}
%%%% Text of bibliographic item
%%
%%\bibitem{}
%%
%%\end{thebibliography}
\end{document}